\documentclass[11pt]{amsart}
\usepackage{amsmath, amsfonts, amscd,  amssymb}
\usepackage[all]{xy}
\def\Si{{\cal D}^{{\rm sign}}}

\def\AAA{{\Bbb A}}
\def\BBB{{\Bbb B}}

\def \ub{{}^b}

\def  \nuint {\raise10pt\hbox{$\nu$}\kern-6pt\int}

\newcommand\Hom{\operatorname{Hom}}

\newcommand\Ch{\operatorname{Ch}}
\newcommand\STR{\operatorname{STR}}

\newcommand\ind{\operatorname{ind}}
\newcommand\Ind{\operatorname{Ind}}
\newcommand\ch{\operatorname{Ch}}
\newcommand\ha{\frac12}

\renewcommand\Im{\operatorname{Im}}

\def \P{\mathcal P}

\newcommand\Q{\mathcal Q}
\newcommand\R{\mathcal R}

\newcommand\B{\mathcal B}
\newcommand\Bi{\B^\infty}
\newcommand\D{\mathcal D}
\newcommand\Di{D\kern-6pt/}
\newcommand\cDi{{\mathcal D}\kern-6pt/}
\newcommand\CC{\mathbb C}
\def \n {\noindent}
\def \s {\smallskip}

\def \m {\medskip}
\def \cal {\mathcal}
\def \BB {\mathbb B}

\newcommand\NN{\mathbb N}
\newcommand\QQ{\mathbb Q}
\newcommand\RR{\mathbb R}
\newcommand\ZZ{\mathbb Z}

\def \Fl{{\cal V}}
\def \Fli{\Fl^\infty}

\newcommand\pa{\partial}

\newcommand\CI{{\mathcal C}^{\infty}}

\newcommand\Id{\operatorname{Id}}

\def \whW {\widehat{W}}

\newtheorem{theorem}{Theorem}[section]
\newtheorem{lemma}[theorem]{Lemma}
\newtheorem{proposition}[theorem]{Proposition}
\newtheorem{definition}[theorem]{Definition}

\newtheorem{corollary}[theorem]{Corollary}

\newtheorem{assumption}[theorem]{Assumption}
{\catcode`@=11\global\let\c@equation=\c@theorem}

\newcommand{\zz}{{\mathbb Z}}

\newcommand{\id}{\operatorname{id}}

\newcommand{\sign}{\operatorname{sign}}
\newcommand{\Ker}{\operatorname{Ker}}

\newcommand{\comment}[1]                      
{
{{\bf Comment: } {\ttfamily #1}}
}

\theoremstyle{remark}
\begin{document}

\title[ ]
{Elliptic operators and   higher signatures}
\author[E. Leichtnam]{Eric Leichtnam}
\address{Institut de Jussieu et CNRS,
Etage 7E,
175 rue du Chevaleret,
75013, Paris,  France}
        \email{leicht@math.jussieu.fr}
\author[P. Piazza]{Paolo Piazza}
        \address{Dipartimento di Matematica G. Castelnuovo,
Universit\`a di Roma ``La Sapienza'', P.le Aldo Moro 2, 00185 Rome, Italy}
       \email{piazza@mat.uniroma1.it}


\begin{abstract} $\;$\\
Building on the theory of elliptic operators, we give  a unified
treatment of the following topics: \begin{itemize} \item the
problem of homotopy invariance of Novikov's higher signatures on
closed manifolds;
\item the problem of  cut-and-paste invariance of Novikov's
higher signatures on  closed manifolds;
\item the problem of  defining higher signatures on manifolds with
boundary and proving their homotopy invariance.
\end{itemize}

\end{abstract}

\maketitle \tableofcontents

\section{{\bf Introduction.}}\label{sect:introduction}

 Let $M^{4k}$ an oriented $4k$-dimensional compact manifold.
Let $g$ be a Riemannian metric on $M$. Let us consider the
Levi-Civita connection $\nabla^g$ and the Hirzebruch $L$-form
$L(M,\nabla^g)$, a {\it closed} form  in $\Omega^{4*}(M)$ with de
Rham class $L(M):=[L(M,\nabla^g)]_{{\rm dR}}\in H^*(M,\RR)$
independent of $g$. Let now $M$ be closed; then 
\begin{equation}\label{intro-lower-closed}
\text{the integral over}\;\;M\;\;\text{of}\;\;L(M,\nabla^g)\;\;\text{is an oriented homotopy invariant
of } \;\; M\,.
\end{equation}
In fact, if $[M]\in H_* (M,\RR)$ denotes the fundamental class of
$M$ then $$\int_M L(M,\nabla^g)=<L(M),[M]>=\ind D^{{\rm
sign}}_{(M,g)}={\rm sign}(M)\,,$$ the last term denoting the
topological signature of $M$, an homotopy invariant of $M$. We
shall call the integral 
$\int_M L(M,\nabla^g)$
the {\it
lower signature of the closed manifold $M$}.

\m
 A second fundamental property of $\int_M
L(M,\nabla^g)\equiv <L(M),[M]>$ is its {\it cut-and-paste
invariance}: if $Y$ and $Z$ are two manifolds with diffeomorphic
boundaries and if $$X_\phi\,:=\,Y\cup_\phi Z^-\,,\quad
X_\psi\,:=\,Y\cup_\psi Z^-\,,\quad Z^-:= (-Z)$$ with $\phi, \psi\,:\pa Y\to \pa Z$
oriented diffeomorphisms, then
$<L(X_\phi),[X_\phi]>=<L(X_\psi),[X_\psi]>$.

\m
A third fundamental property will involve a manifold $M$ with
boundary. Using Stokes theorem we see easily that the integral
of the $L$-form is now metric dependent; in particular it is not
homotopy invariant. However, by the Atiyah-Patodi-Singer index
theorem for the signature operator, we know that there exists a
boundary correction term $\eta(\pa M,g_{|\pa M})$ such that
\begin{equation}\label{intro-lower-boundary}
\int_M L(M,\nabla^g)-\ha\eta(\pa M,g_{|\pa M})\quad\text{is an
oriented homotopy invariant}. \end{equation} In fact, this
difference equals the topological signature of the manifold with
boundary $M$. We call the difference appearing in
(\ref{intro-lower-boundary}) the {\it lower signature of the
manifold with boundary $M$}. The term $\eta(\pa M,g_{|\pa M})$,
i.e. the term  we need to subtract in order to produce a homotopy
invariant out of $\int_M L(M,\nabla^g)$, is a spectral invariant
of the signature operator $D^{{\rm sign}}_{(\pa M,g|_{\pa M})}$ on
$\pa M$; more precisely, this invariant measures the {\it
asymmetry} of the spectrum of this (self-adjoint) operator with
respect to $0\in \RR$.
 We
shall review these basic  facts in Section
\ref{sect:lower-closed.}
and Section \ref{sect:lower-cut}.  

\medskip

Let now $\Gamma$ be a finitely generated discrete group. Let
$B\Gamma$ be the classifying space for $\Gamma$. We shall be
interested in the real cohomology groups $H^*( B\Gamma,\RR)$.
 Let $\Gamma\to
\widetilde{M}\to M$ be a Galois $\Gamma$-covering of an oriented
manifold $M$. For example, $\Gamma=\pi_1 (M)$ and $\widetilde{M}$
is the universal covering of $M$. From the classifying theorem for
principal bundles we know that $\Gamma\to \widetilde{M}\to M$ is
classified by a continuous map $r:M\to B\Gamma$. We shall identify
 $\Gamma\to \widetilde{M}\to M$ with the pair $(M,r:M\to
 B\Gamma)$.
Assume at this point that $M$ is closed. Fix a class $[c]\in H^*
(B\Gamma,\RR)$; then $r^* [c]\in H^* (M,\RR)$ and it makes sense
to consider the number $<L(M)\cup r^*[c],[M]>\in\RR$. The
collection of real numbers $$\{\,{\rm
sign}(M,r;[c])\,:=\,<L(M)\cup r^* [c],[M]>\,,\;\;[c]\in
H^*(B\Gamma,\RR)\,\}$$  are called the {\it Novikov's higher
signatures} associated to the covering $(M,r:M\to B\Gamma)$. It is
important to notice that these number are {\it not} well defined
if $M$ has a boundary; in fact, in this case $L(M)\cup r^* [c]
\,\in\, H^* (M,\RR)$ whereas $[M]\in H_* (M,\pa M,\RR)\,,$ and the
two classes cannot be paired.\\

 One can give a
natural notion of {\it homotopy equivalence} between Galois
$\Gamma$-coverings. One can also give the notion of 2 coverings
being {\it cut-and-paste equivalent.} In this paper we shall
address the following three questions:

\bigskip
\n
{\bf Question 1.} Are Novikov's higher signatures
{\it homotopy invariant} ?

\smallskip
\n
{\bf Question 2.} Are Novikov's higher signatures {\it
cut-and-paste invariant} ?

\smallskip
\n
{\bf Question 3.} If $\pa M\not= \emptyset$, {\it can we define
higher signatures and prove their homotopy invariance ?} Of course
we want these higher signatures on a manifold with boundary $M$ to
generalize the lower signature $$\int_M L(M,\nabla^g) -{1\over 2}
\eta(\pa M,g_{|\pa})\,,$$ which is indeed a homotopy invariant.

\bigskip

 Question 1 is still open and is known
as the {\it Novikov conjecture}. It has been settled in the
affirmative for many classes of groups. In this survey we shall
present two methods for attacking the conjecture, both involving
in an essential way properties of {\it elliptic operators}.

The answer to Question 2 is negative: the higher signatures are
{\it not} cut-and-paste invariants (we shall present a
counterexample). However, one can give {\it sufficient conditions}
 on the group $\Gamma$ and on the separating hypersurface ensuring
that the higher signatures are indeed cut-and-paste invariant.

Finally, under suitable assumption on $(\pa M,r|_{\pa M})$ and on
the group $\Gamma$ one can {\it define} higher signatures on a
manifold with boundary $M$ equipped with a classifying map $r:M\to
B\Gamma$ and prove their homotopy invariance. 
Notice that part of the
problem in Question 3 is to give a meaningful definition.
Our answers to Question 2 and Question 3 will use in
a crucial way properties of elliptic boundary-value
problems.

\smallskip There are several excellent surveys on Novikov's
higher signatures; we mention here the very complete  historical
perspective by Ferry, Ranicki and Rosenberg \cite{FRR}, the
stimulating article  by Gromov \cite{Gromov-macroscopic} and 
the one by Kasparov \cite{Kasparov(1993)}. The {\it novelty} in
the present work is the {\it unified treatment} of closed
manifolds and manifolds with boundary as well as the {\it
treatment of the cut-and-paste problem} for higher signatures
on closed manifolds.

\smallskip 
\n
{\bf Acknowledgements.}
This article will appear in the proceedings of a
conference in honor of Louis Boutet de Monvel. The first author
was very happy to be invited to give a talk at this conference; he
feels that he learnt a lot of beautiful mathematics from Boutet de
Monvel, especially at E.N.S (Paris) during the eighties.

\section{{\bf The lower signature and its homotopy
invariance}}\label{sect:lower-closed.}

\subsection{The L-differential form.}\label{subsect:L-form}$\;$
\medskip

Let $(M,g)$ be an oriented Riemannian manifold of dimension $m$.
We fix a Riemannian connection $\nabla$ on the tangent bundle of
$M$ and we consider $ \nabla^2$, its curvature. In a fixed
trivializing neighborhood $U$ we have $\nabla^2=R$ with $R$ a
$m\times m$-matrix of 2-forms. We consider the $L$-differential
form $L(M,\nabla)\in \Omega^{*} (M)$ associated to $\nabla$.
Recall that $L(M,\nabla)$ is obtained by formally substituting the
matrix of 2-forms $\frac{\sqrt{-1}}{2 \pi} R$ in the power-series expansion
at $A=0$ of the analytic  function $$L(A)={\rm det}^{{1\over 2}}
\left({A\over \tanh A}\right)\,,\quad A\in so(m)\,.$$ Since
$\Omega^*(M)=0$ if $*>\dim M$, we see that the sum appearing in
$L( \frac{\sqrt{-1}}{2 \pi} R )$ is in fact finite. More importantly, since
$L(\cdot)$ is $SO(m)$-invariant, i.e. $$L(A)=L(C^{-1}AC)\,,\quad
C\in SO(m)\,,$$ one can check easily that $L(M,\nabla)$ is a
{\it globally defined}; it is a differential form in $\Omega^{4 *}(M, \RR).$ 
One can
prove the following two fundamental properties of the
$L$-differential form:
\begin{equation}\label{transgression}
d L(M,\nabla)=0,\;\;\;
L(M,\nabla)-L(M,\nabla^\prime)=d T(\nabla,\nabla^\prime)
\end{equation} where $\nabla^\prime$ is any other  Riemannian
connection. 
Consequently the de Rham class $L(M)=[L(M,\nabla)]\in H^* _{{\rm
dR}}(M)$ is well defined; it is called the Hirzebruch $L$-class.

\smallskip

In what follows we shall
always choose the Levi-Civita connection
associated to $g$, $\nabla^g$, as our reference connection.

\subsection{The lower signature on closed manifolds 
and its homotopy invariance.}\label{subsect:hominv-closed}$\;$

\medskip
Assume now that $M$ is {\it closed} ($\equiv$ {\it without
boundary}) and that $\dim M=4k$. Consider
\begin{equation}\label{integral-of-L} \int_M L(M,\nabla^g)
\end{equation}
Because of the  properties 
\eqref{transgression}, this integral does not depend on the choice
of $g$ and is in fact equal to $<L(M),[M]>$, the pairing between
the cohomology class $L(M)$ and the fundamental class $[M]\in
H_{4k}(M;\RR)$. 

\begin{theorem}\label{theo:hirzebruch}(Hirzebruch)
The integral of the L-form $$\int_M L(M,\nabla^g)$$ is an integer
and is an oriented homotopy invariant.
\end{theorem}

\begin{proof}
With some of what follows in mind, we give an  index-theoretic proof of
this theorem, in two steps.\\
{\it First step}:  by the Atiyah-Singer index theorem
$$\int_M L(M,\nabla^g)=\ind D^{{\rm sign},+}$$ where on the right
hand side the index of the signature operator associated to $g$
and our choice of orientation appears \footnote{Let us recall the
definition of the signature operator on a $2\ell$-dimensional
oriented Riemannian manifold. Consider the Hodge star operator
$$\star: \Omega^p (M)\longrightarrow \Omega^{2\ell-p} (M);$$ it
depends on $g$ and the fixed orientation. Let $\tau:=
(\sqrt{-1})^{p(p-1)+\ell} \star$ on $\Omega^{p}_{\CC}(M)$; then
$\tau^2=1$ and we have a decomposition $\Omega_\CC ^* (M)=\Omega^+
(M)\oplus \Omega ^- (M)$. The operator $d+d^*$, extended in the
obvious way to the complex differential forms $\Omega_\CC ^*(M)$,
anticommutes with $\tau$. The signature operator is simply defined
as $$D^{{\rm sign}}\,:=\,
\begin{pmatrix} 0 & {D^{{\rm sign},-}} \cr {D^{{\rm sign},+}}  & 0
\cr
\end{pmatrix}\,\;,\quad D^{{\rm sign},\pm}=(d+d^*)|_{\Omega^{\pm}
(M)}\,.$$ If we wish to be precise, we shall denote the signature
operator on the Riemannian manifold $(M,g)$ by $D^{{\rm
sign}}_{(M,g)}$. }. This proves that $$\int_M
L(M,\nabla^g)\,\in\,\ZZ\,.$$ {\it Second step}: using the Hodge theorem one can
check that $$\ind D^{{\rm sign},+}={\rm sign} (M):=\text{signature
of}\,M \, $$ i.e. the signature of the bilinear form
$H^{2k}(M)\times H^{2k} (M)\longrightarrow \RR$
$$([\alpha],[\beta])\longrightarrow \int_M \alpha\wedge \beta\,.$$
This is clearly an oriented homotopy invariant and the theorem
is proved.
\end{proof}

\s
\n
We shall also call $\int_M L(M,\nabla^g)$ the {\it lower
signature} of the closed manifold $M$.

\s
\noindent {\bf Remark.} The equality $${\rm
sign}(M)\,=<L(M),[M]>$$ is known as
the {\it Hirzebruch signature theorem}. The original proof of this
fundamental result was topological, exploiting the cobordism
invariance of both sides of the equation and the structure of the
oriented cobordism ring. See, for example, Milnor-Stasheff \cite{Milnor-Stasheff}
and  Hirzebruch \cite{Hirzebruch}.

\medskip
\n
{\bf Remark.}
The formulation and the proof of Hirzebruch theorem 
given here is not historically accurate but has the advantage
of introducing  the techniques  
that will be employed later for tackling the homotopy invariance of
higher signatures of a closed  manifold. It is important to single out
informally the two steps in the proof:

\n
(i) connect the lower signature to an index

\n 
(ii) prove that the index is  homotopy invariant.

\subsection{The lower signature on manifolds with boundary and its homotopy invariance
}\label{subsect:lower-boundary}

$\;$

\medskip
Assume now that $M$ has a {\it non-empty boundary}: $\partial
M=N\not= \emptyset$. For simplicity, we assume that the metric $g$
is of product-type near the boundary; thus in a collar
neighborhood $U$ of $\partial M$ we have $g=dx^2 + g_{\partial}$
with $x\in C^\infty (M)$ a boundary defining function. We denote
the signature operator on $(M,g)$ by $D^{{\rm sign}}_{(M,g)}$. We
consider once again
 $\int_M L(M,\nabla^g)$. In contrast with the closed case,
 this integral
 does depend now on the choice of the metric $g$;
 in
 particular it is {\it not} an oriented homotopy invariant.
To understand this point we simply observe that if $h$ is a
different metric, then, by (\ref{transgression}),
 we get
 \begin{equation}\label{transgression2}
 \int_M L(M,\nabla^g)-\int_M L(M,\nabla^h)=\int_{\partial M}
 T(\nabla^g,\nabla^h)|_{\partial M}\,.
 \end{equation}
We
ask ourselves if we can add to $\int_M L(M,\nabla^g)$ a correction
term making it metric-independent and, hopefully,
homotopy-invariant; formula (\ref{transgression2}) shows that it
should be possible to add a term that only depends on the metric
on $\partial M$. \\In order to state the result we need a few
definitions. Consider the boundary $\partial M$ with the induced
metric and orientation. Let $D^{{\rm sign}}_{\partial M,g_{\pa}}$
the signature operator on the odd dimensional Riemannian manifold
$(\partial M,g_{\pa})$; this is the so-called {\it odd signature
operator} and it is defined as follows:
\begin{equation}\label{odd-signature}
D^{{\rm sign}}_{(\partial M,g_{\pa})}\, \phi
\,:=\,(\sqrt{-1})^{2k}(-1)^{p+1}(\epsilon \star d - d \star)\phi
\end{equation} with $\epsilon=1$ if
$\phi\in\Omega^{2p}(\partial M)$ and $\epsilon=-1$ if
$\phi\in\Omega^{2p-1}(\partial M)$. This is a formally
self-adjoint first order elliptic differential operator on the
{\it closed} manifold $\partial M$. We shall sometime denote the
boundary signature operator  by  $D^{{\rm sign}}_{\partial}$.
Thanks to the spectral properties of elliptic differential
operators on closed manifolds, we know that the following series
is absolutely convergent for ${\rm Re} ( s )\gg 0$:
\begin{equation}\label{eta-s}
\eta(s) := \sum \lambda |\lambda|^{-(s+1)}\,,
\end{equation}
with $\lambda$ running over the non-zero eigenvalues of  $D^{{\rm
sign}}_{(\partial M,g_{\pa})}$. One can meromorphically continue
this function to the all complex plane; the points $s_k=\dim (\pa
M)-k$ are poles of the meromorphic continuation. It is a
non-trivial result that the point $s=0$ is regular and one sets
\begin{equation}
\eta(D^{{\rm sign}}_{(\partial M,g_{\pa})})\,:= \, \eta(0)
\end{equation}
This is the {\it eta invariant} associated to $D^{{\rm
sign}}_{(\partial M,g_{\pa})}$; it is a spectral invariant
measuring the {\it asymmetry} of the spectrum of $D^{{\rm
sign}}_{(\partial M,g_{\pa})}$, a subset of the real line, with
respect to the origin. The definition can be given for any
formally self-adjoint elliptic pseudodifferential operator. We can
now state the main theorem of this subsection:

\begin{theorem}\label{prop:hominv} (Atiyah-Patodi-Singer)
The difference
\begin{equation}\label{hominv}
\int_M L(M,\nabla^g)-{1\over 2} \eta(D^{{\rm sign}}_{(\pa
M,g_{\pa})})
\end{equation}
is an integer and is an oriented homotopy invariant of the pair
$(M,\partial M)$.\\
We call the difference $\int_M L(M,\nabla^g)-{1\over 2} \eta(D^{{\rm sign}}_{(\pa
M,g_{\pa})})$ the {\it lower
signature} of the manifold with boundary $M$.

\end{theorem}

\begin{proof}
Following Atiyah-Patodi-Singer \cite{APS1} we give an index theoretic proof
of this theorem, once again in two steps.\\
 There is a
well defined restriction map $$|_{\partial M}: \Omega^+
(M)\rightarrow \Omega^* (\pa M)\,.$$ Next we observe that to the
formally self-adjoint operator $D^{{\rm sign}}_{(\pa M,g_{\pa})}$
we can associate the spectral projection $\Pi_{\geq}$ onto the eigenspaces associated to its
nonnegative eigenvalues. The operator  $D^{{\rm
sign},+}_{(M,g)}\equiv D^{{\rm sign},+}$ on $M$ with boundary
condition $$\omega^+ |_{\partial M} \in \Ker \Pi_{\geq}$$ turns
out to be Fredholm when acting on suitable Sobolev completions
(more on this in the next subsection). The Atiyah-Patodi-Singer
($\equiv$ APS) index formula computes its index as
\begin{equation}\label{apsformula} \ind (D^{{\rm
sign},+},\Pi_{\geq})=\int_M L(M,\nabla^g)-{1\over 2} \eta(D^{{\rm
sign}}_{(\pa M,g_{\pa})})-{1\over 2} \dim \Ker (D^{{\rm
sign}}_{(\pa M,g_{\pa})})\,. \end{equation}
 It should be remarked that
$\Ker (D^{{\rm sign}}_{(\pa M,g_{\pa})})$ has a natural symplectic
structure and it is therefore even dimensional. From
(\ref{apsformula}) we infer that \begin{equation}\label{apsformulabis}
 \int_M L(M,\nabla^g)-{1\over
2} \eta(D^{{\rm sign}}_{(\pa M,g_{\pa})})= \ind (D^{{\rm
sign},+},\Pi_{\geq})+{1\over 2} \dim \Ker (D^{{\rm sign}}_{(\pa
M,g_{\pa})})\,.\end{equation}
This concludes the first step, connecting the lower signature to an index.\footnote{It could be proved that
the right hand side of (\ref{apsformulabis}) is the index of the boundary value problem
corresponding to the projection $\Pi_{>}+\Pi_L$, where $L\subset  \Ker (D^{{\rm sign}}_{(\pa
M,g_{\pa})})$ is the so-called scattering lagrangian.}
Next, using Hodge theory on the
complete manifold $\widehat{M}$ obtained by gluing to $M$ a
semi-infinite cylinder $(-\infty,0]\times \pa M$, one can prove
that $$\ind (D^{{\rm sign},+},\Pi_{\geq})+{1\over 2} \dim \Ker
(D^{{\rm sign}}_{(\pa M,g_{\pa})})={\rm sign} (M):=\text{the
signature of}\, M\,.$$ Since the latter is an oriented homotopy
invariant, the theorem  is proved.
\end{proof}

\n
{\bf Remark.} In contrast with the closed case, there is no purely
topological proof of the homotopy invariance of the difference
$\int_M L(M,\nabla^g)-{1\over 2} \eta(D^{{\rm sign}}_{(\pa
M,g_{\pa})})$ on manifolds with boundary; in this case {\it we do
need to pass through the  Atiyah-Patodi-Singer index theorem}. 


\smallskip
\n
{\bf Remark.} 
Part of the motivation for the Atiyah-Patodi-Singer
 signature index theorem came
from the work of Hirzebruch on Hilbert modular
varieties. For these singular varieties the Hirzebruch signature
formula does not hold; there is a {\it defect} associated
to each cusp.
For Hilbert modular surfaces Hirzebruch computed this defect
and showed that it
was given in terms of the value at $s=1$ of certain $L$-series. 
He then conjectured that a similar result was true for any
Hilbert modular variety. The conjecture was
established by Atiyah-Donnelly-Singer in 
\cite{ADS} \cite{ADSbis} and the proof is based in an essential
way on the Atiyah-Patodi-Singer
index theorem (with the value of the $L$-series corresponding
to the eta-invariant). Hirzebruch's conjecture was also settled
independently and with a different proof by M\"uller
in \cite{MuJDG} (see also \cite{Mucontem}).

\smallskip
\n
{\bf Remark.} It is possible to prove that for the {\it odd}
signature operator $D^{{\rm sign}}_N$ on an odd dimensional
manifold oriented closed manifold $N$ $$\eta(D^{{\rm
sign}}_N)={2\over\sqrt{\pi}}\int\limits_0^\infty {\rm Tr} (D^{{\rm
sign}}_N \,e^{-(t\,D^{{\rm sign}}_N)^2}) dt\,.$$ Notice that the
convergence of this integral near $t=0$ is non-trivial and its
justification requires arguments similar to those involved in the
heat-kernel proof of the Atiyah-Singer index theorem, see
Bismut-Freed \cite{BF1}, \cite{BF2}.

\subsection{More on index theory on manifolds with boundary.}
\label{subsect:moreAPS}$\;$

\medskip
We elaborate further on the analytic features of the above proof.
 Let $M$ be a manifold with boundary. Simple examples (such
as the $\overline{\pa}$-operator on the disc) show that, in
general, elliptic  operators on $M$ are {\it not} Fredholm  on
Sobolev spaces. In order to obtain a finite dimensional kernel and
cokernel it is necessary to impose {\it boundary conditions}.
Among the simplest boundary conditions are those of {\it local
type}, Dirichlet, Neumann or more generally Lopatinski boundary
conditions. It is not at all clear that these classical local
boundary conditions give rise to Fredholm operators. And in fact
Atiyah and Bott  showed that there exist topological obstructions
to the existence of well-posed  {\it local} boundary conditions
for an elliptic operator on a manifold with boundary. When these
obstructions are zero, Atiyah and Bott do prove an index theorem,
see \cite{A-Bott}. The Atiyah-Bott index theorem has been greatly
extended by Boutet de Monvel in \cite{Boutet}. However, precisely
 because of their geometric nature,
the signature operator is among those operators for which these
obstructions are almost always {\it non-zero}. In trying to prove
the  signature theorem on manifolds with boundary, Atiyah, Patodi
and Singer introduced their celebrated  {\it non-local} boundary
condition. This is the boundary condition explained in the proof
of Theorem \ref{prop:hominv}. In a fundamental series of
papers \cite{APS1} \cite{APS2} \cite{APS3}  they investigated the
index theory of such boundary value problems for general
first-order elliptic differential operators; they also gave
important applications to geometry and topology. Their theory
applies to any Dirac-type operator on an even dimensional manifold
with boundary endowed with a Riemannian metric $g$ which is of
product-type near the boundary. The Dirac operators acts between
the sections of a $\ZZ_2$-graded Hermitian Clifford module
$E=E^+\oplus E^-$ endowed with a Clifford connection $\nabla^E$
and it is odd with respect to the grading of $E$:
 $$D=\begin{pmatrix} 0 & D^- \cr D^+ & 0
\end{pmatrix}$$
 Classical {\it examples} of Dirac-type operators are given by the signature
operator $D^{{\rm sign}}$ introduced above, the Gauss-Bonnet
operator $d+d^*$, with $d$ equal to the de Rham differential, {\it
the} Dirac operator $\Di$ on a spin manifold, the
$\overline{\pa}$-operator on a Kaehler manifold. See Berline-Getzler-Vergne
\cite{BGV}
for more on Dirac operators.

 Near the boundary $D$ can be written (up to a bundle
isomorphism)
 as $$\begin{pmatrix} 0 & -\pa/\pa u+D_{\pa M}\cr \pa/\pa
u+D_{\pa M} & 0
\end{pmatrix}$$ with $u$ equal to the inward normal variable to
the boundary and $D_{\pa M}$ the generalized Dirac operator
induced on $\pa M$. For example, in the case of the signature
operator $D^{{\rm sign}}$ the operator induced on the boundary is
simply the odd-signature operator. The boundary operator $D_{\pa
M}$ is an elliptic and essentially self-adjoint operator on the
{\it closed} compact manifold $\pa M$. The $L^2$-spectrum is
therefore discrete and real. Let $\{e_\lambda \}$ be an
$L^2$-orthonormal basis of eigenfunctions for $D_{\pa M}$. Let
$\Pi_{\geq}$ be the spectral projection corresponding to  the
non-negative eigenvalues of $D_{\pa M}$: thus $\Pi_{\geq}
(e_\lambda)=e_\lambda$ if $\lambda\geq 0$ and $\Pi_{\geq}
(e_\lambda)=0$ if $\lambda<0$. Let $$\CI(M,E^+,\Pi_{\geq})=
\{\,s\in\CI(M,E^+)\;|\;\Pi_{\geq}(s|_{\pa M})=0\,\}.$$ Thus a
section $s$ belongs to $\CI(M,E^+,\Pi_{\geq})$ iff $s|_{\pa
M}=\sum_{\lambda <0} s^\lambda_{\pa M}e_\lambda$. The
Atiyah-Patodi-Singer theorem, see \cite{APS1}, states that the
operator $D^+$  acting on the Sobolev completion $H^1
(M,E^+,\Pi_{\geq})$ of $\CI(M,E^+,\Pi_{\geq})$, with range $L^2(
M,E^-)$, is a Fredholm operator with index $$\ind(D^+,
\Pi_{\geq})=\int_M\,{\rm AS}\,-{1\over 2}(\eta( D_{\pa M})
+\dim\Ker D_{\pa M}).$$ Here $\eta( D_{\pa M})$ is the {\it eta
invariant} of the self-adjoint operator $D_{\pa M}$ as introduced
in the previous subsection, whereas the density ${\rm AS}=
\widehat{A}(M,\nabla^g) {\rm ch}^\prime (E,\nabla^E)$ is the local
contribution that would appear in the heat-kernel proof of the
Atiyah-Singer index theorem for Dirac operators.
 In the case $D$ is  Dirac operator acting on the spinor bundle
 of a spin manifold,
one has  ${\rm
AS}= 
\widehat{A}(M,\nabla^g) $. In the case where $D$ is the signature operator
acting on the  bundle of differential forms one has  ${\rm
AS}= 
L(M,\nabla^g) $. 
The $\widehat{A}$-form 
$\widehat{A}(M,\nabla^g)$ is obtained by
substituting  $X$ by $ \frac{\sqrt{-1}}{ 2 \pi} R$ in the analytic
functions $$ \widehat{A} (X) = {\rm det}^{{ 1\over 2}} \left( {X/2
\over \sinh X/2} \right). 
$$ There are nowadays many
alternative approaches to the Atiyah-Patodi-Singer index formula;
we shall mention here the one started by Cheeger, based on conic
metrics (see Cheeger \cite{Cheeger}, Chou \cite{Chou} and also 
Lesch \cite{Lesch})
and the one, fully developed by Melrose, based on manifolds with
cylindrical ends (see Melrose \cite{Melrose} and also Piazza \cite{piazzajfa},
Melrose-Nistor \cite{Melrose-Nistor}). For a proof in the spirit of the embedding
proof of the Atiyah-Singer index formula on closed manifolds see
Dai-Zhang \cite{dai-zhang-embedding}.

\smallskip
\n
{\bf Remark.} Let $P=P^2=P^*$ be a finite rank perturbation of the
projection $\Pi_{\geq}$. Thus, with $\{e_\lambda \}$ still
denoting a $L^2$-orthonormal basis of eigenfunctions for $D_{\pa
M}$, we require that for some $R>0$, $P e_\lambda=e_\lambda$ if
$\lambda > R$ and $P e_\lambda=0$ if $\lambda<-R$. Let
$$\CI(M,E^+,P)= \{\,s\in\CI(M,E^+)\;|\;P(s|_{\pa M})=0\,\}.$$ The
operator $D^+$ with domain $\CI (M,E^+,P)$ extends once again to a
Fredholm operator with $\ind (D^+,P)\in \ZZ$. See, for example
Booss-Bavnbek - Wojciechowski \cite{booss-woj}. 
Moreover: let $P_1$ and $P_2$ be two such
projections and let us consider $H_j= P_j (L^2(\pa M,E|_{\pa
M}))$. One can show easily that the operator $P_2\circ P_1: H_1\to
H_2$ is Fredholm; its index is called the {\it relative index} of
the two projections and  is denoted by $i(P_1,P_2)$.  The
following formula is known as the {\it relative index formula}
(\cite{booss-woj}):  $$\ind (D^+,P_2)-\ind
(D^+,P_1)=i(P_1,P_2)\,.$$ For example: $ \ind(D^+,
\Pi_{>})-\ind(D^+, \Pi_{\geq})=i(\Pi_{\geq},\Pi_>)=\dim \Ker
D_{\pa M}$

\section{{\bf The cut-and-paste invariance of the lower
signature.}}\label{sect:lower-cut}

Let $M$ and $N$ be two compact $4k-$dimensional oriented manifolds
with boundary and let $\phi, \, \psi:
\partial M \rightarrow \partial N$ be two orientation preserving diffeomorphims.
Let $N^-$ be $N$ with the reverse orientation. By gluing we obtain
two closed oriented $4k-$dimensional manifolds, $ M\cup_\phi N^-$
and $ M\cup_\psi N^-$. We shall say that  $$ M\cup_\phi
N^-\quad\text{and}\quad M\cup_\psi N^-\quad\text{are {\it
cut-and-paste equivalent}}\,.$$
Consider the two integers
$<L(M\cup_\phi N^- ), [M\cup_\phi N^-]>$ and $<L(M\cup_\psi N^- ), [
M\cup_\psi N^- ]>\,.$

\begin{proposition}\label{prop:cut-and-paste}
The following equality holds:
\begin{equation}\label{cut-and-paste}
<L(M\cup_\phi N^- ), [M\cup_\phi N^-]>\,=\,<L(M\cup_\psi N^- ), [
M\cup_\psi N^- ]>\,.
\end{equation}
In words, the integral of the $L$-class is a cut-and-paste
invariant.
\end{proposition}

\smallskip
In the next three subsections we shall give three different proofs
of this proposition.

\subsection{The index-theoretic proof.}\label{subsect:index-proof}$\;$

\medskip
We set $$X_\phi:= M\cup_\phi N^- \quad \text{and} \quad X_\psi:=
M\cup_\psi N^-\,.$$ Using the Atiyah-Patodi-Singer index theorem
we shall prove that
\begin{equation}\label{additivity-formula}
<L(X_\phi),[X_\phi]>={\rm sign}(M)-{\rm sign} (N)=
<L(X_\psi),[X_\psi]>\,.\end{equation} Notice that the 2 manifolds
$X_\phi$ and $X_\psi$ are, in general, {\it distinct}. Fix metrics
$g_\phi$ and $g_\psi$ on $ X_\phi $ and $ X_\psi $ respectively.
Since the integral of the $L$-class on closed manifolds in metric-independent, 
we can assume that these metrics are of product
type near the embedded hypersurface $F:=\partial M$. Thus we can
write $$X_\phi= M\cup_{\Id} {\rm Cyl}_\phi \cup_{\Id} N^-$$ with
$$ {\rm Cyl}_\phi := ([-1,0]\times (\partial M)^-) \cup _\phi
([0,1]\times
\partial N)\,.$$
Denoting generically by $\nabla^{{\rm LC}}$ the Levi-Civita
connection associated to the various restrictions of $g_\phi$, we
can write
\begin{align*}
\int_{X_\phi} L(X_\phi,\nabla^{{\rm LC}})\,&=\, \int_M
L(M,\nabla^{{\rm LC}})+\int_{{\rm Cyl}_\phi} L({\rm
Cyl}_\phi,\nabla^{{\rm LC}} )- \int_{N} L(N,\nabla^{{\rm LC}})\\
&= \,\int_M L(M,\nabla^{{\rm LC}})-{1\over 2} \eta(D^{{\rm
sign}}_{(\pa M,g_{\pa})})
\\ & \, +{1\over 2} \eta(D^{{\rm
sign}}_{(\pa M,g_{\pa})}) + \int_{{\rm Cyl}_\phi} L({\rm
Cyl}_\phi,\nabla^{{\rm LC}} ) - {1\over 2} \eta(D^{{\rm
sign}}_{\partial N})
\\ & \, - (\int_{N} L(N,\nabla^{{\rm LC}})  - {1\over 2} \eta(D^{{\rm
sign}}_{\partial N})) \\&= \, {\rm sign}(M) + {\rm sign} ({\rm
Cyl}_\phi) - {\rm sign} (N) \\&= {\rm sign}(M)-{\rm sign}(N)\,.
\end{align*}
We explain why these equalities hold. The first one is obvious; in
the second one, we simply added and substracted the same
quantities; in the third one, we applied the Atiyah-Patodi-Singer
theorem, keeping in mind that the eta invariant is orientation
reversing; in the fourth one, we used the topological invariance
of ${\rm sign} ( \cdot)$ together with the following two
observations:

\smallskip
\noindent (i) the diffeomorphism $\phi$ induces a diffeomeorphism
 between ${\rm Cyl}_\phi$ and $[-1,1]\times \partial M$;

 \noindent (ii) ${\rm sign} ([-1,1]\times \partial M)=0$ (use again
 the APS-formula).

\smallskip
\noindent Since exactly the same argument can be applied to
$X_\psi$, it follows that we have proved
(\ref{additivity-formula}) and thus Proposition
\ref{prop:cut-and-paste}.

\subsection{The topological proof.}\label{subsect:topology-proof}$\;$

\medskip
We start with a simplified situation. Let $X=M\cup N^-$ with $\pa
M=\pa N$; in other words $\phi=\Id$. Using Poincar\'e duality and
reasoning in terms of intersection of cycles one can prove in a
purely topological way the following {\it Novikov gluing formula}
(Hirsch \cite{Hirsch}):
\begin{equation}\label{novikov-add}
{\rm sign}(M\cup N^-)={\rm sign}(M)-{\rm sign}(N)\,.
\end{equation}
Then, using exactly the same reasoning as in the previous section,
one shows that for two different diffeomorphisms $\phi$ and $\psi$
 $${\rm sign}\, (M\cup_\phi N^- ) \,=\, {\rm sign}\, (M)- {\rm sign}\, N \,=\,  {\rm sign}\,
(M\cup_\psi N^- ).$$ By the Hirzebruch signature formula this
implies  $$<L(M\cup_\phi N^-),[M\cup_\phi N^-]>= <L(M\cup_\psi
N^-),[M\cup_\psi N^-]>$$ which is the formula we wanted to prove.

\smallskip
 Following a
suggestion of W. L\"uck,
 we shall now give a more  algebraic proof of this equality. This
 should be considered as a {\it pr\'elude} to the
  arguments of Leichtnam-L\"uck-Kreck \cite{LLK}
that we shall recall in Section \ref{sect:topological-cut} below.
 Since every sub vector space of a real vector space is a direct
summand, one can construct a chain homotopy equivalence $u$
between the cellular chain complex of $\RR-$vector spaces
$C_*(\partial M)$ and a chain complex $D_*$ of finite dimensional
$\RR-$vector spaces  whose $m-$differential  $d_m: D_m \rightarrow
D_{m-1}$ vanishes. With these notations, set
 $\overline{D}_i= D_i$  for
$0\leq i \leq m-1$ and $\overline{D}_i=0$ for $i\geq m$. One then gets a so-called Poincar\'e pair
 $j_*:  D_* \rightarrow \overline{D}_*$ whose boundary is $ {D}_*$.
  By glueing
$j_*:  D_* \rightarrow \overline{D}_*$ and  the Poincar\'e pair
$i_*: C_*(\partial M) \rightarrow C_*(M)$ along their boundaries
with the help of $u$ one gets
a true algebraic Poincar\'e complex denoted $C_*( M \cup_u  \overline{D}) .
$
A reference for these concepts is Ranicki \cite[page 18]{Ranicki(1981)}.
Intuitively
 an algebraic Poincar\' e pair
 $j_*: D_* \rightarrow \overline{D}_*$ is the algebraic analogue
 of the injection $i: \partial M \rightarrow M$ where
 $M$ is an oriented manifold with boundary.
 One can check that the signature ${\rm sign}\, (M \cup_u  \overline{D} )$
 of the  non degenerate quadratic form of $C_*( M \cup_u  \overline{D})$ does
not depend on the choice of $u$ and $D$. Moreover one can prove that
  the signature ${\rm sign}\,( \overline{D}_*,  D_*) $  of the
algebraic Poincare pair $j_*:  D_* \rightarrow \overline{D}_*$ is zero.
\begin{lemma} One has:
$$
{\rm sign}\, (M\cup_\phi N^- ) ={\rm sign}\, M\,-\, {\rm sign}\, N\,=
{\rm sign}\, (M\cup_\psi N^- ).
$$
\end{lemma}
\begin{proof} Of course the second equality is a consequence of the first one.
The  algebraic Poincar\'e complex defined by the cellular chain complex
 $ C_*( M\cup_\phi N^-)$ is (algebraically) cobordant to
the following algebraic Poincar\'e complex:
$$
C_*( M \cup_u  \overline{D}) +
C_*( N^- \cup_{u\circ \phi^{-1}}  \overline{D}^-).
$$ Hence the signature of $M\cup_\phi N^-$ is the sum of the ones of
$C_*( M \cup_u  \overline{D}) $ and $C_*( N^- \cup_{u\circ \phi^{-1}}  \overline{D}^-).$ But one has:
$$
{\rm sign}\,( C_*( M \cup_u  \overline{D}) )=
{\rm sign}\,M + {\rm sign}\,( \overline{D}_* , D_*),
$$
$$
{\rm sign}\,(C_*( N^- \cup_{u\circ \phi^{-1}}  \overline{D}^-))=
{\rm sign}\,N^- + {\rm sign}\,( \overline{D}^-_* , D^-_*).
$$
Since
 the signature of $( \overline{D}_* , D_*)$  is zero
one gets that sign$ M\cup_\phi N^-$ = sign$ M$ - sign$ N$ which
proves the Lemma.
\end{proof}

\subsection{The spectral-flow-proof. }\label{subsect:flow-proof}$\;$

\medskip

Recall that we have set  $$X_\phi:= M\cup_\phi N^- \quad
\text{and} \quad X_\psi:= M\cup_\psi N^-\,.$$ Fix metrics $g_\phi$
on $X_\phi$ and $g_\psi$ on $X_\psi$. We shall assume that these
metrics are of product type near the embedded hypersurface
$F:=\partial M$.
 We shall prove,  {\it analytically and without making use
 of the Atiyah-Patodi-Singer index formula},
 that \begin{equation}\label{equality-indeces}
 \ind (D^{{\rm
sign},+}_{X_\phi}) = \ind (D^{{\rm sign},+}_{X_\psi})\,.
\end{equation} By the Atiyah-Singer index theorem for the signature
operator on closed manifolds, this will suffice in order to
establish $<L(X_\phi),[X_\phi]>=<L(X_\psi),[X_\psi]>$, i.e.
Proposition \ref{prop:cut-and-paste}. The equality of the two
indeces will be obtained exploiting two fundamental properties of
the Atiyah-Patodi-Singer index: the {\it variational formula} and
the {\it gluing formula}.

\subsubsection{{\bf The variational formula for the
APS-index.}}\label{subsubsect:variational}

In contrast with the closed case, the APS-index is not stable
under perturbations. In Subsection \ref{subsect:moreAPS}  we have
defined the APS-boundary value problem for any generalized Dirac
operator on an even dimensional  manifold with boundary, $M$,
endowed with a
 metric $g$ which is of product type near the boundary.
 Assume now that $\{D(t)\}_{t\in [0,1]}$
is a smoothly varying family of such operators. As an important
example we could consider a family of metrics $\{g(t)\}_{t\in
[0,1]}$ on $M$ and the associated family of signature operators
$\{D^{{\rm sign}} (t)\}_{t\in [0,1]}$. Going back to the general
case, consider the family of operators induced on the boundary
$\{D_{\pa M} (t)\}_{t\in [0,1]}\,$; let $\Pi_{\geq} (t)$ the
corresponding spectral projection associated to the non-negative
eigenvalues;
then the following
variational formula for the APS-indeces holds:
\begin{equation}\label{variational}
\ind (D^+ (1),\Pi_{\geq} (1))-\ind (D^+ (0),\Pi_{\geq} (0))={\rm
sf} (\{D_{\pa M} (t) \}_{t\in [0,1]})
\end{equation}
where on the right hand side the {\it spectral flow} of the
1-parameter family of self-adjoint operators $\{D_{\pa M} (t) \}$
appears; this is the net number of eigenvalues changing sign as
$t$ varies from $0$ to $1$ (\cite{APS3}, \cite{Melrose}).
 Formula (\ref{variational}) follows
from the APS-index formula, see \cite{APS3}. It can also be proved
analytically, without making use of the APS-index formula. See for
example Dai-Zhang \cite{Dai-Zhang}.

\subsubsection{{\bf Important remark.}}\label{subsubsect:important}
 If $N$ is odd dimensional and
$\{D^{{\rm sign}}_N (t)\}_{t\in [0,1]}$ is a one-parameter family
of odd signature operators parametrized by a path of metrics $g_N
(t)_{t\in [0,1]}$, then
\begin{equation}\label{zero-spectral-flow}
{\rm sf}(\{D^{{\rm sign}}_N (t)\}_{t\in [0,1]})=0
\end{equation}
In fact, the kernel of the odd signature operator is equal to the
space of  harmonic forms on $N$; from the Hodge theorem we know
that such a vector space is independent of the metric we choose;
thus there are not eigenvalues changing sign and the spectral flow
is zero.

\subsubsection{{\bf The gluing
formula.}}\label{subsubsect:additive}

We start with a simplified situation: $X$ is a {\it closed}
compact manifold which is the union of two manifolds with
boundary. Thus there exists an embedded hypersurface $F$ which
separates $M$ into two connected components and such that
$$X=M_+\cup_F M_-\,,\quad\text{with}\quad \pa M_+=\pa M_-=F\,.$$
We assume that the metric $g$ is of product type near the
hypersurface $F$, i.e. near the boundaries of $M_+$ and $M_-$. Let
$D_X$ be a Dirac-type operator on $X$; then we obtain in a natural
way two Dirac operators on $M_+$ and $M_-$. The following gluing
formula holds:
\begin{equation}\label{additivity-simple}
\ind (D_X)=\ind (D_{M_+},\Pi_{\geq})+ \ind (D_{M_-},1-\Pi_{\geq})
\,. \end{equation} The discrepancy in the spectral projections
come from the orientation of the normals to the two boundaries (if
one is inward pointing, then the other is outward pointing).
\footnote{Notice that $ 1-\Pi_{\geq}$ is not exactly the
APS-projection associated to the non-negative eigenvalues of
$D_{\pa M_-}$; to be precise  $ 1-\Pi_{\geq}=\Pi^{\pa M_-}_{>},$
the projection onto the {\it positive} eigenvalues of $D_{\pa
M_-}$}

\smallskip

Formula (\ref{additivity-simple}) can be proved directly, in a
purely analytical fashion, see Bunke \cite{Bunke}, Leichtnam-Piazza
\cite{LPCUT}. Of
course it is also a consequence of the APS-index theorem.

\smallskip

\subsubsection{{\bf Proof of formula (\ref{equality-indeces}).}}\label{subsubsect:proof}
The gluing formula (\ref{additivity-simple}) can be generalized to
our more complicated situation, where $X_\phi$ is a closed
manifold obtained by {\it gluing} two manifolds with boundary {\it
through a diffeomorphism}. Using  this gluing formula on $X_\phi$
(with metric $g_\phi$) and on $X_\psi$ (with metric $g_\psi$),
applying then  the variational formula for the APS index on $M$
with respect to a path a metrics connecting $g_\psi |_M$ to
$g_\phi |_M$ and then doing the same on $N$ (with a path of
metrics connecting $g_\psi| _N$ and $g_\phi| _N$), one proves that
$ \ind (D^{{\rm sign},+}_{X_\psi})- \ind (D^{{\rm
sign}}_{X_\phi})= {\rm sf} (\{D^{\rm sign}_{{\rm odd}}
(\theta)\}_{\theta\in S^1}$. The spectral flow appearing in this
formula is associated to a $S^1$-family of odd signature operators
acting on the fibers of the mapping torus $F\to M(\phi^{-1}
\psi)\to S^1$ and parametrized by a family of metrics. As remarked
in \ref{subsubsect:important} this spectral flow is zero because
of the cohomological significance of the zero eigenvalue for the
signature operator. References for this material are, for example,
the book \cite{booss-woj} and the survey Mazzeo-Piazza \cite{piazza2}. {\it
Summarizing}, the equality of $\ind(D^{{\rm
sign}}_{X_\phi})=\ind(D^{{\rm sign}}_{X_\phi})$ has been obtained
through the following two equalities
\begin{equation}\label{zero-with-mapping-torus}
\ind (D^{{\rm sign},+}_{X_\psi})- \ind (D^{{\rm sign}}_{X_\phi})=
{\rm sf} (\{D^{\rm sign}_{{\rm odd}} (\theta)\}_{\theta\in S^1})
=0\,.\end{equation}

\n
{\bf Remark.} It should be remarked that in this third proof we
have not used the APS-index formula; only the analytic properties
of the APS boundary value problem were employed. This will be
important later, when we shall consider higher signatures.

\section{{\bf Summary.}}\label{sect:summary1}

Let us summarize what we have seen so far. Let $(M,g)$ be an
oriented Riemannian manifold of dimension $4k$ and let $D^{{\rm
sign}}_{(M,g)}$ be the associated signature operator.

\begin{itemize}
\item {\it If $M$ is closed then $\int_M L(M,\nabla^g)$ is an oriented
homotopy invariant}. In fact $$\int_M L(M,\nabla^g)=
<L(M),[M]>=\ind D^{{\rm sign},+}_{(M,g)}= {\rm sign} (M)\,,$$ with
$L(M)=[L(M,\nabla^g)]\in H_{{\rm dR}}^*(M)\,,\,[M]\in H_* (M,\RR)
$ and ${\rm sign}(M)=$ signature of $M$.
\item If $M$ has a boundary, $\pa M\not= \emptyset$, then we can
define a correction term $\eta(D^{{\rm sign}}_{(\pa M,g_\pa)})$
such that $$\int_M L(M,\nabla^g)-{1\over 2}\eta(D^{{\rm
sign}}_{(\pa M,g_\pa)})$$ is an oriented homotopy invariant of the
pair $(M,\pa M)$. In fact
 $$\int_M L(M,\nabla^g)-{1\over 2}\eta(D^{{\rm
sign}}_{(\pa M,g_\pa)})= \ind (D^{{\rm sign}}_{(M,g)},\Pi_{\geq})+
{1\over 2}\dim\Ker(D^{{\rm sign}}_{(\pa M,g_\pa)})={\rm sign}
(M)\,.$$
\item
Let $X_\phi=M\cup_\phi N^-$ and $X_\psi=M\cup_\psi N^-$ be two
{\it cut-and-paste equivalent} closed  manifolds. Then
$$<L(X_\phi),[X_\phi]>\,=\,<L(X_\psi),[X_\psi]>\,.$$
\end{itemize}

\section{{\bf Novikov higher
signatures.}}\label{sect:higher-closed}

\subsection{Galois coverings and classifying maps.}$\;$

\m

Let $\Gamma$ be a discrete finitely presented group.
 Let $\Gamma\longrightarrow \widetilde{M} \longrightarrow M$
be a Galois $\Gamma$-covering (the term {\it normal} is also in
common usage). For example $\Gamma:=\pi_1(M)$ and $\widetilde{M}$=
universal covering of $M$. As a particular example to keep in
mind, let $\Sigma_g$ be a closed connected Riemann surface of
genus $g\geq 2$ and let $\Gamma_g$ be its fundamental group, then
$\Sigma_g\simeq \mathcal{H}/\Gamma_g$ where $ \mathcal{H}$ denotes
the Poincar\'e upper halfplane ($\{ z\in \CC,\;/ {\rm Im}\, z>
0\} ).$ The projection map $p: \mathcal{H}\rightarrow
\mathcal{H}/\Gamma_g$ defines the universal covering of
$\Sigma_g$.\\ {\it From now on all our $\Gamma$-coverings will be
Galois.} Recall that $\Gamma$-coverings are, in particular,
principal $\Gamma$-bundles. Thus, thanks to the classification
theorem for principal bundles, see Lawson-Michelson \cite{LM}, we know that there exist
topological spaces $B\Gamma$, $E\Gamma$, with $E\Gamma$
contractible,   and a $\Gamma$-covering $E\Gamma \rightarrow B
\Gamma$ such that the following statement holds:

\medskip
\n
{\it there is a natural bijection between the set of isomorphism
classes of $\Gamma-$coverings on $M$ and the set of homotopy
classes of continuous maps $r:M \rightarrow B \Gamma$.}

\m
\n
The bijection is realized by the map that associates to $(M,r:M\to
B\Gamma)$ the $\Gamma$-covering  $r^* E\Gamma$. The space $B
\Gamma$ is uniquely defined up to homotopy equivalences and is
called the {\it classifying space} of $\Gamma$. The map $r$ is
called the classifying map. In the example above one has
$E\Gamma_g=\mathcal{H}, B\Gamma_g=\Sigma_g$ and $r=$ identity. As
a different example: $B\ZZ^k=(S^1)^k,$ $E \ZZ^k=\RR^k$ with
covering map: $$ (x_1, \ldots , x_k) \in \RR^k \rightarrow ( e^{i
x_1}, \ldots ,e^{i x_k} ). $$
\\ From
now on we shall identify a $\Gamma$-covering with the
corresponding pair $(M,r:M\to B\Gamma)$.

\smallskip
\begin{definition}\label{def:homotopy-equivalent}
Let $M$ and $M^\prime$  be  closed oriented manifolds. We shall
say that  two $\Gamma$-coverings $$(M,r:M\to B\Gamma)
\quad\text{and} \quad(M^\prime,r^\prime:M^\prime\to B\Gamma)$$ are
{\it oriented homotopy equivalent} if there exists an oriented
homotopy equivalence $h:M^\prime\to M$ such that $r\circ h \simeq
r^\prime$, where $\simeq$ means homotopic.
\end{definition}

\smallskip
\begin{definition}\label{def:cut-and-paste}
Let  $M$ and $N$ be two oriented compact 
manifolds with boundary and let $\phi,\psi: \partial M \to
\partial N$ be orientation preserving diffeomorphisms. Let
$r: M\cup_{\phi} N^- \to B\Gamma$ and $s: M\cup_{\psi} N^- \to
B\Gamma$ be two reference maps. We say that they define   {\it
cut-and-paste equivalent} $\Gamma$-coverings if  $r|_M \simeq
s|_M$ and $r|_N \simeq s|_N$ holds, where  $\simeq$ means
homotopic.
\end{definition}

\n
Geometrically this means that $r^* E\Gamma\to M\cup_{\phi} N^-$
and $r^* E\Gamma\to M\cup_{\psi} N^-$ give rise to isomorphic
bundles when restricted to $M$ and $N$ respectively.

\subsection{The definition of higher signatures.}$\;$

\medskip
Let  $\Gamma\to \widetilde{M} \to M$ be a $\Gamma$-covering of a
closed oriented manifold and let $r:M\to B\Gamma$ be a classifying
map for such a covering. Consider the cohomology of $B\Gamma$ with
real coefficients $H^*(B\Gamma,\RR)$. It can be proved that there
is a natural isomorphism $$H^* (B\Gamma,\RR)\cong H^*
(\Gamma,\RR)$$ where on the right hand side we have the algebraic
cohomology of the group $\Gamma$. We recall that  $H^*(\Gamma ,
\RR)$ is by definition the graded homology group associated to the
complex $\{ \mathcal{C}^*(\Gamma) , d\}$ whose $p-$cochains are
functions $c:\Gamma^{p+1} \rightarrow \RR$ satisfying the
invariance condition $$ c(g\cdot g_0,\dots, g\cdot g_p)= c(g_0,
\dots, g_p)\quad\forall g,g_0,\dots, g_p\in \Gamma, $$ and with
coboundary given by the formula $$ (d c)(g_0, \cdots, g_{p+1})=
\sum_{i=0}^{p+1} (-1)^i c(g_0, \cdots, g_{i-1}, g_{i+1},\cdots,
g_{p+1}). $$ Since we deal with real coefficients, the above
complex can be replaced by the subcomplex of antisymmetric
cochains: $$ \forall \tau \in S_{p+1},\; c( g_{\tau(0)}, \cdots,
g_{\tau(p+1)}) = {\rm sgn}\,(\tau)\, c(g_0, \cdots; g_{p+1}). $$
Let us fix a class $[c]\in H^* (B\Gamma,\RR)$. We take the
pull-back $r^* [c]\in H^* (M,\RR)$ and consider
\begin{equation}\label{def-of-higher-sign}
{\rm sign} (M,r;[c])\,:=\,<L(M)\cup r^* [c],[M]>
\end{equation}
This real number is called the {\it  Novikov higher signature}
associated to $[c]\in H^* (B\Gamma,\RR)$ and the classifying map
$r$. Using the de Rham isomorphism we can equivalently write
\begin{equation}\label{def-of-higher-sign2}
{\rm sign} (M,r;[c])\,:=\,\int_M [L(M,\nabla^g)]\wedge r^* [c]\,.
\end{equation}
If $\dim M=4k$ and $[c]=1\in H^0 (B\Gamma,\RR)$, then $$ {\rm
sign} (M,r;1)=\int_M L(M)(={\rm sign} (M))$$ and we reobtain the lower
signature.

\medskip
\n
{\bf Remark.} We have defined the Hirzebruch $L$-class as the de
Rham class of the $L$-form $L(M,\nabla^g)$.
In fact, using a more
topological approach to characteristic classes, one can define the
$L$-class in $H^* (M,\QQ)$; consequently the higher signatures
${\rm sign}(M,r;[c])$ can be defined for each $[c]\in H^*
(B\Gamma,\QQ)$.


\medskip
For motivation and historical remarks concerning Novikov  higher
signatures the reader is referred to the survey by Ferry-Ranicki-Rosenberg
\cite{FRR}.

\section{{\bf Three fundamental
questions.}}\label{sect:3questions}

Having defined the higher signatures $$\{\,\;{\rm
sign}(M,r;[c])\,,\;\;[c]\in H^* (B\Gamma,\RR)\,\;\}$$ and keeping in mind 
the properties of the lower signature, we can ask
the following three fundamental questions.

\bigskip
\n
{\bf Question 1.} Are the higher signatures
{\it homotopy
invariant} ?

\smallskip
\n
{\bf Question 2.} Are the higher signatures {\it cut-and-paste
invariant} ?

\smallskip
\n
{\bf Question 3.} If $\pa M\not= \emptyset$, {\it can we define
higher signatures and prove their homotopy invariance ?} Of course
we want these higher signatures on a manifold with boundary $M$ to
generalize the lower signature $$\int_M L(M,\nabla^g) -{1\over 2}
\eta(D^{{\rm sign}}_{(\pa M,g_\pa )})\,,$$ which is indeed a
homotopy invariant by Theorem \ref{prop:hominv}.

\bigskip
We anticipate our answers: Question 1 is still open and is known
as the {\it Novikov conjecture}. It has been settled in the
affirmative for many classes of groups. For instance, the
following groups satisfy the Novikov conjecture:  virtually
nilpotent groups and more generally  amenable groups, any discrete
subgroup of $GL_n(F)$ where $F$ is a field of characteristic zero,
Artin's braid groups $B_n$, one-relator groups, the discrete
subgroups of Lie groups with finitely many path components,
$\pi_1(M)$ for a complete Riemanniann manifold with non-positive
sectional curvature. The Novikov conjecture has also been proved
for hyperbolic groups  and, more generally, for groups acting
properly on bolic spaces (see the recent work of Kasparov and
Skandalis). A few relevant references are  Mishchenko
\cite{Mish1},
Kasparov \cite{kasparov-conspectus},
\cite{kasparov-topological}, \cite{Kasparov(1993)}, Weinberger
\cite{Weinberger(survey)},
Connes-Moscovici
\cite{CM}, Connes-Gromov-Moscovici \cite{CGM},
\cite{CGM2}, Ferry-Ranicki-Rosenberg \cite{FRR}, Gromov
\cite{Gromov-macroscopic},
Higson-Kasparov, \cite{Higson-Kasparov(2001)}, Kasparov-Skandalis \cite{Ka-Ska},
Guentner-Higson-Weinberger \cite{GHW}. For related material see also
Lafforgue \cite{lafforgue}, Cuntz \cite{Cuntz}, Mathai
\cite{Mathai}, L\"uck-Reich \cite{Reich}, Schick \cite{schick-survey}.

The answer to Question 2
is negative: the higher signatures are {\it not} cut-and-paste
invariants (we shall present a counterexample below). However, one
can give sufficient conditions on the separating hypersurface $F$
and on the group $\Gamma$ ensuring that the higher signatures are
indeed cut-and-paste invariant.

Finally, under suitable assumption on $(\pa M,r|_{\pa M})$ and on
the group $\Gamma$ one can {\it define} higher signatures on a
manifold with boundary $M$ equipped with a classifying map $r:M\to
B\Gamma$ and {\it prove their homotopy invariance}. 

Relevant references for the solution to the last 2 questions will be given 
along the way.


\section{{\bf The Novikov conjecture on closed manifolds: the
K-theory approach.}}\label{sect:k-theory-approch}

In this section we shall describe one of the approaches that have
been developed in order to attack, and sometime solve, the Novikov
conjecture. We begin by introducing important mathematical objects
associated to $M$, $\Gamma$ and $r:M\to B\Gamma$.

\subsection{The reduced group $C^*$-algebra $C^*_r \Gamma$.} $\;$

\m

We consider the group ring $\CC\Gamma$. It can be identified with
the complex-valued functions on $\Gamma$ of compact support. Any
element $f\in \CC\Gamma$ acts on $\ell^2 (\Gamma)$ by left
convolution. The action is bounded in the $\ell^2$ operator norm
$\|\cdot\|_{\ell^2 (\Gamma)\to \ell^2 (\Gamma)}\,.$ The reduced
group $C^*$-algebra, denoted $C^*_r \Gamma$, is defined as the
completion of $\CC\Gamma$ in  $B(\ell^2 (\Gamma))$. Let us give an
example: if $\Gamma=\ZZ^k$ then using Fourier transform one can
prove that there is a natural isomorphism of $C^*$algebras:
$$C^*_r \ZZ^k \longleftrightarrow C^0 (T^k)$$ with $T^k={\rm
Hom}(\ZZ^k,U(1))$ the dual group associated to $\ZZ^k$ (a
$k$-dimensional torus).

\subsection{K-Theory.}$\;$

\m

Let $A$ be a unital $C^*$-algebra,  such as $C^*_r \Gamma$. We
recall that $K_0 (A)$ is defined as the  group generated by the
{\it stable} isomorphism classes of
 finitely generated projective left $A-$modules; more precisely
such a module is the range of a projection $p$ in a matrix algebra
$M_n(A)$  and one identifies two pairs of projections $(p,q) \in
M_n(A)^2$ and $(p^\prime, q^\prime)\in M_{n^\prime} (A)^2$ if for
suitable $k, k^\prime \in \NN$, $$p\oplus  q^\prime \oplus {\rm
Id}_k \oplus 0_{k^\prime}\;\;\text{is conjugate to}\;\; p^\prime
\oplus q \oplus {\rm Id}_k \oplus 0_{k^\prime}\;\;\text{in}\;\;
M_{n+n^\prime+k+k^\prime}(A ).$$ One then denotes by $[p-q]$
(=$[p^\prime-q^\prime]$) the class of $(p,q)$; similarly, if $E$
and $F$ are finitely generated projective left $A$-modules, then
we denote by $[E-F]$ the associated class in $K_0 (A)$. $K_0 (A)$
is an additive group. When $A$ is a non unital $C^*-$algebra one
introduces the unital $C^*-$algebra $\widetilde{A}= A \oplus \CC$
obtained by adding the unit element $0 \oplus 1$ to $A$; one
considers the morphism $\epsilon: \widetilde{A} \rightarrow \CC$
defined by $\epsilon (a \oplus \lambda) = \lambda.$ One then
defines $K_0(A)$ to be equal to the kernel of the map $\epsilon_*:
K_0(\widetilde{A}) \rightarrow K_0(\CC)$ induced by $\epsilon$.
Observe that $K_0(\CC)= K_0( M_n(\CC) ) = \ZZ.$ We also define
$K_1(A)$ to be equal to $K_0( A \otimes C_0(\RR) )$ where $ A
\otimes C_0(\RR)$ is the {\it suspension of $A$}. For instance
$K_1(\CC)= K_1(M_n(\CC) )= 0.$ Alternatively,  $K_1(A)$ can be
identified with the set of connected components of $GL_\infty(A).$
We  recall that for any compact
Hausdorff space 
$M,$ the 
$K-$theory group $ K^0(M)$  is defined as the set of formal 
differences of isomorphism classes of complex vector bundles 
over $M.$ Then Swann's theorem states that   $K_0(C^0(M) )$ is 
isomorphic to  $ K^0(M).$  Thus, from the
previous sub-section one gets an isomorphism:
 $K_0 (C^*_r \ZZ^k) \simeq K^0 (T^k)$.

\subsection{The index class of the signature operator in $K_* (C^*_r \Gamma)$ .}\label{subsect:index-class-closed}$\;$

\subsubsection{{\bf  $C^*_r
\Gamma$-linear operators.}}$\;$

\m

Let $(M,g)$ be a closed, compact and oriented Riemannian manifold.
Let $\pi:\widetilde {M} \to M$ be a Galois $\Gamma$-covering. Let
$r:M\to B\Gamma$ be a classifying map for this covering. We
consider a Hermitian Clifford module $E\rightarrow M$, endowed 
with a Clifford connection $\nabla^E,$ and let $D$ be the associated
 Dirac-type operator acting on $C^\infty(M,E).$
 For example we could consider
the signature operator associated to $g$ and our choice of
orientation. 
Notice that we can lift the operator $D$ to a $\Gamma$-invariant
differential operator $\widetilde{D}$ on $\widetilde{M}$;
$\widetilde{D}$ acts on the section of the $\Gamma$-equivariant
bundle $\widetilde{E}:=\pi^* E$. 
Consider now $C^*_r \Gamma$. The
group $\Gamma$ acts in a natural way on $C^*_r \Gamma$ by right
translation. It also act on $\widetilde{M}$ (on the left) by deck
transformations: we can therefore consider the associated bundle
$$\mathcal{V} := C^*_r \Gamma\times_\Gamma\widetilde{M} $$ which
is a vector bundle with typical fiber $C^*_r \Gamma$. We shall be
interested in the space of sections $C^\infty (M,E\otimes
\mathcal{V})$. If ${\rm rank \, E}=N$ and $s\in C^\infty (M,E\otimes
\mathcal{V})$, then in a trivializing neighborhood $U$ we can
identify $s|_U$ with a $N$-tuple of $C^*_r \Gamma$-valued
functions $(s^1,\dots,s^N)$. This shows that $C^\infty (M,E\otimes
\mathcal{V})$ is in a natural way a left $C^*_r \Gamma$-module.
 Moreover, using the Hermitian metric $h(\cdot,\cdot)$ on $E$
we can define a $C^*_r \Gamma$-valued inner product
$<\cdot,\cdot>$: if $s,t\in C^\infty_0 (U,
(E\otimes\mathcal{V})|_{U})$ then $$<s,t>:=\int_U \sum h_{ij}s^i
t^j {\rm dvol}_g\,\,\in \,C^*_r \Gamma\,$$ The general case is
obtained by using a partition of unity. $C^\infty (M,E\otimes
\mathcal{V})$ equipped with the above $C^*_r \Gamma$-valued inner
product is a pre-Hilbert $C^*_r \Gamma$-module, in the sense that
it satisfies the following properties: $
 \forall a\in C^*_r \Gamma, \forall s,t,u\in C^\infty (M,E\otimes
 \mathcal{V})$:
$$ <s,t+u>=<s,t>+<s,u>, \;\;\;<a\cdot s,t>=a\cdot <s,t>, \;\;\; <s
, a\cdot t > = a^* <s , t >. $$  The completion of $C^\infty
(M,E\otimes \mathcal{V})$ with respect to the norm $$
\|s\|=\sqrt{\|<s,s>\|_{C^*_r \Gamma}}$$ is denoted by $L^2_{C^*_r
\Gamma} (M,E\otimes \mathcal{V})$; it is a {\it Hilbert $C^*_r
\Gamma$-module}.

\smallskip
The product bundle
 $ C^*_r \Gamma \times \widetilde{M} \rightarrow \widetilde{M}$ is 
 endowed with the trivial flat connection. It induces a 
 (non trivial) flat connection $\nabla^{\mathcal{V}}$ on 
 the $C^*_r \Gamma-$bundle $\mathcal{V}.$ Then the bundle 
 $E\otimes \mathcal{V} \rightarrow M$ is endowed with the 
 connection $\nabla^E \otimes \Id + \Id \otimes \nabla^{\mathcal{V}}.$
 We denote by $\D_{(M,r)}$ the associated twisted  
 Dirac type operator. Directly from the definition we see that:  
\begin{equation}\label{C*-linearity}
\D_{(M,r)}: C^\infty (M,E\otimes \mathcal{V})\rightarrow C^\infty
(M,E\otimes \mathcal{V})\;\;\;\text{is}\;\;\;C^*_r
\Gamma-\text{linear}
\end{equation} A good reference for seeing
the details of this approach is Schick \cite{schick03}.
 We also remark that it is possible to
introduce Sobolev $C^*_r \Gamma$-modules $H^m_{C^*_r \Gamma}
(M,E\otimes \mathcal{V})$ and $\D_{(M,r)}$ extends to a bounded
$C^*_r \Gamma$-linear operator from $H^1_{C^*_r \Gamma}
(M,E\otimes \mathcal{V})$ to $L^2_{C^*_r \Gamma} (M,E\otimes
\mathcal{V})$.

 If $M$ is even dimensional, then $E$ is $\ZZ_2$-graded,
$E=E^+\oplus E_-$; thus $$\D_{(M,r)}=
\begin{pmatrix} 0 & \D_{(M,r)}^- \cr \D_{(M,r)}^+ & 0
\end{pmatrix} \,,$$ with $\D_{(M,r)}^+ $ and $\D_{(M,r)}^-$ both
$C^*_r \Gamma$-linear.

\subsubsection{{\bf The index class in $K_*( C^*_r \Gamma)$.}}$\;$

\m

From (\ref{C*-linearity}) we infer that $\Ker \D_{(M,r)}^+$ and
${\rm coker} \D_{(M,r)}^+$  are both $C^*_r \Gamma$-modules. In general,
the modules $\Ker \D_{(M,r)}^+$ and ${\rm coker} \D_{(M,r)}^+$ are not finitely generated 
and projective so they cannot be used directly to define 
the index class $\Ind (\D_{(M,r)}^+)\in K_0 (C^*_r \Gamma)$ as
$[\Ker \D_{(M,r)}^+]- [{\rm coker} \D_{(M,r)}^+]$.
\footnote{This is similar to the problem one encounters in defining the
index class of a family $\mathcal{F}:= (F_\theta)_{\theta\in T}$ of Fredholm 
operators parametrized by
a space $T$: the kernel-bundle and the cokernel-bundle do  not
in general vary continuously.}
However this is true up to a smoothing perturbation $\mathcal{R}$;
one defines the index class as
$$\Ind (\D_{(M,r)}^+):= [\Ker \D_{(M,r)}^+ +\R]-[{\rm coker} \D_{(M,r)}^+ +\R]\in K_0 (C^*_r \Gamma)\,$$
(and the definition does not depend on the choice of $\R$).
Let us see the details. 



\smallskip
\n
{\bf The Mishchenko-Fomenko pseudodifferential calculus.} 
One can define a space of $C^*_r \Gamma$-linear
{\it differential} operators ${\rm Diff}^*_{C^*_r \Gamma}
(M;E\otimes \mathcal{V},E\otimes \mathcal{V})$; these are simply
operators locally given by a $N\times N$-matrix $A_{ij}$, $N={\rm
rk}E$, with $$A_{ij}=\sum_{|\alpha|\leq k} a(ij)_{\alpha}
{\pa^\alpha\over
\pa_{x_1}^{\alpha_1}\cdots\pa_{x_n}^{\alpha_n}}\,,\quad
\text{with}\quad  a(ij)_{\alpha}\in C^\infty (U,C^*_r \Gamma)\,.$$
In a very natural way we can give the notion of ellipticity in
${\rm Diff}^*_{C^*_r \Gamma} (M;E\otimes \mathcal{V},E\otimes
\mathcal{V})$. From the definitions, we discover first of all that
$\D_{(M,r)}\in {\rm Diff}^1_{C^*_r \Gamma} (M;E\otimes
\mathcal{V},E\otimes \mathcal{V})$; moreover the ellipticity of
$D$ implies that $\D_{(M,r)}$ is elliptic in ${\rm Diff}^1_{C^*_r
\Gamma} (M;E\otimes \mathcal{V},E\otimes \mathcal{V})$ .
Mishchenko and Fomenko have developed a pseudodifferential
calculus for $C^*_r \Gamma$-linear operators $$\Psi^*_{C^*_r
\Gamma} (M;E\otimes \mathcal{V},E\otimes \mathcal{V})\supset {\rm
Diff}^*_{C^*_r \Gamma} (M;E\otimes \mathcal{V},E\otimes
\mathcal{V})\,.$$ Using this calculus one can prove that given an
elliptic operator  $\P\in {\rm Diff}^k_{C^*_r \Gamma} (M;E\otimes
\mathcal{V},E\otimes \mathcal{V})$, it is possible to find an
inverse  $\Q\in \Psi^{-k}_{C^*_r \Gamma} (M;E\otimes
\mathcal{V},E\otimes \mathcal{V})$ modulo elements in
$\Psi^{-\infty}_{C^*_r \Gamma}(M;E\otimes \mathcal{V},E\otimes
\mathcal{V})$. Notice that the smoothing operators in the
Mishchenko-Fomenko calculus are simply the integral operators with
a Schwartz kernel on $M\times M$ locally given by a smooth
function  with values in $M_{N\times N}(C^*_r \Gamma)$. \\Let in
particular $M$ be even dimensional and let $E=E^+ \oplus E^-$
be the $\ZZ_2$-graded bundle appearing in the definition of our
Dirac operator. The operator $$\D^{+}_{(M,r)}\in {\rm
Diff}^*_{C^*_r \Gamma} (M;E^+\otimes \mathcal{V},E^-\otimes
\mathcal{V})$$ is elliptic and there exists therefore a parametrix
$\Q\in \Psi^{-1}_{C^*_r \Gamma} (M;E^-\otimes
\mathcal{V},E^+\otimes \mathcal{V})$ such that
\begin{equation}\label{parametrix}
\D^{+}_{(M,r)}\circ \Q={\rm Id}-\R_-\,,\;\; \Q\circ
\D^{+}_{(M,r)}={\rm Id}- \R_+
\end{equation}
with $R_\pm\in \Psi^{-\infty}_{C^*_r \Gamma} (M;E^\pm\otimes
\mathcal{V},E^\pm\otimes \mathcal{V})$.
 This part of the  theory runs
quite  parallel to the usual case, when the $C^*$-algebra is equal
to $\CC$; the main differences arise in the functional analytic
consequences of (\ref{parametrix}). The point is that doing
functional analysis on a Hilbert $A$-module, with $A$ a
$C^*$-algebra, is a more delicate matter than doing functional
analysis on a Hilbert space (see  Wegge-Olsen \cite{WO} and
Higson \cite{Higson} for more on this delicate point).

 \m\n
{\bf The Mishchenko-Fomenko decomposition theorem.}
On a Hilbert $A$-module there exists a natural
 notion of $A$-compact operator:
   using (\ref{parametrix}), elliptic regularity and the fact that elements
  in $\Psi^{-\infty}_{C^*_r \Gamma}$ are $C^*_r
  (\Gamma)$-compact on $L^2_{C^*_r \Gamma}$,
  one can prove a decomposition of
the space of sections
  of $E\otimes \mathcal{V}$ with respect to $\D^+_{(M,r)}$, i.e.
  \begin{equation}\label{Mishchenko_decomposition_closed}
    C^\infty(M,E^+\otimes \mathcal{V}) = \mathcal{I}_+ \oplus
    \mathcal{I}_+^{\perp},\quad C^\infty(M,E^-\otimes
    \mathcal{V})=\mathcal{I}_- \oplus
    \D_{(M,r)}^+(\mathcal{I}_+^{\perp})\,,
  \end{equation}
  with  $\mathcal{I}_+$ and $\mathcal{I}_-$
   {\it finitely generated projective} $C^*_r \Gamma$-modules.
Notice that the second decomposition is not, a priori, orthogonal.
  However, $\D_{(M,r)}^+$ induces an isomorphism (in the Fr\'echet
  topology) between $\mathcal{I}_+^{\perp}$ and
  $\D_{(M,r)}^+(\mathcal{I}_+^{\perp})$.
Intuitively $\mathcal{I}_+$ should be thought of
as the kernel of $\D_{(M,r)}^+$ and $\mathcal{I}_-$ as the
cokernel.

\s
\n
{\bf The index class.}
 The index class of $\D^+_{(M,r)}$, \`a la
Mishchenko-Fomenko, is precisely given by
\begin{equation}\label{index class closed}
\Ind (\D^+_{(M,r)})= [\mathcal{I}_+] - [\mathcal{I}_-]\;\in\;K_0
(C^*_r \Gamma )\,.
\end{equation}
Although the decomposition (\ref{Mishchenko_decomposition_closed})
is not unique, the index class is uniquely defined in $K_0 (C^*_r
\Gamma)$ . The main reference for  this material is the original
article of Mishchenko and Fomenko \cite{MF}; see also
\cite[Appendix A]{LPMEMOIRS}. Working a little bit more one can
show that the orthogonal projection $\Pi_+$ onto $\mathcal{I}_+$
and the projection $\Pi_-$ onto $\mathcal{I}_-$ along
$\D^+(\mathcal{I}_+^\perp)$ are elements in $\Psi^{-\infty}_{C^*_r
\Gamma}$ (see \cite{LPMEMOIRS}, Appendix A) . Thus
$$\D^+_{(M,r)}-\Pi_- \D^+_{(M,r)}\Pi_+$$ is a {\it smoothing}
perturbation of $\D^+_{(M,r)}$ with the property that its kernel
and cokernel are finitely generated and projective. \\{\it Summarizing}:
there exists a smoothing perturbation $\R$ of
 $\D^+_{(M,r)}$ such that $\Ker (  \D^+_{(M,r)}+\R)$ and
 ${\rm coker} ( \D^+_{(M,r)}+\R) $ are finitely generated
 projective $C^*_r \Gamma$-modules; the index class can be  defined as
$$\Ind ( \D^+_{(M,r)}):=[\Ker (  \D^+_{(M,r)}+\R)]-[{\rm coker} ( \D^+_{(M,r)}+\R)]
\;\in\; K_0 (C^*_r \Gamma)$$
and it is not difficult to prove that it does not depend
 on the choice of $R\in \Psi^{-\infty}_{C^*_r \Gamma}$.

If $M$ is odd dimensional, then the Clifford module $E$ will be
ungraded; we obtain in this case an index class $\Ind
\D_{(M,r)}\in K_1 (C^*_r \Gamma)$. We shall not give the details
here.

\subsubsection{{\bf The example
$\Gamma=\ZZ^k$.}}\label{subsubsect:abelian-class}$\;$

Let $N$ be a closed oriented manifold with $\pi_1(N)=\ZZ^k$. Let
$r$ be the classifying map.
 In this case the
higher index class $\Ind(\D^{{\rm sign}}_{(N,r)})$ has, thanks to
 Lustzig \cite{Lu},
 a geometric description 
Details for the material that follows
can be found in \cite{Lu} and Lott \cite{Lott II}. As already remarked the space
$B\ZZ^k$ is a $k$-dimensional torus; more precisely, it is the
dual torus $(T^k)^*$ to $T^k=\widehat{\ZZ^k}=\Hom(\ZZ^k,U(1))$. Using the duality
between the two tori it is easy to see that on
the product $(T^k)^*\times T^k$ there is a canonical Hermitian
line bundle $H$ with a canonical Hermitian connection $\nabla^H$.
The bundle $H$ is flat when restricted to any fibre of the
projection $(T^k)^*\times T^k \rightarrow T^k$.  Using the map
$r\times\id:N\times T^k \rightarrow (T^k)^*\times T^k$ we obtain a
line bundle $F$ on $N\times T^k$ with a natural Hermitian
(pulled-back) connection $\nabla^{F}$. In this way we have
obtained a fibration of  closed manifolds $\phi:N\times
T^k\rightarrow T^k$ and a Hermitian line bundle $F$ over the total
space with a flat structure in the fibre directions. Let
$\theta\in T^k$ and let $F_\theta$ be the restriction of $F$ to
$N\times \{\theta\}$. Since $F_\theta$ is flat, the de Rham
differential can be extended to act on $\Lambda^*(M)\otimes
F_\theta$; we obtain a twisted de Rham differential $d_\theta$.
Let $D^{{\rm sign}}_\theta$ be the corresponding twisted signature
operator on $N$. As $\theta$ varies in $T^k$, we obtain a smoothly
varying family of twisted signature operators. Thus, according to
Atiyah and Singer \cite{AS IV}, we obtain an index class
$\Ind(\{D^{{\rm sign}}_\theta \}_{\theta\in T^k})\in K^* (T^k),$
with  $*={\rm dim} N.$ It can be proved that $$\Ind(\D^{{\rm
sign}}_{(N,r)})\in K_* (C^*_r (\ZZ^k))\quad\text{and}\quad
\Ind(\{D^{{\rm sign}}_\theta \}_{\theta\in T^k})\in K^* (T^k)$$
{\it corresponds} under the isomorphisms $ K_* (C^*_r (\ZZ^k))
\simeq K_*(C^0(T^k))\simeq K^* (T^k)\,.$

\subsection{The symmetric signature of Mishchenko.}\label{subsect:symmetric-closed}$\;$

\m

Let $A$ be
an involutive algebra and let us introduce  $L^0 (A),$ the Witt
group of non singular Hermitian forms on $ A$: it classifies
Hermitian forms $Q$ on finitely generated  left projective modules
on $ A.$ Given $E$ a finitely generated left projective module
over $A$, a Hermitian form $Q$ on $E$ is a sesquilinear form
$E\times E \rightarrow A$ such that: $$ \forall \xi, \eta \in E,\;
\forall a,b\in A,\;\;\; Q(a\cdot \xi , b\cdot \eta )= a\cdot  Q(
\xi , \eta ) \cdot b^*\,, \quad Q( \xi ,  \eta )^*=Q( \eta ,  \xi
). $$ The form $Q$ is said to be invertible when the map from $E$
to ${\rm Hom}_A( E, A)$ given by $\xi \rightarrow Q(\cdot, \xi)$
is invertible. {\it The Witt group $L^0 (A)$
 is the group generated by the isomorphism classes of
invertible Hermitian forms with the relations: $[Q_1 \oplus Q_2]=
[Q_1] + [Q_2]$, $0=[Q] +[-Q].$} When $A$ is a $C^*-$algebra with
unit then each finitely generated left projective module over $A$
admits an invertible Hermitian form $Q$ satisfying the positivity
condition $Q(\xi, \xi) \geq 0$ for any $\xi \in E$. (Recall 
that an element $x$ of the $C^*-$algebra $A$ is positive if and only if it is 
of the form $x= y y ^*,$ or equivalently, $x$ is self-adjoint 
and its spectrum lies in $[0, +\infty[$ ). Moreover, on
$E$ all such positive Hermitian forms are pairwise isomorphic so
that there is a well defined map $K_0(A) \rightarrow L^0(A)$
sending $E$ to $(E,Q)$ with $Q$ an invertible positive Hermitian
form on $E$; it turns out that this map is an isomorphism.

Let $M$ be an  oriented $2n-$dimensional manifold and let $r: M
\rightarrow B \Gamma$ be a (continuous) reference map. We are
going to recall, following Mishchenko \cite{Mish1}, \cite{Mish2} (see also
Kasparov \cite{Kasparov(1993)}, \cite{kasparov-conspectus}),
 the construction of the Mishchenko symmetric signature  $\sigma_{\CC\Gamma}(M,r)\in L^0 (\CC\Gamma)$.

Denote by $\widetilde{M}\rightarrow M$ the associated Galois $\Gamma-$covering.
Take a (suitably nice) triangulation of $M$ and pull it back
to $ \widetilde{M}$ to a $\Gamma-$invariant triangulation of $ \widetilde{M}.$
 Let $(C_*, \partial_*)$ and $(C^*, \delta_*)$ denote the associated simplicial chain  complex and
cochain complex:
$
\delta_j :C^j \rightarrow C^{j+1},
$
$
\partial_j: C_{j+1} \rightarrow C_j,
$  for $0\leq j \leq 2n.$ The $C^j, C_j$ are finitely generated free left
$\CC[\Gamma]-$modules.
There is a chain map $\xi_j: C^j \rightarrow C_{2n-j},$ $0\leq j \leq 2n,$
 defining Poincar\'e duality, which satisfies $\partial_{2n-j} \xi_j = (-1)^j
\xi_{j+1} \delta_j$ and induces a chain homotopy equivalence. It
can be arranged that $\xi^*_j =(-1)^j \xi_{2n-j}$ where $\xi^*_j$
denotes the adjoint of the left $\CC[\Gamma]-$linear map
$\xi^*_j.$ We can add a
 $(\CC[\Gamma])^k$,  for a suitable $k\in \NN$, to
both $C^{2n-1}$ and $C_1$ to make $\delta_{2n-1}: C^{2n-1} \rightarrow C^{2n}$
surjective.  Then we add  $((\CC[\Gamma])^k)^*$ to both $C^1$ and
$C_{2n-1}$ in order to preserve Poincar\'e duality and modify accordingly
$\delta_0$ and $\partial_{2n}.$ Now we may split off $\xi_0$ and
$\xi_{2n}$ and repeat this algebraic
 surgery process so as to come down to a complex concentrated in middle
dimension: $\xi_{n}: C^{2n}\rightarrow C_{2n},$ $(i^n \xi_n)^*= i^n \xi_n. $
Since $ i^n \xi_n$ defines a non degenerate Hermitian form
one gets an element, denoted $\sigma_{\CC\Gamma} (M,r)$,
 of $L^0( \CC[\Gamma]). $ {\it Mishchenko has shown that  $\sigma_{\CC\Gamma} (M,r)$
depends only on the oriented homotopy type of $(M,r).$}

\s

Consider now the natural homomorphism $L^0 (\CC\Gamma)\to L^0
(C^*_r \Gamma)$ induced by the inclusion $\CC\Gamma\hookrightarrow
C^*_r \Gamma$. We compose it with the inverse isomorphism $  K_0
(C^*_r \Gamma)\leftarrow L^0 (C^*_r \Gamma)$ and get a well
defined homomorphism $$J:L^0 (\CC\Gamma)\rightarrow K_0 (C^*_r
\Gamma)\,.$$ Let $$\sigma (M,r)\,:=\,J(\sigma_{\CC\Gamma}
(M,r))\in K_0 (C^*_r \Gamma);$$ $\sigma(M,r)$ is the {\it $C^*_r
\Gamma$-valued Mishchenko symmetric signature.} It is a homotopy
invariant of the pair $(M,r:M\to B\Gamma)$.

\subsection{ Homotopy invariance of the index
class.}\label{subsect:hominv-index-closed}$\;$

\m

The following theorem will play a crucial role both in the
treatment of the Novikov conjecture and of the cut-and-paste
invariance of higher signatures. It is due to Mishchenko and
Kasparov, \cite{Mish2}, \cite{kasparov-conspectus}:

\begin{theorem}\label{theo:sign=ind}
Let $(M,r:M\rightarrow B\Gamma)$ be an oriented manifold with
classifying map $r$. Then the index class $\Ind (\D^{{\rm
sign}}_{(M,r)})\in K_* (C^*_r \Gamma)$, $*=\dim M$, is equal to
$\sigma (M,r)$, the $C^*_r \Gamma$-valued Mishchenko symmetric
signature:
\begin{equation}\label{sign=ind}
\Ind (\D^{{\rm sign}}_{(M,r)})=\sigma (M,r)\;\;\;\text{in}\;\;\;
K_* (C^*_r \Gamma)\,.
\end{equation}
\end{theorem}

As a corollary we then get the following fundamental information:

\begin{corollary}\label{cor:hominv-closed}
The index class $\Ind (\D^{{\rm sign}}_{(M,r)})\in K_* (C^*_r
\Gamma)$ is an oriented homotopy invariant.
\end{corollary}

\n
{\bf Remark.} It is possible to give a purely analytic proof of
Corollary \ref{cor:hominv-closed}. This important result is due to
Hilsum and Skandalis \cite{H-S}. See also the work of
Kaminker-Miller \cite{KM}.

\n
{\bf Remark.} When $\Gamma=\ZZ^k$, Lusztig was the first to
establish the homotopy invariance of $\Ind (\D^{{\rm
sign}}_{(M,r)})\in K_* (C^*_r \ZZ^k).$ The proof of
Kaminker-Miller \cite{KM} cited above is an extension of  Lusztig's proof to
the noncommutative context.

\smallskip
\n
{\bf Important remark.} Although Theorem \ref{theo:sign=ind} and
Corollary \ref{cor:hominv-closed} are extremely interesting
results, they still do not settle in anyway the Novikov
conjecture. In fact, these  results should be viewed as  the
higher analogue of {\it only one} out of the two steps we used in
order to prove that $\int_M L(M)$ is an homotopy invariant. This
step is, more precisely, the homotopy invariance of the signature
and its equality with the index. What we are still missing in the
present higher case is the first step, the one relating $\int_M
L(M)$ to the index. The problem we face now is therefore quite
clear: \bigskip

\n
{\bf Fundamental Problem:} {\it how can one use the homotopy invariance of the
index class $\Ind (\D^{{\rm sign}}_{(M,r)})$ in order to prove the
homotopy invariance of the higher signatures $<L(M)\cup
r^*[c],[M]>, [c]\in H^* (B\Gamma,\RR)$ ?}

\bigskip
\n
 Alternatively: \\{\it
how can we connect the index class and its homotopy invariance to
the higher signatures ?}

\medskip
\n We shall present below two answers to this question. The first
one, due to Kasparov, employs the K-homology of $B\Gamma$, $K_*
(B\Gamma)$, and a natural map $\mu: K_* (B\Gamma) \to K_* (C^*_r
\Gamma)$; the second one, due to Connes and Moscovici, employs
cyclic cohomology.

\subsection{The  assembly map and the Strong Novikov Conjecture.}\label{subsect:operator-k}$\;$

\m

We are considering a closed oriented manifold $M$ and a
classifying map $r:M\to B\Gamma$. Let $L(M)\cap [M]$ be the
Poincar\'e dual to $L(M)$ and consider $r_* (L(M)\cap [M])\in H_*
(B\Gamma,\RR)$. One can check, using some basic algebraic
topology, that $${\rm sign}(M,r;[c])=<[c],r_* (L(M)\cap [M])>\,.$$
Thus the homotopy invariance of the real homology class $r_*
(L(M)\cap [M])$ implies the homotopy invariance of all the higher
signatures $\{ {\rm sign}(M,r;[c]), [c]\in H^* (B\Gamma,\RR)\}$.

It is well known that $K$-theory is a generalized cohomology
theory; it thus admits a dual theory, {\it K-homology}, and there
is a homological Chern character map ${\rm Ch}: K_* (\;\;,\ZZ)\to
H_* (\;\;,\ZZ)$ which is an isomorphism modulo torsion.
Summarizing: the K-homology of $B\Gamma$ is well defined
 and there is an isomorphism ${\rm Ch}^{-1}:H_*
 (B\Gamma,\RR)\rightarrow K_* (B\Gamma)\otimes_\ZZ \RR$.
 Thus we
 are led to consider the following K-homology class
 $${\rm Ch}^{-1}(r_* (L(M)\cap [M])\in K_* (B\Gamma)\otimes_\ZZ
 \RR\,.$$
Clearly: {\it if this class is homotopy invariant, then the
Novikov conjecture is true.}\\ In order to understand why we  wish
to pass from homology  to K-homology we shall simply mention that
besides the abstract definition (a dual theory), there are other
characterizations of K-homology, directly connected to elliptic
operators. Historically, Atiyah was the first to realize that
cycles in $K_* (X)$ should be thought of as "abstract elliptic
operators" \cite{Atiyah-global}. His ideas were further pursued by
Kasparov \cite{kasparov-old} and Baum-Douglas-Fillmore \cite{BDF}.
At the same time, Baum and Douglas \cite{BD} proposed a purely
topological definition of $K$-homology and showed that it was
compatible with the analytic one of Atiyah. We shall present this
topological definition, since it is the easiest to explain and
leads directly to the map $\mu:K_* (B\Gamma)\rightarrow K_* (C^*_r
\Gamma)$ that was mentioned at the end of the previous section. We
shall concentrate on the even dimensional case and pretend that
$B\Gamma$ is compact (the general case is obtained by taking an
inductive limit).\\Cycles in the (topological) K-homology groups
$K_0 (X)$ of a compact topological Hausdorff space $X$ are given
by triples $(M,r:M\to X,E)$ where $M$ is an
even dimensional oriented manifold, 
$r:M\to X$ is continuous, and $E$ is a $\ZZ_2$-graded vector
bundle over $M$ which can be given the structure of graded
Clifford module. \footnote{The original definition of Baum-Douglas
was slightly different: it assumed  $M$ to be spin$_c$ but left
$E$ arbitrary; Keswani has proved, see \cite{Kes3}, that the two
definitions are equivalent.} One then introduces an equivalence
relation on this triples given by {\it bordism}, {\it direct sum}
and  {\it vector bundle modification}. We do not enter into the
details here. The quotient turns out to be the $K_0$-homology
group of $X$. For example $[M,\id,\Lambda^{{\rm sign}}_\CC (M)]$,
with $\Lambda^{{\rm sign}}_\CC
(M)=\Lambda^+_\CC(M)\oplus\Lambda^-_\CC (M)$ the Clifford module
defining the signature operator, is a class in $K_0 (M)$.
Similarly, if $r:M\to B\Gamma$ is a classifying map, then
$[M,r:M\to B\Gamma,\Lambda^{{\rm sign}}_\CC (M)]$ defines an element in $K_0
(B\Gamma)$.

Let now $[M,r:M\to B\Gamma,E^+ \oplus E^-]$ be an element in
$K_0(B\Gamma)$: we define a map
\begin{equation}\label{assembly} \mu:K_0 (B\Gamma)\longrightarrow
K_0 (C^*_r \Gamma)
\end{equation}
by sending $[M,r:M\to B\Gamma,E^+\oplus E^-]$ to the index class,
in $K_0 (C^*_r \Gamma)$,  associated to the $C^*_r \Gamma$-linear
Dirac operator  associated to the Clifford module $E$ and the
classifying map $r:M\to B\Gamma$. Thus if
 $D^{E}$ is the Dirac operator associated to $E$ on $M$ and if, as usual, we
denote by $\D^{E}_{(M,r)}$ the operator $D^{E}$ twisted by the
flat bundle $\mathcal{V}=r^* E\Gamma\times_\Gamma C^*_r \Gamma$,
then the map (\ref{assembly}) is given by $$K_0 (B\Gamma)\ni
[M,r:M\to B\Gamma,E]\xrightarrow{\mu} \Ind \D^{E,+}_{(M,r)}\in K_0
(C^*_r \Gamma)\,.$$ As a fundamental example we have:
$$\mu\,[M,r:M\to B\Gamma,\Lambda^{{\rm sign}}_\CC (M)]=\Ind
\D^{{\rm sign},+}_{(M,r)}\in K_0 (C^*_r \Gamma).$$

\smallskip
A similar map, from $K_1 (B\Gamma)$ to $K_1 (C^*_r \Gamma)$, can
be defined in the odd case, considering odd dimensional manifolds
and ungraded Clifford modules in the definition of the cycles of
$K_1 (B\Gamma)$. We shall denote by $\mu_\RR$ the map induced from
$K_* (B\Gamma)\otimes_\ZZ \RR$ to $K_* (C^*_r \Gamma)\otimes_\ZZ
\RR$. The map $\mu$ is called the {\it assembly map};  it is also
referred to as the {\it Kasparov map}. If $\Gamma$ is torsion free
then it also known as the {\it Baum-Connes map}. One can check,
unwinding the definitions, that $${\rm Ch}^{-1}(r_* (L(M)\cap
[M]))=[M,r:M\to B\Gamma,\Lambda^{{\rm sign}}_\CC (M)]\in K_*
(B\Gamma)\otimes_\ZZ \RR\;.$$  Hence
\begin{equation}\label{mu-poincare=index}
\mu_\RR ({\rm Ch}^{-1}(r_* (L(M)\cap [M])))= \Ind \D^{{\rm
sign}}_{(M,r)}\in K_* (C^*_r \Gamma)\otimes_\ZZ \RR\,.
\end{equation}

\smallskip
We thus arrive at the following fundamental conclusion:

\begin{theorem}\label{theo:snc}
 If the map $\mu_\RR$ is injective then the Novikov
conjecture is true.
\end{theorem}

\begin{proof}
If $(M,r:M\to B\Gamma)$ and $(N,s:N\to B\Gamma)$ are homotopy
equivalent, then by Corollary \ref{cor:hominv-closed} we have
$\Ind \D^{{\rm sign}}_{(M,r)}= \Ind \D^{{\rm sign}}_{(N,s)}$.
Using (\ref{mu-poincare=index}), the injectivity of $\mu_\RR$ and
the bijectivity of $\Ch$ we get $r_* (L(M)\cap [M])=s_* (L(N)\cap
[N])$, which implies the equality of all the higher signatures.
\end{proof}

For later use we notice that the conclusion we can draw is
slightly more general:

\begin{proposition}\label{prop:from-index-to-higher} If $\mu_\RR$ is injective then the equality of the index
classes $\Ind \D^{{\rm sign}}_{(M,r)}= \Ind \D^{{\rm
sign}}_{(N,s)}$  in $K_* (C^*_r \Gamma)\otimes_\ZZ \RR$ implies
the equality of all the higher signatures: $$  < L(M) \cup
r^*(c),[M]>\, =\, <L(N) \cup s^*(c),[N]>\,,\quad\forall c \in
H^*(B \Gamma, \RR)\;. $$
\end{proposition}

This remark
 will be important in treating the cut-and-paste problem for higher
 signatures. Of course, since $\Ind \D^{{\rm sign}}_{(M,r)}=\sigma
 (M,r)$ (the $C^*_r \Gamma$-valued  symmetric signature of
 Mishchenko), we can also state the following
\begin{proposition}\label{from-mish-to-higher} If $\mu_\RR$ is
 injective, then the equality of the $C^*_r \Gamma$-valued
 symmetric signatures, $\sigma(M,r)=\sigma(N,s)$,
implies the equality of all the higher signatures.
\end{proposition}

The injectivity of $\mu_\RR$ (in fact, of $\mu_{\QQ}$)
is known as the {\it Strong Novikov
Conjecture} ($\equiv$ SNC) ; it is still open. Most of the
groups for which the Novikov conjecture has been verified,
satisfy the SNC as well.


We refer to the nice survey of Kasparov
\cite{Kasparov(1993)} for seeing, informally, the technique of the
Dirac-dual Dirac for constructing a left inverse of $\mu_\RR$. See
also 
\cite{kasparov-topological}.

For the connection between the Strong Novikov Conjecture and the
existence of metrics of positive scalar curvature (an important
topic that will be left out of the present survey) we refer, for
example,  to Rosenberg  \cite{Rosenberg0}, \cite{Rosenberg1},
\cite{Rosenberg}, Stolz \cite{St}, Joachim-Schick \cite{JS}, \cite{schick-survey}.

\subsection{The Baum-Connes conjecture} $\;$

\m
The {\it Strong Novikov Conjecture} states that 
the rational assembly map 
$$
\mu_{\QQ}:=\mu \otimes_\ZZ {\id}_\QQ: K_0( B \Gamma )\otimes_\ZZ \QQ \rightarrow 
K_0( C^*_r \Gamma)\otimes_\ZZ \QQ 
$$ is injective. We have just seen that the Strong Novikov conjecture
 implies the Novikov conjecture.

\s

If the discrete group $\Gamma$ is {\it torsion free} then the 
{\it Baum-Connes conjecture} states that the assembly map 
$\mu : K_0( B \Gamma ) \rightarrow K_0( C^*_r \Gamma) $ 
is {\it bijective}. 

When the group $\Gamma$ has torsion then, in general 
the map $\mu$ is not an isomorphism. For instance, 
if $\Gamma = \frac{\ZZ}{ 2 \ZZ }$ then 
$K_0( B \frac{\ZZ}{ 2 \ZZ }) \otimes \QQ =\QQ$ whereas 
$K_0( C^*_r ( \frac{\ZZ}{ 2 \ZZ }) )= \ZZ \oplus \ZZ$ so that
$\mu \otimes_\ZZ {\id}_\QQ$ (and thus $\mu$) cannot be surjective.

\s
In the general non-torsion-free case, 
Baum, Connes and Higson have introduced the space 
$\underline{E} \Gamma$, classifying the 
proper actions of $\Gamma.$ Such a space $\underline{E} \Gamma$ 
is uniquely defined up to $\Gamma-$equivariant homotopy 
(see \cite{BCH}). Baum, Connes and Higson have then constructed an 
assembly map:
$$
\widehat{\mu}: 
K_0^\Gamma ( \underline{E} \Gamma ) \rightarrow K_0( C^*_r \Gamma).
$$ The {\it Baum-Connes conjecture} states that $\widehat{\mu}$ is an {\it isomorphism}.
Notice that there is a natural map 
$\sigma: E \Gamma \rightarrow \underline{E} \Gamma$ 
which induces a map $$\sigma_*: K_0 (B \Gamma)\simeq 
K_0^\Gamma ( {E}\Gamma)\longrightarrow K_0^\Gamma ( \underline{E} \Gamma ).$$ 
The  map $\mu $ is given 
by $\mu =\widehat{\mu} \circ \sigma_*$. One can prove \cite[Section 7]{BCH} that 
the injectivity of $\widehat{\mu}$ implies the rational injectivity of $\mu$. 
In other words, the Baum-Connes conjecture  implies the Strong 
Novikov Conjecture.

For more on the Baum-Connes conjecture 
we also refer the reader to Valette \cite{valette}, Lafforgue
\cite{lafforgue}, L\"uck-Reich \cite{Reich}, Schick \cite{schick-survey}.

\section{{\bf The cyclic-cohomology approach to the Novikov
conjecture.}}\label{sect:connes-moscovici-lott}

Let $M\xrightarrow{r} B\Gamma$ be a closed  oriented manifold with
classifying map $r$. In the previous subsection we have explained
one way to link the index class $\Ind (\D^{{\rm sign}}_{(M,r)})\in
K_* (C^*_r \Gamma)$ (and its homotopy invariance), to the higher
signatures $<L(M)\cup r^* [c], [M]>$, $[c]\in H^*(B\Gamma,\RR)$.
This link is provided by the assembly map $\mu:K_*
(B\Gamma)\rightarrow K_* (C^*_r \Gamma)$. In this subsection we
shall explain a different approach for establishing such a link;
this method, due to Connes and Moscovici \cite{CM}, will use
cyclic cohomology. Our presentation will heavily  employ results
by Lott \cite{Lott I}. In order to understand the main ideas, we
begin by the abelian case, $\Gamma=\ZZ^k$, thus explaining the
{\it seminal work of Lusztig}.

\subsection{The abelian case: family index theory.}\label{subsect:lusztig}$\;$

\m

Let us assume $\Gamma=\ZZ^k$. In subsection \ref{subsect:lusztig} 
we recalled the construction of  the index class $\Ind
(\D^{{\rm sign}}_{(M,r)})\in K_* (C^*_r \ZZ^k)$  in
terms of the {\it index bundle} associated to the Lusztig's family
$\{D^{{\rm sign}}_\theta\}_{\theta\in T^k}$, $T^k= {\rm
Hom}(\ZZ^k,U(1)).$  We briefly
denote this family by $\{D^{{\rm sign}}_\theta\}$.
 The index bundle lives in
$K^*(T^k)$ and we can therefore consider its Chern character $\Ch
(\Ind\,\{D^{{\rm sign}}_\theta \}\,)\in H^*_{{\rm dR}} (T^k).$ An
application of the Atiyah-Singer {\it family index theorem} gives
\begin{equation}\label{chern-lustzig}
\Ch (\Ind\,\{D^{{\rm sign}}_\theta \}\,)=[\int_M L(M)\wedge
\omega]\,\;\;\in\,H^* (T^k),\end{equation} with $\omega$ an
explicit closed form in $\Omega^* (M\times T^k)$.
\smallskip
Let now   $[c]\in H^\ell
(\ZZ^k,\RR)=H^\ell(B\ZZ^k,\RR)=H^\ell((T^k)^*,\RR)$; starting from
$[c]$, Lusztig defines in a natural way $[\tau_c]\in H_\ell
(T^k,\RR)$ so that
\begin{equation}\label{pullback}
r^* [c]= {1\over C(\ell)} <\omega,[\tau_c]>\,,\quad C(\ell)\in
\RR\setminus \{0\}\,.
\end{equation}
Consequently  $$\int_M L(M)\cup r^* [c]={1\over C(\ell)} <\Ch
\Ind(\{D^{{\rm sign}}_\theta \}_{\theta\in T^k}),[\tau_c]>\,,\;\;
C(\ell)\not= 0\,.$$  Lusztig  settled the Novikov conjecture in
the abelian case by using the last formula and
 showing furthermore that the index bundle $\Ind\,\{D^{{\rm sign}}_\theta
\}$ is a {\it homotopy invariant}.

\subsection{Bismut's proof of the family index theorem.}$\;$

\m
  It is clear from what has been just explained that
Lusztig's treatment of the Novikov conjecture is heavily based on
the Atiyah-Singer {\it family index formula.} Besides the original
K-theoretic proof of Atiyah and Singer, see \cite{AS IV}, there is
 a {\it heat kernel proof} of the family index theorem, due to
Bismut \cite{Bismut}. Bismut's theorem applies to any family of
Dirac operators along the fibers of a fiber bundle $X\to B$;
notice that in the present case this fiber bundle is nothing but
$M\times T^k\to T^k$. We briefly explain Bismut's approach, as we
shall need it later.

\subsubsection{{\bf The superconnection heat-kernel.}}
Consider the bundle $\mathcal{E}$ on $T^k$ whose fiber
$\mathcal{E}_\theta$ at $\theta\in T^k$ is $C^\infty
(M,\Lambda_\CC^{{\rm sign}}(M)\otimes F_\theta)$. The Levi-Civita
connection on $M\times T^k$ and the connection $\nabla^F$ on the
vector bundle $F$ on $M\times T^k$ (see subsection
\ref{subsubsect:abelian-class}) define together a connection on
$\mathcal{E}$:
\begin{equation}\label{conn-abelian}
\nabla^\mathcal{E}:C^\infty (T^k, \mathcal{E})\longrightarrow
\Omega^1(T^k,\mathcal{E}):= C^\infty (T^k, \Lambda^1 (T^k)\otimes
\mathcal{E})\,.
\end{equation}
 The sum $$\AAA:=\{D^{{\rm
sign}}_\theta\}+\nabla^{\mathcal{E}}$$ is called a {\it
superconnection}; its curvature,
 $\AAA^2$, turns out to be a $T^k$-family of
differential operators on $M$ with coefficients in $\Omega^*
(T^k)$. \footnote{The concept of superconnection is due to 
D. Quillen who moreover suggested that it could be used 
in a heat kernel proof of the family index theorem. 
Quillen heuristic arguments were rigorously developed
by Bismut.} Thus $\exp( - \AAA^2)$ is a $T^k$-family $\{K(\theta)
\}_{\theta\in T^k}$ of {\it smoothing operators} on $M$ with
coefficients differential forms on $T^k$.

\subsubsection{{\bf The fiber-supertrace.}} Let $\Lambda^*_\theta
(T^k)$ the Grassmann algebra of the cotangent space to $T^k$ in
$\theta$. One can see more precisely that the Schwartz kernel of
$K (\theta)$ restricts to the diagonal $\Delta\leftrightarrow  M$
in $M\times  M$ as a section of the bundle $\Lambda^*_\theta (T^k)
\otimes {\rm End}\,( \Lambda_\CC^{{\rm sign}}(M)\otimes F_\theta
)$ over $M$. Let ${\rm str}_\theta$ denote the natural supertrace
on ${\rm End}\,( \Lambda_\CC^{{\rm sign}}(M)\otimes F_\theta )$;
we can extend this supertrace to $\Lambda^*_\theta (T^k) \otimes
{\rm End}\,( \Lambda_\CC^{{\rm sign}}(M)\otimes F_\theta )$ by
letting it act on the first factor as the identity. Thus $$ {\rm
str}_\theta ( K(\theta)|_{\Delta}) \in \Lambda^*_\theta
(T^k)\otimes C^\infty (M) \,.$$ We conclude that if $\{K (\theta)
\}_{\theta\in T^k}$ denotes the family of Schwartz kernels
associated to $\exp (-\AAA^2)$, then $$\int_M {\rm str}_\theta K
(\theta)|_{\Delta} \in \Lambda^*_\theta (T^k)$$ and  as $\theta$
varies in $T^k$ we obtain a differential form. {\it Summarizing}
we can give the following

\begin{definition}
 The functional analytic fiber-supertrace
$\STR(\exp( - \AAA^2))$ is the differential form on $T^k$ defined
by the equality $$\STR(\exp( - \AAA^2))(\theta)= \int_M {\rm
str}_{\theta} K(\theta)|_\Delta\,.$$
\end{definition}

\subsubsection{{\bf Bismut's theorem.}}
Consider the so-called rescaled superconnection
$\AAA_s:=s\{D^{{\rm sign}}_\theta\}+\nabla^{\mathcal{E}}$ for
$s>0$. Bismut's theorem, in this special case, states that
\begin{itemize}
\item  for each $s>0$
the differential form $\STR(\exp( - \AAA_s^2))$ is closed in $\Omega^* (T^k)$;
\item
for each $s>0$ it represents the Chern character of the index
bundle: $$\Ch (\Ind \,\{D^{{\rm sign}}_\theta\}\,)=[ \STR(\exp( -
\AAA_s^2))]\quad\text{in}\text \quad H^*_{{\rm dR}}(T^k)\,;$$
\item the
short-time limit  can be computed, giving
$$\lim_{s\downarrow 0}
\STR(\exp( - \AAA_s^2))=\int_M L(M,\nabla^g)\wedge \omega \,.$$
\end{itemize}

The notion of superconnection can be given for any  family of
Dirac operators $\{D_b\}_{b\in B}$ acting on the sections of  a
vertical Clifford module $E$ on a non-trivial fiber bundles
$Z\rightarrow M \xrightarrow{\pi} B$. \footnote{It is an operator
$\AAA:C^\infty (B,\mathcal{E})\longrightarrow C^\infty
(B,\Lambda^* (B)\otimes \mathcal{E})$, with $\mathcal{E}_b=
C^\infty (\pi^{-1} (b), E|_{\pi^{-1} (b)})$, which is odd with
respect to the total grading defined by $E$ and $\Lambda^* (B)$,
satisfies Leibnitz rule and can be written as $$\AAA=\{D_b\}
+\nabla^{\mathcal{E}}+\sum_{j=2}^k
\AAA_{[k]}\;\;\;\;\text{with}\;\;\;\;\AAA_{[k]}: C^\infty
(B,\mathcal{E})\longrightarrow C^\infty (B,\Lambda^k (B)\otimes
\mathcal{E})\,.$$ The first two results in Bismut's theorem are
true for any superconnection; the short-time limit, on the other
hand, only holds for a specific superconnection, nowadays called
the {\it Bismut's superconnection}; its rescaled version  can be
written as
\begin{equation}\label{bismut-super}
\AAA^{{\rm
Bismut}}_s=s\{D_b\}+\nabla^{\mathcal{E}}+\frac{1}{s}\AAA_{[2]}
\end{equation}
with $\AAA_{[2]}$ an additional term involving the curvature of
the fiber bundle $Z\rightarrow M \xrightarrow{\pi} B$. In
particular, if the fiber bundle is trivial, as in the Lustzig's
family, this additional term is zero.}

\n
Besides the original article of Bismut, \cite{Bismut}, the reader
is also referred  to \cite{BGV}, Chapter 9 and 10.

\medskip It is clear that if we wish to generalize Lustzig's
approach to a noncommutative group $\Gamma$ then  we need to bring
to the noncommutative context the notion of Chern character,
defined on $K_* (C^*_r \Gamma)$, and prove a {\it noncommutative
family index theorem}. In order to do so we need the notion of
cyclic (co)homology.

\subsection{Cyclic (co)homology.}\label{subsect:cyclic-algebraic}$\;$

\m

Let $A$ be a unital $k-$algebra over $k= \RR$ or $k=\CC$. The
cyclic cohomology groups $HC^*(A)$ Connes \cite{Connes-ihes} ( see also
Tsygan \cite{tsygan}) are the cohomology groups of the complex
$(C^n_\lambda, b)$ where $ C^n_\lambda$ denotes the space of
$(n+1)-$linear functionals $\varphi$ on $A$ satisfying the
condition: $$ \varphi(a^1, a^2, \ldots, a^n,a^0)= (-1)^n
\varphi(a^0, \ldots, a^{n+1})\,,\;\; \forall a^i \in A $$ and
where $b$ is the Hochschild coboundary map given by
\begin{align*} (b \varphi) (a^0, \ldots, a^{n+1})& = \sum_{j=0}^n
(-1)^j \varphi (a^0, \ldots, a^j a^{j+1}, \ldots, a^{n+1}) +\\&
(-1)^{n+1} \varphi(a^{n+1}a^0,\ldots , a^n). \end{align*} Set $
\overline{ C^0_\lambda}=C^0_\lambda$ and, for any $n\in \NN^*$,
denote by $\overline{ C^n_\lambda}$ the sub vector space  of $
C^n_\lambda$ formed by the $(n+1)-$linear functionals $ \varphi$
such that $ \varphi(a^0, a^1, \ldots, a^n)=0$ if $a^i=1$ for some
$i\in \{1,\ldots,n\}.$  $(\overline{ C^n_\lambda}, b)$ is then a
subcomplex of $ (C^n_\lambda, b)$ whose cohomology groups are
called the {\it reduced cyclic cohomology groups}
 $\overline{HC^*}(A).$
 \m

Of particular importance to us will be the cyclic cohomoly group
$HC^* (\CC\Gamma)$. Let $c \in H^k(\Gamma,\CC)$ be a group cocycle.
Connes has associated to c a cyclic cocycle
$\tau_c$ and thus a cyclic class $[\tau_c] \in HC^k(\CC\Gamma):$
if $\gamma_0,\ldots, \gamma_k \in \Gamma$ then set
\begin{align*} &\tau_c ( \gamma_0,\ldots, \gamma_k) = c ( 1_\Gamma,
\gamma_0,
 \ldots, \gamma_0 \cdots \gamma_{k-1}) \;\;\;
{\rm if}\, \gamma_0 \cdots \gamma_k=1_\Gamma, \\& \tau_c (
\gamma_0,\ldots, \gamma_k)=0  \;\;\;{\rm if}\, \gamma_0 \cdots
\gamma_k\not=1_\Gamma. \end{align*} If $k\geq 1$, then, using the
fact that $c$ is antisymmetric, one checks that $\tau_c$ is a
reduced cyclic cocycle.

\m

We can also introduce cyclic {\it homology.} Denote by
$A^{\otimes, n+1}$ the tensor product over $k$ of $n+1$ copies of
$A$ and consider the endomorphism $t$ of $A^{\otimes, n+1}$
defined by: $$ t( a_0\otimes a_1\otimes \cdots \otimes a_n) =
(-1)^n t( a_n\otimes a_0\otimes \cdots \otimes a_{n-1}). $$
Consider also the map $b: A^{\otimes, n+1}\rightarrow A^{\otimes,
n}$ defined by: \begin{align*} b( a_0\otimes a_1\otimes \cdots
\otimes a_n) & = \sum_{i=0}^{n-1} (-1)^i a_0 \otimes \cdots
\otimes a_i a_{i+1} \otimes \cdots \otimes a_n +\\ & (-1)^n a_n
a_0 \otimes a_1 \otimes \cdots \otimes a_{n-1}. \end{align*} Then
set
 $C^\lambda_n(A) = \frac{ A^{\otimes, n+1}}{{\rm Im}\, {\rm Id}-t}.$
 The  cyclic homology groups ${HC_*}(A)$ are then
the homology groups of the complex $( C^\lambda_n(A), b).$ Next,
denote by $\overline{C_n^\lambda} (A)$ the quotient of
$C^\lambda_n(A)$ by the sub $k-$module generated by the tensor
products $a_0\otimes a_1 \otimes \cdots \otimes a_n$ where $a_i=0$
for some $i\in \{1, \ldots, n\}.$ Then the reduced cyclic homology
groups $\overline{HC_*}(A)$ are defined to be the homology groups
of the complex $( \overline{C_n^\lambda}(A), b).$

\subsection{Noncommutative de Rham homology and the Chern
character.}\label{subsect:derham-algebraic}$\;$

\m
 We follow
Karoubi \cite{Karoubi}. Recall that $k$ is $\RR$ or $\CC$. Let $A$ be a
unital $k-$algebra and
 consider a
graded algebra
$$
\Omega_*(A)= \Omega_0(A)\oplus   \Omega_1(A) \oplus\Omega_2(A) \ldots
$$ with $\Omega_0(A)=A,$ endowed with a $k-$linear derivation of degree $1$,
$d=d_j:  \Omega_j(A)\rightarrow \Omega_{j+1}(A)$ satisfying
$d^2=0$ and
$$
d(\omega_j \cdot \omega_l) = d \omega_j \cdot \omega_l + (-1)^j
\omega_j \cdot d \omega_l, \;
\forall \omega_j \in  \Omega_j(A),\; \omega_l \in  \Omega_l(A).
$$ Denote by $[\Omega_*(A), \Omega_*(A)]_t$ the $k-$module
generated by the graded commutators $[\omega_j, \omega_l] =
\omega_j \cdot \omega_l - (-1)^{j l}\omega_l \cdot\omega_j$
where $j+l=t$ and $  \omega_j \in  \Omega_j(A),\; \omega_l \in  \Omega_l(A).$
We then set:
$$
\overline{\Omega_t(A)} = \frac{\Omega_t(A)} {[\Omega_*(A), \Omega_*(A)]_t}.
$$ It is clear that the derivation $d$ induces a $k-$linear differential,
still denoted $d$, on the graded $k-$vector space $ \overline{\Omega_t(A)}$,
we then denote by $\overline{H}_*(A)$ the homology of this quotient complex and call
 it the non commutative de Rham homology of $\Omega_*(A).$

\s
Now let $E$ be a finitely generated projective left $A-$module, a
connection $D$ on $E$ is a $k-$linear homomorphism $$ D: E
\rightarrow \Omega_1(A) \otimes_A E $$ satisfying Leibniz's rule
$$ \forall (a,s) \in A \times E,\; D(a\cdot s)= d a \otimes s + a
\otimes D(s). $$ Set $\Omega_{even}(A)= \oplus_{k\in \NN}
\Omega_{2 k}(A).$ One then checks that $D^2$ extends a left linear
$\Omega_{even}(A)-$endomorphism of $ \Omega_{even}(A) \otimes_A E$
sending  $ \Omega_{2k}(A) \otimes_A E$ into  $ \Omega_{2k+2}(A)
\otimes_A E$ for each $k\in \NN.$ Since $E$ is assumed to be
finitely generated and projective, there is a natural trace map:
$$ {\rm TR}:   \Omega_{even}(A) \otimes_A E\rightarrow
\frac{\Omega_{even}(A) }{[\Omega_{even}(A) , \Omega_{even}(A) ]}
\rightarrow \overline{\Omega_{even}}(A) $$ where the last
$\rightarrow$ is the obvious one. The Chern character is then
defined by $$ {\rm Ch}: K_0(A) \rightarrow \overline{H}_{even}(A)
$$ $$ E \rightarrow {\rm TR}\, e^{-D^2}. $$ It is indeed a theorem
(see Section 1 of \cite{Karoubi}) that ${\rm TR}\, e^{-D^2} $
defines a cycle in $ \overline{H}_{even}(A)$ which does not depend
on the choice of $D$.

\subsection{Cyclic (co)homology and noncommutative de Rham homology.}\label{subsect:derham-cyclic}$\;$

\m

We recall that Connes has constructed an operator $B$ from
  $\overline{HC_*}(A)$ (resp. $HC_*(A) $) to the Hochschild homology group
$H_{*+1}(A,A)$ where $B$ is a non commutative analogue of the de
Rham exterior derivative. In Section 2 of \cite{Karoubi} the
following is proved. For $*>0$, $\overline{H}_*(A)$ is isomorphic
to the kernel of $B$ acting on $\overline{HC_*}(A),$ while
 $\overline{H}_0(A)$ is isomorphic to the kernel
of $B$ acting on $HC_0(A)$. We shall not give the details here but
only retain the information that, for $*>0,$ there is a pairing
between noncommutative de Rham {\it homology} $\overline{H}_*(A)$
and the reduced cyclic {\it cohomology} group
$\overline{HC^*}(A).$ For $*=0$ there is  is a pairing between
noncommutative de Rham {\it homology} $\overline{H}_0(A)$ and
cyclic {\it cohomology} $HC^0 (A)$.

\subsection{Topological cyclic
(co)homology.}\label{subsect:derham-topological} $\;$

\m

Now let $A$ be a unital Fr\'echet locally  convex $k$-algebra;
i.e. a  Fr\'echet locally convex topological vector space for which
the product is continuous.  The {\it topological} cyclic
cohomology groups $HC^n(A)$ are defined as above but by
considering only continuous linear $(n+1)-$linear functionals.
Similarly, the topological cyclic homology groups $HC_n(A)$ are
defined as above but considering completed projective tensor
products. Moreover, one can define a completion
$\widehat{\Omega}_*(A)$ of ${\Omega}_*(A)$ which is a Fr\'echet
differential graded algebra. The noncommutative topological de
Rham homology $\widehat{H_*}(A)$ is defined as the homology of the
complex $$ \left( \widehat{\Omega}_*(A)/
\overline{[\widehat{\Omega}_*(A) , \widehat{\Omega}_*(A)]},\; d
\right)\,; $$ it pairs with the topological cyclic cohomology
$HC^*(A)$. In fact, if $*>0$, it pairs with the reduced
 topological cyclic cohomology.

\subsection{Smooth subalgebras of
$C^*$-algebras.}\label{subsect:smooth}$\;$

\m

In general the  {\it topological} cyclic homology of a
$C^*$-algebra is too poor. For instance on a smooth manifold $M$
$$HC^{2p}(C^0 (M) ) \simeq HC^{0}(C^0 (M) ) \;\;\;\text{and}\;\;\;
HC^{2p+1}(C^0 (M) )=0\;\; \forall p \in \NN \,.$$ In fact the
right algebra to consider in order to recover the (co)homology of
a smooth manifold $M$ is the algebra of {\it smooth functions} on
$M$, as there  are many more interesting cyclic cocycles on
$C^\infty (M)$ than on $C^0(M)$ \footnote{ For example, the
following interesting $2-$cyclic cocycle on $C^\infty (S^2)$, $
(a^0, a^1,a^2) \rightarrow \int_{S^2} a^0 d a^1 \wedge d a^2 $
does not extend to $C^0(S^2)$.}. In order to  further clarify this
point let us recall that Connes has defined a periodicity operator
$$ S: HC^k(A)\rightarrow HC^{k+2}(A), $$ and introduced the two
periodic cyclic cohomology groups $$P HC^{even} (A)= \lim_{+\infty
\leftarrow S} HC^{2k}(A)\,, \quad P HC^{odd} (A)= \lim_{+\infty
\leftarrow S} HC^{2k+1}(A).$$ The relationship between the
homology of $M$ and  cyclic cohomology 
is then the following: $$P HC^{even} (C^\infty(M))= \oplus_{k\in
\NN} H_{2k}(M ; \CC)\,,\quad\quad P HC^{odd} (C^\infty(M))=
\oplus_{k\in \NN} H_{2k+1}(M ; \CC).$$ Notice now that $C^\infty
(M)$ is a {\it dense subalgebra} of $C^0(M)$ which is furthermore
{\it closed under holomorphic functional calculus.} In general, if
$A$ is a $C^*$-algebra and $\B\subset A$ is a (Fr\'echet locally
convex) dense subalgebra closed under holomorphic functional
calculus, then $K_* (A)\simeq K_* (\B)$; such a subalgebra is
usually referred to as a {\it smooth subalgebra}. Thus, for
example,  $K_* ((C^\infty(M))\simeq K_* (C^0 (M)).$ So,
considering a {\it smooth subalgebra} $\B$ of a $C^*$-algebra $A$
allows us on the one hand to leave the $K$-theory unchanged
 and, on the other hand,   to consider an
interesting {\it topological} cyclic cohomology and thus, from the
previous subsection, an {\it interesting Chern character
homomorphism}: $$ {\rm Ch}\,: K_0 (\B) \rightarrow
\widehat{H_*}(\B). $$

\subsection{The smoothing of the index
class.}\label{subsect:smoothing} $\;$

\m
On the basis of our discussion so far, it is clear that in order
to apply an interesting  Chern character to our index class $\Ind
(D^{{\rm sign}}_{(M,r)})$, we need to {\it fix a subalgebra
$\mathcal{B}^\infty$ of $C^*_r \Gamma$ which is dense and closed
under holomorphic functional calculus}. As we have explained, it
is only by fixing such a subalgebra that we can hope to land, via
the Chern character, into
an interesting noncommutative de Rham homology.

Such an algebra does exist and it is called the Connes-Moscovici
algebra. Let us see the definition.
 Fix a word metric $\| \cdot \|$ on $\Gamma$.
Define an unbounded operator $D$ on $\ell^2(\Gamma)$ by setting
$D(e_\gamma)= \|\gamma\| e_\gamma $ where $(e_\gamma)_{\gamma \in
\Gamma}$ denotes the standard orthonormal basis of
$\ell^2(\Gamma)$. Then consider the unbounded derivation  $
\delta( T)= [D, T]$ on ${B} (\ell^2(\Gamma))$ and set
$$\mathcal{B}^\infty=\{ T \in C^*_r(\Gamma)/\; \forall k \in \NN,
\; \delta^k(T) \in {B}(\ell^2(\Gamma)) \}.$$ One can prove that
$\mathcal{B}^\infty$ is dense in $C^*_r \Gamma$ and closed under
holomorphic functional calculus. Thus $K_* (C^* _r \Gamma)\simeq
K_* (\mathcal{B}^\infty)$; the image of $\Ind (D^{{\rm
sign}}_{(M,r)})\in K_* (C^*_r \Gamma))$ in $K_*
(\mathcal{B}^\infty)$ under this isomorphism
 should be thought of as a "smoothing" of the index class,
 since in the commutative context it is nothing but the passage from a {\it continuous} index bundle
 for the Lustzig's family to a {\it smooth} index bundle.
Since $\Bi$ is a smooth subalgebra 
one may define $\widehat{\Omega}_\ast(\Bi)$ and
$\widehat{H}_*(\Bi)$ as above.

The smoothing of the index class can in fact be achieved directly
and explicitly. We wish to explain this point, for it will be
important in the next subsection. We do it directly for the
signature operator but it is clear that what we explain will hold
for any Dirac-type operator. Let $\Bi$, $\CC\Gamma\subset
\Bi\subset C^*_r \Gamma$, be {\it any} smooth subalgebra of $C^*_r
\Gamma$. Thus $\Bi$ is dense and holomorphically closed in $C^*_r
\Gamma$. Consider the flat $\Bi$-bundle $ \mathcal{V}^\infty = \Bi
\times_\Gamma \widetilde{M} \rightarrow M$ and set
$$\mathcal{E}^\infty:= \mathcal{V}^\infty \otimes_\CC
\Lambda^{{\rm sign}}_\CC (M)\,,\quad \mathcal{E}^{\infty,\pm}:=
\mathcal{V}^\infty \otimes_\CC \Lambda^{{\rm sign},\pm}_\CC (M)
\,.$$ Proceeding as in subsection \ref{subsect:index-class-closed}
we see that the signature operator on $M$ defines in a natural way
an odd $\Bi$-linear signature operator $$D^{{\rm
sign},\infty}_{(M,r)}\,:\, \CI(M,\mathcal{E}^\infty)\to
\CI(M,\mathcal{E}^\infty).$$ For simplicity, we keep the notation
$D^{{\rm sign}}_{(M,r)}$ for this operator. It is possible to
develop a $\Bi$-pseudodifferential calculus
$\Psi^*_{\Bi}(M,\mathcal{E}^\infty)$ and construct a parametrix
for $\D^{{\rm sign}}_{(M,r)}$ with rests
$\Psi^{-\infty}_{\Bi}(M,\mathcal{E}^\infty)$.  Starting from a
$\Bi$-parametrix one can prove a decomposition theorem analogous
to the one appearing in (\ref{Mishchenko_decomposition_closed});
thus
  \begin{equation}\label{Mishchenko_decomposition_closed bi}
    C^\infty(M, \mathcal{E}^{\infty,+}) = \mathcal{I}_+ (\infty) \oplus
    \mathcal{I}_+^{\perp} (\infty),\quad C^\infty(M,
    \mathcal{E}^{\infty,-})=\mathcal{I}_- (\infty) \oplus
    \D_{(M,r)}^+(\mathcal{I}_+^{\perp} (\infty))\,,
  \end{equation}
  with  $\mathcal{I}_+ (\infty)$ and $\mathcal{I}_- (\infty)$
   {\it finitely generated projective} $\Bi$-modules and $\D_{(M,r)}^{{\rm sign},+}$ inducing an isomorphism
   (in the Fr\'echet
  topology) between $\mathcal{I}_+^{\perp} (\infty)$ and
  $\D_{(M,r)}^+(\mathcal{I}_+^{\perp} (\infty))$.
The proof of this $\Bi$-decomposition theorem rests ultimately on
the fact that $\Bi$ is dense and closed under holomorphic
functional calculus in $C^*_r \Gamma$. For the proof see
Leichtnam-Piazza \cite[Appendix A] {LPMEMOIRS} and also Lott \cite[Section 6]{Lott III}.
Summarizing, the index class can be defined directly in $\Bi$:
\begin{equation}\label{bi index class closed}
\Ind (\D^+_{(M,r)})= [\mathcal{I}_+ (\infty)] - [\mathcal{I}_-
(\infty) ]\;\in\;K_0 (\Bi )\,.
\end{equation}

\subsection{The higher index theorem of Connes-Moscovici
(following Lott).}\label{subsect:higher-closed}$\;$

\m
  One can prove that the
heat operator associated to the Dirac laplacian on $\widetilde{M}$
defines a heat operator  $\exp (-(s\D^{{\rm sign}}_{(M,r)})^2)$
which is a $\Bi$-smoothing operator, i.e.  $\exp (-(s\D^{{\rm
sign}}_{(M,r)})^2)\in \Psi^{-\infty}_{\Bi}$. Inspired by Bismut's
heat-kernel proof of the family index theorem, Lott has defined in
\cite{Lott I} a certain noncommutative  connection on
$\mathcal{E}^\infty\equiv \mathcal{V}^\infty\otimes_\CC
\Lambda^{{\rm sign}}_\CC (M) $:
\begin{equation}\label{lott-connection}
 \nabla:
C^\infty(M,\mathcal{E}^\infty) \rightarrow C^\infty(M,
\widehat{\Omega}_1(\Bi)\otimes_{\Bi} \mathcal{E}^\infty).
\end{equation}
This is the analogue, in the nonabelian case, of
(\ref{conn-abelian}). He has then considered the rescaled
superconnection $\AAA_s:=s \D^{{\rm sign}}_{(M,r)}  + \nabla$ and,
using Duhamel expansion, the heat operator $$
e^{-\AAA_s^2}\,: C^\infty(M, \mathcal{E}^\infty ) \rightarrow
C^\infty(M, \widehat{\Omega}_*(\Bi)\otimes_{\Bi}
\mathcal{E}^\infty).$$ For any real $s>0$, this is, in a sense
that can be made precise, a $\Bi$-smoothing operator with
coefficients in $\widehat{\Omega}_* (\Bi) $. The restriction of
the superconnection heat kernel  $\mathcal{K}(e^{-\AAA^2_s})$ to
the diagonal $\Delta\leftrightarrow M$ in $M\times M$ is an element in
 $$ \widehat{\Omega}_* ({\cal {B}}^\infty) \otimes_{
{\cal {B}}^\infty} C^\infty(M,  {\cal{V}}^\infty \otimes_{\CC} {\rm End}\, E );
$$ taking the vector bundle supertrace ${\rm str}_E$ we get a
supertrace $$\STR ( e^{-\AAA^2_s})\,:= \int_M {\rm str}_E
\mathcal{K} (e^{-\AAA^2_s})|_{\Delta} d{\rm vol}_g\;\;\;\text{with
values
in}\;\;\;
\widehat{\Omega}(\Bi)/\overline
{[\widehat{\Omega}(\Bi),\widehat{\Omega}(\Bi)]}\,.$$ Notice that
since the algebra of non commutative differential forms
$\widehat{\Omega}_*(\Bi)  $ is not commutative, the super trace  $
{\rm STR}$ must take values in the quotient space $$
\widehat{\Omega}_*(\Bi)/ \overline{[ \widehat{\Omega}_*(\Bi),
\widehat{\Omega}_*(\Bi)]} $$ (i.e. modulo the closure of the space
of graded commutators; we take the closure so as to ensure that
the quotient space is Fr\'echet). Using Getzler's rescaling
\cite{Getzler-cmp} and adapting to the noncommutative context
Bismut's proof of the family index theorem, Lott proves in
\cite{Lott I} that
\begin{itemize}
\item the noncommutative differential form  ${\rm STR} (\,e^{-\AAA_s^2})$ is closed
\item its homology class is equal to the Chern character of the index:
 $${\rm Ch}\, \Ind ( D^{{\rm
sign}}_{(M,r)})= [{\rm STR}
\,e^{-\AAA_s^2}]\quad\text{in}\quad\widehat{H}_* (\Bi) .$$
\item there exists  a certain closed biform $\omega_{(M,r)}\in \Omega^*(M) \otimes
\widehat{\Omega}_* (\Bi)$ such that  $$\lim_{s\downarrow 0} {\rm
STR} \,e^{-\AAA_s^2}=\int_M L(M,\nabla^g) \wedge \omega_{(M,r)} $$
with the limit taking place in $\widehat{\Omega}(\Bi)/
\overline{[\widehat{\Omega}(\Bi),\widehat{\Omega}(\Bi)]}\,.$
\end{itemize}

In this way, we have explained how Lott has proved  {\it the
higher index theorem on Galois covering}:
\begin{equation}\label{lottformula} {\rm Ch}\, \Ind ( \D^{{\rm
sign}}_{(M,r)}) = [\int_M L(M,\nabla^g) \wedge \omega_{(M,r)}] \in
\widehat{H}_*(\Bi)
\end{equation}
In fact, one can prove that $\omega_{(M,r)}$ is an element in $
\Omega^*(M) \otimes \Omega_* (\CC\Gamma)$; however, we do point
out that the equality (\ref{lottformula}) only  makes sense in
$\Bi$.

\subsection{The Novikov conjecture
for hyperbolic groups.} $\;$

\m

Let $[c]\in H^\ell(B\Gamma,\CC)\equiv H^\ell (\Gamma,\CC)$ and let
$[\tau_c]\in HC^\ell (\CC\Gamma)$ the corresponding cyclic
cocycle. Lott has also proved \cite{Lott I}  that, in general,
there exists a nonzero constant $C(\ell)$ such that
$$\frac{1}{C(\ell)} < [\int_M L(M) \wedge \omega_{(M,r)}] ;
[\tau_c]
> = \int_M L(M) \wedge r^*(c)$$ where on the left-hand-side the
pairing between noncommutative de Rham {\it homology} and cyclic
{\it cohomology} has been used.

By formula (\ref{lottformula}), this means that if $[\tau_c]\in
HC^\ell (\CC\Gamma)$ {\it extends} to $HC^\ell (\Bi)$ then
\begin{equation}\label{cm-higher}
\frac{1}{C(\ell)} <{\rm Ch}\, \Ind ( \D^{{\rm
sign}}_{(M,r)}),[\tau_c]> = \frac{1}{C(\ell)} <[\int_M
L(M,\nabla^g) \wedge \omega_{(M,r)}],[\tau_c]>=\int_M L(M) \wedge
r^*(c)\,.
\end{equation}

The equality of the first and last term is due to Connes and
Moscovici and it is known as the {\it Connes-Moscovici higher
index theorem on Galois coverings}. We anticipate that the extra
information given by Lott's heat-kernel proof will be crucial on
manifolds with boundary.  Thus, for cyclic cocycles that {\it
extends from $HC^* (\CC\Gamma)$   to $HC^* (\Bi)$} we have
expressed the higher signatures in terms of the index class: $$
\int_M L(M) \wedge r^*(c)= \frac{1}{C(\ell)} <{\rm Ch}\, \Ind (
\D^{{\rm sign}}_{(M,r)}),[\tau_c]>\,.$$
 Since
the index class is a homotopy invariant, we conclude  that the
Novikov conjecture is established for all those groups having the
extension property for all the cocycles $\tau_c$. Connes and
Moscovici have shown that  Gromov hyperbolic groups do satisfy
this fundamental property; their proof exploits results by
Haagerup, de la Harpe and Jolissaint. We shall not give here the
definition of Gromov hyperbolic group but refer the reader instead
to Gromov \cite{gromov}, Ghys \cite{Ghys}, Connes-Moscovici
\cite{CM} and Connes \cite{Connes}.
 Basic examples of hyperbolic groups
are provided by  fundamental groups of a compact connected Riemann
surfaces of genus $g>1$ or more generally by  fundamental groups
of compact, negatively curved manifolds. Summarizing:

\begin{theorem} (Connes-Moscovici \cite{CM})
If $\Gamma$ is Gromov hyperbolic, then the Novikov conjecture is
true.
\end{theorem}

In fact, Connes and Moscovici even proved that
the Strong Novikov conjecture holds for Gromov hyperbolic groups.
 It should be remarked that there are now K-theoretic proofs of
 this
 result: after the work of Connes-Moscovici appeared,
 Ogle has  proved \cite{Ogle} by $K$-theoretic methods
that $\mu_\RR: K_* (B\Gamma)\otimes\RR\rightarrow
 K_* (C^*_r \Gamma)\otimes \RR$ is {\it injective} for Gromov
 hyperbolic groups.
In fact recently
Mineyev and Yu have proved that 
the Baum-Connes conjecture holds for Gromov-hyperbolic groups,
see \cite{Yu}.

\subsection{Groups having the extension property.}\label{subsect:extension}$\;$

\m
We can slightly generalize the content of the previous subsection
as follows. Let $\Gamma$ be a finitely generated group. We shall
say that $\Gamma$ {\it has the extension property} if there exists
a subalgebra $\Bi$ of $C^*_r \Gamma$, $\CC\Gamma \subset \Bi
\subset C^*_r \Gamma$, with the following 2 properties:
\begin{itemize}
\item $\Bi$ is dense and holomorphically closed in $C^*_r \Gamma$.
\item  Each class $[c]\in H^*(\Gamma; \CC)$ has
a cocycle representative whose corresponding cyclic cocycle
$\tau_c \in ZC^*(\CC \Gamma)$ extends to a continuous cyclic
cocycle on $\mathcal{B}^\infty$.
\end{itemize}
Examples of groups satisfying the extension property are Gromov
hyperbolic groups and also virtually nilpotent groups, see
\cite{J}. For this latter example it suffices to recall that by a
result of Gromov a group $\Gamma$ is virtually nilpotent if and
only if is of polynomial growth with respect to a (and thus any)
word metric ; the smooth subalgebra for such a group is simply
given by $$\Bi:=\{f:\Gamma\to \CC\,| \, \forall N\in
\NN\,\,\sup_{\gamma\in\Gamma} (1+\|\gamma\|)^N
|f(\gamma)|<\infty\}\,.$$

The following theorem, again due to Connes and Moscovici, is the
main result of this entire section
\ref{sect:connes-moscovici-lott}:

\begin{theorem}
If $\Gamma$ has the extension property, then the Strong Novikov
conjecture is true.
\end{theorem}

\section{{\bf The cut-and-paste problem for higher
signatures.}}\label{sect:cut-example}

Let $M$ and $N$ be two oriented compact  manifolds with boundary
and let $\phi,\psi: \partial M \to \partial N$ be orientation
preserving diffeomorphisms. We consider the closed oriented
manifolds  $$M\cup_{\phi} N^- \quad\text{and}\quad  M\cup_{\psi}
N^- \,;$$ we shall sometime use the notation $X_\phi:=
M\cup_{\phi} N^- $ and $ X_\psi:= M\cup_{\psi} N^- \,.$ Let
$r:M\cup_{\phi} N^-\to B\Gamma$, $s:M\cup_{\psi} N^- \to B\Gamma$
be reference maps and  assume that the two coverings are cut-and-paste equivalent,
see definition \ref{def:cut-and-paste}.

The cut-and-paste problem for higher signatures can be then stated
as follows: {\it for any $c\in H^*(B\Gamma, \QQ)$, compare the two
higher signatures: } $$ \int_{M\cup_{\phi} N^-} L(M\cup_{\phi}
N^-) \cup r^*(c),\quad \int_{ M\cup_{\psi} N^-} L(M\cup_{\psi}
N^-) \cup s^*(c). $$

The problem (raised by J. Lott and S. Weinberger, see \cite[Section 4.1]{Lott 3} and
\cite{Weinberger(1997)}) is then
to determine  which higher signatures of closed manifolds are cut
and paste invariant;  we refer to \cite[Section 4.1]{Lott 3} for
further discussion.

As remarked by Lott in \cite[Section 4.1]{Lott 3}, it is
implicitly established in \cite{Karras-Kreck-Neumann-Ossa (1973)}, \cite{Neumann}.
that, in general, higher signatures of closed manifolds are not
cut and paste invariant.  We shall
describe below a recent counterexample constructed in
\cite[Example 1.10]{LLK} to which we refer for the details.

\smallskip
\n
{\bf Example.} Let $s: \CC P^2 \times S^1 \rightarrow B\ZZ=S^1$ be
the reference map given by $s(z,e^{i\theta}) = e^{i\theta} .$ Then
there exists a compact oriented $4-$dimensional manifold $F$
endowed with an orientation preserving diffeomorphism $h$ such
that $( \CC P^2 \times S^1, s)$ is cobordant to $M( (F,h), T)$
where $M(F,h)$ denotes the mapping torus obtained from
$[0,1]\times F$ by identifying $(0,x)$ with $(1,h(x)$.  It is
shown by M. Kreck in \cite{LLK} that $F$ may be choosen of the
form $ (S^1 \times S^3) \# ( \CC P^2 \times \overline{ \CC P }^2)
\# m (S^2 \times S^2)$ for a suitable $m\in \NN$. The reference
map $T: M(F,h) \rightarrow B\ZZ$ induces a map $r: F\rightarrow
B\ZZ$ such that $r$ and $r\circ h$ are homotopic as (continuous)
maps from $F$ to $B\ZZ$. M. Kreck has shown that one may assume
that $r:F\rightarrow B\ZZ$ is two-connected. Moreover there exists
a manifold with boundary $W$ such that
 $\partial W=F$ and there are two maps $R, R^\prime$ from $W$
to $B\ZZ$ such that $r=R_{|\partial W}$ and $r\circ h=
R^\prime_{\partial W}$. Therefore, $(M(F,h), T)$ (and thus $( \CC
P^2 \times S^1, s)$) is cobordant to: $$ (W\cup_{id} W, R\cup
R)\,-\,(W\cup_h W, R\cup R^\prime ). $$ Thus, $( \CC P^2 \times
S^1 \times S^1, s\times \id_{S^1} )$ is cobordant to $$
(\,(W\cup_{id} W)\times S^1, (R\cup R)\times \id_{S^1} \,) -
(\,(W\cup_h W)\times S^1, (R\cup R^\prime) \times \id_{S ^1} \,)
$$ where $ s \times \id_{S^1}: \CC P^2 \times S^1\times S^1
\rightarrow B\ZZ \times B\ZZ$.

Now, let $\omega_1$ denote the fundamental class of $S^1$. Then,
since the signature of of $\CC P^2$ is not zero, one checks
immediately that $$ \int_{\CC P^2 \times S^1 \times S^1} L( \CC
P^2 \times S^1 \times S^1)\wedge (s \times
\id_{S^1})^*(\omega_1\times \omega_1) \not=0. $$
 Then, by cobordism invariance, it is clear that $(\,(W\cup_{id} W)\times S^1, (R\cup
R)\times \id_{S^1}\,)$ and $(\,(W\cup_h W)\times S^1, (R\cup
R^\prime)\times \id_{S^1}\,)$ do not have the same higher
signatures. {\bf End of example.}

\medskip Despite the negative result explained in the previous
example, we can ask whether we can give sufficient conditions (on
$\Gamma$ and on the two coverings defined by $r|_{\pa M}$ and
$s|_{\pa M}$) ensuring that the higher signatures are indeed
cut-and-paste invariant. This would answer, at least partially,
Question 2 in subsection \ref{sect:3questions}. Now, for the lower
signature $\int_M L(M)$ we have explained 3 different ways for
treating the cut-and-paste problem; the first method makes use of
the Atiyah-Patodi-Singer index formula for the signature of a
manifold with boundary , the second method employs a purely
topological argument, whereas the third method uses a spectral
flow argument (based, ultimately, on a gluing formula for the
index and a variational formula for the Atiyah-Patodi-Singer
index).

As we shall now see,  these 3 methods can be pursued in the higher
case too. We shall begin by the first method and in fact explain a
general theory of {\it higher signatures on manifolds with
boundary}, thereby answering simultaneously to Question 2 and
Question 3 of section \ref{sect:3questions}.

\section{{\bf Higher signatures on manifolds with
boundary.}}\label{sect:higher-boundary}

\subsection{Introduction and main strategy for the definition.} 

\subsubsection{{\bf Introduction.}}
Let $M$ be an oriented manifold {\it with boundary} and let
$r:M\to B\Gamma$ be a classifying map. Let $[c]\in H^*
(B\Gamma,\RR)$. Since the expression $ \int_M L(M, \nabla^g)\cup r^*[c]$ 
depends on the choice of the metric $g$, 
it is not obvious how to define  the higher
signature ${\rm sign}(M,r;[c])$ associated to $r$ and $[c]\in H^*
(B\Gamma,\RR)$. Still, Theorem \ref{prop:hominv} shows that
the difference $$\int_M L(M,\nabla^g)-{1\over 2} \eta(D^{{\rm
sign}}_{(\partial M,g_\pa)})$$ is an oriented homotopy invariant
of the pair $(M,\partial M)$; in other words by {\it subtracting}
a suitable {\it boundary correction term} to the metric-dependent
integral $\int_M L(M,\nabla^g)$ we have produced a homotopy
invariant of the pair $(M,\pa M)$.
In the higher case, having observed that  on a manifold with boundary $M$ the higher
analogue of $\int_M L(M,\nabla^g)$ is the metric-dependent
integral \begin{equation}\label{integral-omega} \int_M
L(M,\nabla^g)\wedge \omega_{(M,r)}\in \widehat{\Omega}( \Bi)
\end{equation}
 appearing
in (\ref{cm-higher}), we ask ourselves the following

\m
\n
{\bf Fundamental question:} {\it which boundary correction term
should we subtract to (\ref{integral-omega}) in order to
obtain a {\it homotopy invariant} noncommutative de Rham class in
$\widehat{H}_* (\Bi)$} ?

\m
To have a feeling on the strategy we shall follow, let us 
recall 
how Lustzig
managed to prove the Novikov conjecture for $\Gamma=\ZZ^k$ in the
closed case. The proof was in four steps:

\m\n
(i) Define a suitable family of twisted signature operators 
$\{D^{{\rm sign}}_\theta\}_{\theta\in T^k}$, $T^k:= {\rm Hom}(\ZZ^k, U(1))$,
and its index class in $K^0 (T^k)$.

\medskip
\n
(ii) Prove the homotopy invariance of the index class $\Ind
(\{D^{{\rm sign}}_\theta\}_{\theta\in T^k})\in K^0 (T^k)$.

\medskip
\n
(iii) Apply the family index formula, thus computing
the Chern character of the index class as $[\int_M L(M)\wedge
\omega]$$\in$ $H^*(T^k,\RR)$, $\omega\in \Omega^*(M\times T^k)$.

\medskip
\n
(iv) Express the higher signatures in terms of the pairing
between this cohomology class and a homology class $[\tau_c]\in
H_* (T^k,\RR)$ naturally defined by $[c]\in H^*(B\ZZ^k,\RR)\equiv
H^* ((T^k)^*,\RR)$.

\subsubsection{{\bf Strategy in the commutative case.}}
Let now $M$ have a boundary, $\pa M\not= \emptyset$. We  assume
again  $\pi_1 (M)=\ZZ^k$. The Lustzig's family $\{D^{{\rm
sign}}_\theta\}_{\theta\in T^k}$ is still perfectly defined and is
a family of Dirac-type operators on the manifold with boundary
$M$. Keeping in mind Lustzig's approach and our discussion in the
case of a single manifold (Theorem \ref{prop:hominv}), we
would like to

\medskip
\n
(i) define a Atiyah-Patodi-Singer ($\equiv$ APS) {\it index class}, in $K^*
(T^k)$, for the Lustzig's family.

\medskip
\n
(ii) establish the homotopy invariance of this index class.

\medskip
\n
(iii) prove a {\it family index formula} for its Chern character
in $H^* (T^k,\RR)$; this formula will involve the {\it boundary
correction term} we alluded to in the fundamental question raised
above

\medskip
\n
(iv) {\it define} the higher signatures by coupling the Chern
character with $\tau_c\in H_* (T^k,\RR)$.

\subsubsection{{\bf Strategy in the noncommutative case.}}
Let us pass to the noncommutative case and consider $(M,r:M\to
B\Gamma)$. Keeping in mind the analogy between higher index theory
and family index theory, we would like to

\medskip
\n
(i)  define a APS {\it index class} associated to
$D^{{\rm sign}}_{(M,r)}$; this class will live in $K_* (\Bi)=K_*
(C^*_ \Gamma)$.

\medskip
\n
(ii) establish the homotopy invariance of the this index class.

\medskip
\n
(iii) prove a {\it higher index formula} for its Chern character in
$\widehat{H}_* (\Bi)$; this formula will have to involve the {\it
boundary correction term} we alluded to in the fundamental question.

\medskip
\n
(iv) {\it define} the higher signatures ${\rm sign}(M,r;[c])$ {\it
for a group satisfying the extension property}, by coupling the
Chern character in  $\widehat{H}_* (\Bi)$ with the {\it extended}
cyclic cocycle $\tau_c\in HC^* (\Bi)$ defined by $[c]\in H^*
(\Gamma,\CC)\equiv H^* (B\Gamma,\CC)$.

\m
The details of this program, which was conceived by Lott in
\cite{Lott II}, shall now be explained. We begin once again by the
commutative case.

\subsection{The Bismut-Cheeger eta form.}\label{subsect:bismut-cheeger} $\;$

\m

Let $M$ be an even dimensional oriented manifold with boundary
with $\pi_1(M)=\ZZ^k$, as in the previous subsection. Consider an
odd Dirac-type operator $D:C^\infty (M,E)\to C^\infty (M,E)$
acting on the sections of a $\ZZ_2$-graded Clifford bundle. For
each $\theta \in T^k$, one has a twisted  operator $D_\theta$
acting on $C^\infty ( M ; E\otimes F_\theta)$ where $F_\theta$ is
the flat complex line bundle of $M$ associated with $\theta\in
T^k:= {\rm Hom}(\ZZ^k, U(1))$ (see Section 7.3.3). Let us consider
the family $\D:=\{D_\theta\}_{\theta\in T^k}$ on $M$ parametrized
by the torus
 $T^k$ . From the
variational formula for the Atiyah-Patodi-Singer index, see
(\ref{variational}), one realizes immediately that the family of
Atiyah-Patodi-Singer boundary value problems associated to the
family of Dirac-type operators $\D:=\{D_\theta\}_{\theta \in T^k}$
{\it is not continuous} in $\theta \in T^k$, unless the boundary
family $\D_{\pa}:=\{D_{\theta,\pa M}\}_{\theta \in T^k}$ is
invertible (notice that in the latter case there would not be any
spectral flow).
Under this additional assumption Bismut and Cheeger defined an
index class $\Ind(\D,\Pi_{\geq})\in K^0 (T^k)$
 and proved a
family index formula for its Chern character in $H^{{\rm even}}
(T^k,\RR)$:
\begin{equation}\label{bismut-cheeger}
\Ch (\Ind (\D,\Pi_{\geq}))= [\int_{M}  \widehat{A}(M,\nabla^g){\rm
Ch}\,^\prime (E,\nabla^E)\wedge \omega -\frac{1}{2}
\widetilde{\eta} (\D_{\pa})]\in H^{{\rm even}} (T^k,\RR) \,.
\end{equation}
In this formula
 $\omega
$ is the bi-form in $\Omega^*(M\times T^k)$ we met in subsection
\ref{subsect:lusztig}, whereas $\widetilde{\eta} (\D_{\pa})\in
\Omega^* (T^k)$ is the Bismut-Cheeger {\it eta form} associated to
$\D_{\pa}$. This is our {\it boundary correction term}. The eta
form is defined as
\begin{equation}\label{etaform}
\widetilde{\eta} (\D_{\pa})={2\over\sqrt{\pi}}\int_0^\infty
 {\rm STR}_{{\rm Cl}(1)}\left( {d\BBB_s \over ds} e^{-\BBB_s^2} \right)ds\;\;\in\;\;
\Omega^{{\rm even}} (T^k)
\end{equation}
with $\BB_s$ the superconnection induced on the boundary by the
rescaled Bismut superconnection $\AAA_s$. 
The supertrace appearing in this formula is the {\it odd}
fiber-supertrace on the odd-dimensional boundary; it is
defined using the isomorphism $E|_{\pa M}\simeq E^+|_{\pa M}\otimes {\rm Cl}(1)$,
with ${\rm Cl}(1)$ denoting the complex Clifford algebra generated by 1 and $\sigma$.\\ As an example, the
0-degree part of this differential form, computed at $\theta \in
T^k$, is simply the eta invariant of $D_{\theta,\pa}$: $$
\widetilde{\eta} (\D_{\pa})_{[0]} (\theta )=\eta
(D_{\theta,\pa})\,.$$  
Notice that the operator ${d\BBB_s \over ds} e^{-\BBB_s^2}$ 
is again a smoothing operator with differential form coefficients.
The convergence of the
$s$-integral in (\ref{etaform})  near zero is non-trivial and
requires  Bismut's local index theorem for families.
 The convergence at $\infty$ depends heavily on the assumption
 that the family is invertible. The
Bismut-Cheeger eta form can be defined for any {\it invertible}
family of Dirac-type operators, not necessarily arising as a
boundary family. It is more generally defined for any {\it
invertible} family $\{D_b\}_{b\in B}$ acting on the sections of a
vertical Clifford bundle on a  fiber bundle $Z\to X\to B$ with odd
dimensional fiber.

\subsection{Lott's higher eta invariant in the invertible
case.}$\;$

\m

We now pass to the noncommutative case.  Let $(N,r:N\to B\Gamma)$
be closed and  odd dimensional (for example the boundary of an
even dimensional manifold with boundary). We fix a Riemannian
metric $g$ on $N$ and  consider a Dirac-type operator $D$ on $N$
acting between the sections of an ungraded Clifford module $E$. We
consider $E\otimes {\rm
Cl}(1)\simeq E\oplus E$. Let $\D_{(N,r)}:C^\infty (N,E\otimes
\mathcal{V}) \to C^\infty (N,E\otimes \mathcal{V})$ be the
associated $C^*_r (\Gamma)$-linear operator, with $\mathcal{V}=
 C^*_r \Gamma \times_\Gamma r^*E\Gamma $.
Fix now a {\it smooth subalgebra} $\Bi\subset C^*_r \Gamma$; for
example the Connes-Moscovici algebra.
 We still denote by
$\D_{(N,r)}$ the operator acting on $C^\infty (N,E\otimes
\mathcal{V}^\infty)$, with
$\mathcal{V}^\infty=\B^\infty\times_\Gamma r^* E\Gamma $. Let
$s\sigma\D_{(N,r)}+\nabla$ be the rescaled Lott superconnection.
The Schwartz kernel $\mathcal{K}(t)$ of the operator  $ \sigma
\D_{(N,r)} \exp( -(\nabla + t \sigma \D_{(N,r)} )^2 $, which is a
smoothing operator with coefficients in $\widehat{\Omega}_*(\Bi)$,
can be restricted to the diagonal in $N\times N$, giving $\mathcal{K}(t)|_{\Delta}$,
an element
in
$$\widehat{\Omega}_*(\Bi)\otimes_{\Bi} C^\infty(N,\mathcal{V}^\infty
\otimes_\CC {\rm End} (E\otimes {\rm Cl}(1)))$$ where we identify
$N\leftrightarrow \Delta$. As in the previous section
there is an odd-supertrace ${\rm
Str}_{{\rm Cl (1)}}$ acting on the endomorphisms of $E\otimes {\rm Cl}(1)$; 
using this vector-bundle odd supertrace
we can define the odd supertrace $ {\rm STR}_{{\rm Cl(1)}}$ of the
smoothing operator $ \sigma \D_{(N,r)} \exp( -(\nabla + t \sigma \D_{(N,r)}
)^2 $; this is the noncommutative differential form defined by
$${\rm STR}_{{\rm Cl(1)}}[ \sigma \D_{(N,r)} \exp( -(\nabla + t
\sigma \D_{(N,r)} )^2 )] := \int_N {\rm Str}_{{\rm Cl (1)}}
\mathcal{K}(t)|_{\Delta}\,d{\rm vol}_g
$$ Once again, since $\widehat{\Omega}_*(\Bi)  $ is not
commutative, the odd super trace  $ {\rm STR}_{{\rm Cl}(1)}$ must
take values in the quotient space $$\widehat{\Omega}_*(\Bi)/
\overline{[ \widehat{\Omega}_*(\Bi), \widehat{\Omega}_*(\Bi)]}$$
(i.e. modulo the closure of the space of graded commutators).
Summarizing, for each $t\in (0,\infty)$ we can consider the
following noncommutative differential form
\begin{equation}\label{higher-integrand}
\widetilde{\eta}(\D_{(N,r)})(t):= \frac{ 2} { \sqrt{\pi}}\,
 {\rm STR}_{{\rm Cl}(1)}\,
[ \sigma \D_{(N,r)} \exp( -(\nabla + t \sigma \D_{(N,r)} )^2 ) ]
\end{equation}
The following theorem is due to J. Lott:

\begin{theorem}\label{theo:convergence-eta}
Assume that  $\D_{(N,r)}$ is invertible in the Mishchenko-Fomenko
calculus. Then
\begin{equation}\label{higher-eta}
\widetilde{\eta}(\D_{(N,r)})=\int_0^\infty
\widetilde{\eta}(\D_{(N,r)})(t)dt
\end{equation}
converges in $$\widehat{\Omega}_*(\Bi)/
\overline{[\widehat{\Omega}_*(\Bi), \widehat{\Omega}_*(\Bi)]}. $$
 The
resulting form is called the {\it higher eta invariant of}
$\D_{(N,r)}$.
\end{theorem}

\n
{\bf Remarks.}

\s
\n
(i) The invertibility of $\D_{(N,r)}$ in the Mishchenko-Fomenko
calculus is equivalent to the existence of a full gap at
$\lambda=0$ in the $L^2$-spectrum of the operator $\widetilde{D}$
on $\widetilde{N}$, i.e. to the $L^2$-invertibility of
$\widetilde{D}$.

\smallskip
\n
(ii) Theorem \ref{theo:convergence-eta} is proved by Lott in
\cite{Lott II} for virtually nilpotent groups and implicitly in
\cite{Lott III} in the general case. For additional details on the
general case see also \cite[Theorem 4.1]{LPBSMF}.

\smallskip
\n
(iii) The higher eta invariant of Lott is the noncommutative
analogue of the Bismut-Cheeger eta form; the convergence of the
integral near $t=0$ follows from the local index theory developed
by Lott, in the same way that the convergence of the eta-form for
families is due to Bismut's local index theory. On the other hand,
the convergence for $t\rightarrow +\infty$ is much more delicate.
Once again, {\it the proof depends heavily on the invertibility of
$\D_{(N,r)}$}.

\subsection{Higher Atiyah-Patodi-Singer index theory in the invertible case.}
$\;$

\m

Let $(M,g)$ be a compact even-dimensional Riemannian manifold with
boundary. We assume $g$ to be of product type near $\pa M$ and we
let $D$ be a generalized Dirac operator acting on the sections of
a $\ZZ_2-$graded Hermitian Clifford module $E$. As a fundamental
example we could consider the signature operator $D^{{\rm sign}}$.
Let $\Gamma$ be a finitely generated discrete group and let
$\Bi\subset C^*_r \Gamma$ be a smooth subalgebra. Let $r:
M\rightarrow B \Gamma$ be a continuous map defining a
$\Gamma-$covering $\widetilde{M} \rightarrow M.$  We denote by
$\widetilde{D}$ the lift of $D$ to $\widetilde{M}$. We denote by
$\D_{(M,r)}: C^\infty (M,E\otimes \mathcal{V}^\infty)\to C^\infty
(M,E\otimes \mathcal{V}^\infty)$ the $\Bi-$left linear operator
induced by $\widetilde{D}$. The boundary operator associated to
$D$ will be denoted, as usual, by $D_{\pa}$. Making use of
$\widetilde{D}_{\pa}$ we also get an operator $\D_{(\partial
M,r|_{\pa M})}$ which is nothing but the boundary operator of
$\D_{(M,r)}$. We set $r|_{\pa M}:=r_\pa$. Assume now that
$\D_{(\partial M,r_\pa)}$ is invertible in the Mishchenko-Fomenko
calculus; equivalently, we assume that $\widetilde{D}_{\pa}$ is
$L^2-$invertible. Let $$\Pi_{\geq}=\ha \left( 1+ {\D_{(\pa
M,r_\pa)}\over |\D_{(\pa M,r_\pa)}|} \right)\,;$$ this is a 0th
order $\Bi$-pseudodifferential operator and we can consider the
domain $$C^\infty (M,E^+\otimes
\mathcal{V}^\infty,\Pi_{\geq})=\{s\in C^\infty ( M,E^+\otimes
\mathcal{V}^\infty)\,\,|\,\, \Pi_{\geq} (s|_{\pa M})=0\}\,.$$
 The following Theorem is conjectured in
 \cite{Lott II} and proved in   \cite{LPMEMOIRS}, \cite[Appendix]{LPBSMF}.

\begin{theorem}\label{theo:memoirs} Assume that $\D_{(\partial M,r_\pa)}$ is invertible in
the Mishchenko-Fomenko calculus. Then the operator $\D_{(M,r)}$
with domain  $C^\infty (M,E^+\otimes
\mathcal{V}^\infty,\Pi_{\geq})$ gives rise to a well defined
APS-index class ${\rm Ind}\, (\D_{(M,r)},\Pi_{\geq})$ in $K_0(\Bi)
\simeq K_0(C^*_r(\Gamma))$. The following formula holds in the non
commutative topological de Rham homology of $\Bi$:
\begin{equation}\label{hif-aps} {\rm Ch}\, {\rm Ind}\,
(\D_{(M,r)},\Pi_{\geq}) = \left[ \int_M {\rm AS}\wedge \omega
-\frac{1}{2}\widetilde{\eta}(\D_{(\partial M,r_\pa)})\right] \in
\widehat{H}_*(\Bi) \end{equation} with ${\rm AS}=
\widehat{A}(M,\nabla^g) \wedge {\rm Ch}^\prime (E,\nabla^E).$
\end{theorem}

\n
In particular: {\it under the invertibility assumption we have
proved that Lott's higher eta invariant is the boundary correction
term we have been looking for.}

\s
The proof of the theorem rests ultimately on the heat-kernel proof
of the higher index theorem given by Lott and on an extension to
Galois coverings of Melrose's $b$-pseudodifferential calculus on
manifolds with boundary. For the latter, the reader is referred to
the book by Melrose \cite{Melrose} and also to the surveys
\cite{piazza2} Grieser \cite{grieser}.

\subsection{Higher signatures on manifolds with
$L^2$-invisible boundary.} $\;$

\m
Let $(M,g)$ be a Riemannian manifold with boundary; we assume the
metric to be of product type near the boundary. Let
$\widetilde{M}\to M$ be a Galois $\Gamma$-covering; let $r:M\to
B\Gamma$ be a classifying map. 
We shall assume that the operator $\D^{{\rm sign}}_{(\pa
M,r_{\pa})}$ is invertible in the Mishchenko-Fomenko calculus.
Equivalently, the operator $\widetilde{D}^{{\rm sign}}_{\pa}$ is
$L^2$-invertible, or, again equivalently, the differential-form
Laplacian $\Delta_{\pa\widetilde{M}}$ is $L^2$-invertible in each
degree. {\it We shall say that the boundary $\pa \widetilde{M}$ is
$L^2$-invisible.} Recent results of Farber and Weinberger show
that there do exist coverings having a $L^2$-invisible boundary,
see \cite{Farber-Weinberger}. See also the subsequent paper
\cite{HRS} from which the term $L^2$-invisible is borrowed. Since
$(\pa M,r_{\pa})$ is $L^2$-invisible, the higher eta invariant of
Lott, $\widetilde{\eta}(\D^{{\rm sign}}_{(\pa N,r_\pa)})$, is well
defined. We set
\begin{equation}\label{notation-eta}
\widetilde{\eta}(\D^{{\rm sign}}_{(\pa
N,r_\pa)})\,:=\, \widetilde{\eta}_{(\pa N,r_\pa)}
\end{equation}

\begin{definition}\label{def:higher-class}
We define the higher signature class in $\widehat{H}_* (\Bi)$
of a covering $(M,r:M\to B\Gamma)$ with $L^2$-invisible
boundary as
\begin{equation}\label{higher-class}
\widehat{\sigma}(M,r) \,:=\, \left[ \int_M L(M,\nabla^g)\wedge
\omega_{(M,r)}- \ha\widetilde{\eta}_{(\pa N,r_\pa)} \right]\,\in\,
\widehat{H}_* (\Bi)\,.
\end{equation}
Notice that in this formula $\Bi$ is any {\it smooth} subalgebra
of $C^*_r \Gamma$, for example the Connes-Moscovici algebra.
\end{definition}

Let now $N$ be a manifold with boundary and let $(N,s:N\to
B\Gamma)$ be a Galois covering. Let $h:N\to M$, with $h(\pa
N)\subset \pa M$, a homotopy equivalence between $(N,s:N\to
B\Gamma)$ and $(M,r:M\to B\Gamma)$. A fundamental result of Gromov
and Shubin \cite{gromov-shubin} states that $(\pa N, s|_{\pa N}:
\pa N\to B\Gamma)$ is then also $L^2$-invisible. The following
result is conjectured in Lott \cite{Lott II} and proved in
Leichtnam-Piazza \cite{LPBSMF}:

\begin{theorem}
Let $M$ be an oriented manifold with boundary, let $r:M\to
B\Gamma$ be a classifying map and  assume that $(\pa M,r_{\pa}:
\pa M\to B\Gamma)$ is $L^2$-invisible. Then the higher signature
class $\widehat{\sigma}(M,r)$ is a oriented homotopy invariant of
the pair $(M,\pa M)$ and of the map $r:M\to B\Gamma$.
\end{theorem}

\begin{proof}
Following techniques of Kaminker-Miller \cite{KM}, one proves that
the APS-index class introduced in Theorem \ref{theo:memoirs},
$\Ind\, (\D^{{\rm sign}}_{(M,r)},\Pi_{\geq})\in K_0 (\Bi)$, is a
homotopy invariant. The Theorem follows at once from the higher
index formula (\ref{hif-aps}) applied to the signature operator.
\end{proof}

\begin{definition}\label{def:higher-invisible}
Let $[c]\in H^\ell (\Gamma,\CC)$.
 Let $\Gamma$ {\it have the extension property}, see subsection
 \ref{subsect:extension},
  and let $\tau_c\in ZC^\ell (\Bi)$
 be the extended cyclic cocycle associated to $c$. We  {\it
define} higher signatures ${\rm sign} (M,r,[c])\in \CC$ on a
manifold with $L^2$-invisible boundary by setting
\begin{equation}\label{higher-invisible}
 {\rm sign} (M,r,[c])\,:=\,<\widehat{\sigma} (M,r),[\tau_c ]>\,\equiv\,<[\int_M
L(M,\nabla^g)\wedge \omega_{(M,r)}- \ha\widetilde{\eta}_{(\pa
N,r_\pa)}],[\tau_c ]>\,.
\end{equation}
\end{definition}
If the boundary is empty then, up to the constant $C(\ell)$
appearing in (\ref{cm-higher}), we reobtain the Novikov higher
signatures.

\begin{corollary}\label{coro:hominv-invertible}
Let $\Gamma$ be a finitely generated discrete group having  the
extension property. On manifolds $(M,r:M\to B\Gamma)$ with
$L^2$-invisible boundary 
the higher signatures (\ref{higher-invisible}) are oriented
homotopy invariants for each $[c]\in H^* (\Gamma,\CC)$.
\end{corollary}

\n
The result hold more generally for certain {\it twisted} higher
signatures, manufactured out of the index class of twisted
signature operators. See \cite{LPBSMF}.


\s
\n
{\bf Remark.} Corollary \ref{coro:hominv-invertible} should be
seen as a topological application of the {\it higher}
Atiyah-Patodi-Singer index theorem \ref{theo:memoirs} . For
applications in the realm of  positive scalar curvature metrics
see Leichtnam-Piazza \cite{LPPSC}.

\subsection{Non-invertibility, perturbations and index
classes.}\label{subsect:perturbations}$\;$

\m

Let $(M,r:M\to B\Gamma)$ be a covering with non-empty boundary. We
set, as usual, $\widetilde{M}=r^* E\Gamma$, $r_\pa:=r|_{\pa M}$.
The invertibility assumption on $\D^{{\rm sign}}_{(\pa M,r_\pa)}$,
or, equivalently, the $L^2$-invertibility assumption on
$\Delta_{\pa\widetilde{M}}$, is very strong. In fact,  until the
recent work of Farber-Weinberger \cite{Farber-Weinberger}, it was
an open question whether for  a Galois covering $\Gamma\to
\widetilde{N}\to N$ it is always the case that the operator
$\Delta_{\widetilde{N}}$ is  not $L^2$-invertible (see
\cite{Lott-zero}).
 \\ For example, when $\Gamma=\ZZ^k$ the
invertibility condition requires  the cohomology groups of $\pa M$
with coefficients in the flat bundle $F_\theta$ to vanish {\it for
all} $\theta\in T^k$. Although this is indeed a strong hypothesis,
there is no way to avoid it if one wants to set up a {\it
continuous} family of APS-boundary value problems for the
Lustzig's family or if one wants to prove the large time
convergence of the integral defining the eta form. Similarly, in
the noncommutative context, we do need the invertibility of
$\D^{{\rm sign}}_{(\pa M,r_\pa)}$ for the projection
$$\Pi_{\geq}=\ha\left( 1+\frac{\D^{{\rm sign}}_{(\pa
M,r_\pa)}}{|\D^{{\rm sign}}_{(\pa M,r_\pa}|} \right)$$ to make
sense as $C^*_r \Gamma$-linear operator \footnote{Notice that in
the context of $C^*$-algebras Hilbert modules we only have a {\it
continuous} functional calculus; in particular the operator
$\chi_{[0,\infty)} (\D^{{\rm sign}}_{(\pa M,r_\pa)})$ does not
make sense  as $\Bi$-linear or $C^*_r \Gamma$-linear operator. It
is only by going to a Von Neumann context that one can make sense
of the operator $\chi_{[0,\infty)} (\D^{{\rm sign}}_{(\pa
M,r_\pa)})$.}. 
The
invertibility is also necessary in order to prove the convergence
of the higher eta invariant. The question then arises as whether
it is possible to lift the invertibility assumption on the
boundary operator and still develop a family index theory or a
higher index theory on manifolds with boundary. This problem was
tackled for the first time  by Melrose and Piazza in \cite{MP I}
\cite{MP II} and subsequently extended to the noncommutative
context in Leichtnam-Piazza \cite{LPGAFA}, \cite{LPCUT}. We shall now
give a very short account of this theory, concentrating on the
results leading to the definition of a (generalized)
Atiyah-Patodi-Singer index class.

\subsubsection{{\bf Spectral sections.}}\label{subsubsect:spectral-sections}
 Let $D$ be a Dirac-type operator
acting between the sections of a Hermitian Clifford module $E$.
Let  $(M,r:M\to B\Gamma)$ be a Galois covering of a manifold with
boundary $M$. We shall concentrate on the even-dimensional case;
thus $\pa M$ is odd-dimensional. Let $\D_{(M,r)}$ be the $C^*_r
\Gamma$-linear operator associated to $D$ and $(M,r:M\to
B\Gamma)$. The starting point in Melrose-Piazza \cite{MP I} is the observation
that although the boundary operator $\D_{(\pa M,r_\pa)}$ is not
invertible, its index class in $K_1 (C^*_r \Gamma)$ is equal to
zero (by cobordism invariance). In order to define a higher
APS-index class in $K_0 (C^*_r \Gamma)$ we need a projection $\P$
playing the role of the non-existing projection $\Pi_{\geq}$. Of
course, we need somewhat special projections; these are nowadays
called spectral sections. Let $(N,s:N\to B\Gamma)$ be a
odd-dimensional closed Galois covering  (we shall eventually
choose $(N,s:N\to B\Gamma)=(\pa M,r_\pa: \pa M\to B\Gamma)$). A
{\it spectral section} associated to $\D\equiv \D_{(N,s)}$ is a
self-adjoint $C^*_r \Gamma$-linear projection $\P$ with the
additional property that there exists smooth functions
$$\chi_j:\RR\to [0,1]\;\;\text{with}\;\;\chi_j
(t)=0\;\text{for}\;t\ll 0\;,\;\; \chi_j (t)=1\;\text{for}\;t\gg
0\,,$$ $\chi_2\equiv 1$ on a neighborhood of  the support of
$\chi_1 $, and such that $$\Im \chi_1(\D_{(N,s)}) \subset \Im \P
\subset \Im \chi_2(\D_{(N,s)})\,.$$ Intuitively, $\P$ is equal to
$1$ on the large positive part of the spectrum and equal to $0$ on
the large negative part of the spectrum, precisely as $\Pi_{\geq}$
when the latter is defined. In fact, we have already encountered
spectral sections in this paper; see the Remark at the end of
subsection \ref{subsect:moreAPS}. The main result is then the
following

\begin{theorem}\label{theo:main-theorem-MP} (\cite{MP I} \cite{LPCUT})
A spectral section for $\D_{(N,s)}$ exists if and only if $\Ind
(\D_{(N,s)})=0$ in $K_1 (C^*_r \Gamma)$.
\end{theorem}

\subsubsection{{\bf Index classes and relative index theorem.}}

The cobordism invariance of the numeric index can be extended to
index classes Rosenberg \cite{Rosenberg} \cite[Proposition 2.3]{LPGAFA}.
Thus $\Ind (\D_{(\pa M,r_\pa)})=0$ in $K_1 (C^*_r \Gamma)$; hence
there exists a spectral section $\P$ for $\D_{(\pa M,r_\pa)}$. We
can use this $C^*_r \Gamma$-linear projection in order to define
the domain $$C^\infty (M,E^+\otimes\mathcal{V},\P)=\{s\in C^\infty
(M,E^+\otimes\mathcal{V}) \,|\, \P(s|_{\pa M} )=0 \}.$$ One can
prove that $\D_{(M,r)}$ with domain  $C^\infty
(M,E^+\otimes\mathcal{V},\P)$ gives rise to an index class $\Ind
(\D_{(M,r)},\P)\in K_0 (C^*_r \Gamma)$ \`a la
Atiyah-Patodi-Singer, see Wu \cite{Wu I} \cite{LPCUT} for the proof.
Different choices of spectral sections produces different index
classes; however there is a relative index theorem describing how
these index classes are related: given spectral sections $\P$,
$\Q$ there is a difference class $[\P - \Q]\in K_0 (C^*_r \Gamma)$
such that
\begin{equation}\label{relative-K}
\Ind (\D_{(M,r)},\Q) \,-\, \Ind
(\D_{(M,r)},\P)=[\P-\Q]\;\;\text{in}\;\; K_0 (C^*_r \Gamma)\,.
\end{equation}
This relative index theorem, first proved in \cite{MP I} and then
extended in \cite{LPAGAG} \cite{LP04}, is the higher analogue of
the last formula  of subsection
\ref{subsect:moreAPS}.

\subsubsection{{\bf Perturbations.}}
Let $(N,s:N\to B\Gamma)$ odd dimensional and $D$ a Dirac type
operator. Assume that $\Ind (\D_{(N,s)})=0$. Fix a spectral section
$\P$. Using  $\mathcal{P}$ one can construct a smoothing operator
$\mathcal{C}_{N,\mathcal{P}}\in \Psi^{-\infty}_{C^*_r \Gamma}$
 such that $\D_{(N,s)}+\mathcal{C}_{N,\mathcal{P}}$
is invertible in the Mishchenko-Fomenko calculus. Moreover
\begin{equation*}
\mathcal{P}=\frac{1}{2}\left(\Id+
\frac{\D_{(N,s)}+\mathcal{C}_{N,\mathcal{P}}}{|
\D_{(N,s)}+\mathcal{C}_{N,\mathcal{P}}|} \right).
\end{equation*}
In words: $\P$ is the positive spectral projection for the
perturbed operator $\D_{(N,s)}+\mathcal{C}_{N,\mathcal{P}}$. We
call such an operator $\mathcal{C}_{N,\mathcal{P}}$ a {\it
trivializing perturbation}. Let us go back to the case where
$(N,s)=(\partial M,r_\pa) $, with $(M,r:M\to B\Gamma)$ an even
dimensional Galois covering with boundary. Fix a spectral section
$\P$ for $\D_{(\pa M,r_\pa)}$; fix a trivializing perturbation
$\mathcal{C}_{\pa M, \mathcal{P}}$. One can extend the operator
$\mathcal{C}_{\pa M,\mathcal{P}}$ to the whole manifold with
boundary $M$. The resulting operator $\mathcal{C}_{M,\mathcal{P}}$
gives {\it  a perturbation  $\D_{(M,r)}+
\mathcal{C}_{M,\mathcal{P}} $ which has, by construction, an
invertible boundary operator.} It turns out that the index class
$\Ind (\D_{(M,r)},\P)$ \`a la Atiyah-Patodi-Singer can also be
described as an $L^2$-index class for the perturbed operator
$\D_{(M,r)}+ \mathcal{C}_{M,\mathcal{P}} $ extended to the
manifold obtained by adding a cylindrical end to $M$. Cylindrical
index theory is also referred to as $b$-index theory, because of
the exhaustive treatment given by Melrose using the
$b$-pesudodifferential calculus. See \cite{Melrose}. {\it
Summarizing}: the index class $\Ind (\D_{(M,r)},\P)$ is equal to
the $b$-index class $\Ind_b (\D_{(M,r)}+
\mathcal{C}_{M,\mathcal{P}})$. The advantage in considering the
latter index class comes from the invertibility of the boundary
operator: this allows us to consider the higher eta invariant of
the boundary operator, $\widetilde{\eta}( (\D_{(\pa M,r_\pa)}+
\mathcal{C}_{\pa M,\mathcal{P}}))$ and prove a higher index
formula similar to \ref{theo:memoirs}. The higher eta invariant,
denoted \begin{equation}\label{phigher} \widetilde{\eta}_{(\pa
M,r_\pa),\P} \,: = \, \widetilde{\eta}( (\D_{(\pa M,r_\pa)}+
\mathcal{C}_{\pa M,\mathcal{P}})) \end{equation} only depends on
$\P$ (and not on the particulat choice of perturbation) {\it
modulo exact forms}. This program is achieved in \cite{MP I}
\cite{MP II} in the family case and in \cite{LPGAFA} in the Galois
covering case. Recent topological applications of this general
theory are given in Piazza-Schick \cite{PS}.

\subsection{Middle-degree invertibility and a perturbation of the signature complex.}
\label{subsect:middle-degree}$\;$

\subsubsection{{\bf The middle-degree assumption.}}

 Let us now go back to the signature operator $\D^{{\rm
sign}}_{(M,r)}$ on a covering with boundary $(M,r:M\to B\Gamma)$
and to the problem of defining higher signatures when the operator
$\D^{{\rm sign}}_{(\pa M, r_\pa)}$ is {\it not} invertible. The
above subsection shows how to extend the Atiyah-Patodi-Singer
index theory developed in the invertible case to this general
case: crucial to this extension is the notion of 
{\it trivializing perturbation}. Unfortunately, the relative index
formula (\ref{relative-K}) shows very clearly that the resulting
index classes will depend on the choice of trivializing
perturbation. This is not very encouraging if our  goal is to
produce a {\it homotopy invariant APS index class.} In his
fundamental paper \cite{Lott II} Lott  points out an heuristic
cancellation mechanism indicating why the following assumption
might be sufficient for defining a {\it canonical} signature
class.

\m
Let $(N,s:N\to B\Gamma)$ be an odd dimensional Galois covering of
a closed oriented manifold. For example $(N,s:N\to B\Gamma)=(\pa
M,r_\pa: \pa M\to B\Gamma)$. Let $2m-1=\dim N$. Let $d$ denote the
de Rham differential on $\widetilde{N}$. Endow $\widetilde{N}$
with a $\Gamma-$invariant Riemannian metric.

\begin{assumption} \label{GreatLott}
The differential form Laplacian acting on $L^2(\widetilde{N},
\Lambda^{m-1}(\widetilde{N})) / \ker d$ has a strictly positive
spectrum.
\end{assumption}
 If $\mathcal{V}=C^*_r
\Gamma\times_\Gamma s^* E\Gamma$ and if $d_{\mathcal{V}}$ denotes
the twisted de Rham differential, then it is proved in \cite{LLP}
that  Assumption (\ref{GreatLott}) for $(N,s)$ is equivalent to
the following:

\m
\n
\begin{assumption}\label{GreatLott2}
Let $\Omega_{(2)}^{\ell}(N,{\cal V})$ denote the $L^2_{C^*_r
\Gamma}$ -completion of  $\Omega^{\ell}(N,{\cal V})$. The operator
$$d_{\mathcal{V}}:\Omega_{(2)}^{m-1}(N,{\cal V}) \rightarrow
\Omega_{(2)}^{m}(N,{\cal V}),$$ with domain equal to the $C^*_r
\Gamma$-Sobolev space $H^1_{C^*_r \Gamma}$, has {\it closed
image}.
\end{assumption}

\n
These equivalent assumptions  are for example satisfied  when $N$
has a cellular decomposition without any cell of dimension $m$.
Thanks to a deep result of  Gromov-Shubin \cite{gromov-shubin} we
know that these are {\it homotopy invariant conditions}. Notice
that if Assumption \ref{GreatLott} is satisfied, then necessarily
the index class of the signature operator in $K_1 (C^*_r \Gamma)$
is equal to zero.

\n
Since the index class of the signature operator is concentrated in
middle degree,  Assumption \ref{GreatLott} makes us guess that it
should be possible to find a set of {\it symmetric} trivializing
perturbations of the boundary operator producing first of all a
well defined higher eta invariant and, secondly,   a well defined
index class, both independent of the perturbation chosen. This is
indeed the case. There are in fact two equivalent ways to proceed:
one, proposed by Lott in \cite{Lott 3} and fully developed in
\cite{LLP} contructs perturbations of the {\it signature
complexes}, on $\pa M$ and on $M$, with the right symmetry
property for making the eta invariant and the index class well
defined. This is the approach we shall explain below. The other
approach, developed in Leichtnam-Piazza \cite{LPAGAG}, makes use of a special set
of spectral sections for the boundary signature operator; these
spectral sections have a certain symmetry property with respect to
forms of  degree $(m-1)$. We mention the approach through {\it
symmetric spectral sections} because we shall use it later, in
conjunction with the cut-and-paste problem for higher signatures.
 Thus, following John Lott 
\cite{Lott 3} and Leichtnam-Lott-Piazza \cite{LLP}
 we shall now explain how it is
possible to {\it add} a finitely $\Bi-$generated perturbation to
the {\it complex} of $\Bi-$differential forms on $\pa M$ and
consequently  perturb $\Si_{(\pa M , r_\pa)}$ into an {\it
invertible} (generalized) signature operator. To this aim we have
to recall, in the next sub-section, how to express $\Si_{(\pa M ,
r)}$ in terms of the $\Bi-$flat exterior derivative $d$ and the
Hodge duality operator $\tau$ acting on the Hermitian complex of
differential forms.

\subsubsection{{\bf More on the signature complex on closed
manifolds.}} First of all we recall the following

\begin{definition}  \label{Hcomplex}
A graded regular $n$-dimensional Hermitian complex consists of

1. A $\ZZ$-graded cochain complex $({\cal E}^*, D)$ of
finitely-generated projective left ${\cal B}^\infty$-modules,

2. A nondegenerate quadratic form $Q : {\cal E}^* \times {\cal
E}^{n-*} \rightarrow {\cal B}^\infty$ and

3. An operator $\tau \in Hom_{{\cal B}^\infty} \left( {\cal E}^*,
{\cal E}^{n-*} \right)$ such that
\begin{enumerate}
\item[1.] $Q(bx, y) = b Q(x,y)$.
\item[2.] $Q(x,y)^* = Q(y,x)$.
\item[3.] $Q(Dx, y) + Q(x, Dy) = 0$.
\item[4.] $\tau^2 = I$.
\item[5.] $<x, y> \equiv Q(x, \tau y)$ defines a Hermitian metric on ${\cal
E}$.
\end{enumerate}

\end{definition}

Let $M$ be a closed oriented $n$-dimensional Riemannian manifold
and let  $r:M\to B\Gamma$ be a reference map. We set ${\cal
V}^\infty=\Bi\times_\Gamma r^* E\Gamma$. Let $\Omega^*(M; {\cal
V}^\infty)$ denote the vector space of smooth differential forms
with coefficients in ${\cal V}^\infty$. The twisted de Rham
differential will be still denoted by $d$. If $n=\dim(M)
> 0$ then $\Omega^*(M; {\cal V}^\infty)$ is not finitely-generated
over ${\cal B}^\infty$, but we wish to show that it still has all
of the formal properties of a graded regular $n$-dimensional
Hermitian complex. If $\alpha \in \Omega^*(M; {\cal V}^\infty)$ is
homogeneous, denote its degree by $|\alpha|$. In what follows,
$\alpha$ and $\beta$ will sometimes implicitly denote homogeneous
elements of $\Omega^*(M; {\cal V}^\infty)$. Given $y \in M$ and
$(\lambda_1 \otimes e_1), (\lambda_2 \otimes e_2) \in
\Lambda^*(T^*_y M) \otimes {\cal V}^\infty_y$, we define
$(\lambda_1 \otimes e_1) \wedge (\lambda_2 \otimes e_2)^* \in
\Lambda^*(T^*_y M) \otimes {\cal B}^\infty$ by $$(\lambda_1
\otimes e_1) \wedge (\lambda_2 \otimes e_2)^* = (\lambda_1 \wedge
\overline{\lambda_2}) \otimes <e_1, e_2>.$$ Extending by linearity
(and antilinearity), given $\omega_1, \omega_2 \in \Lambda^*(T^*_y
M) \otimes {\cal V}^\infty_y$, we can define $\omega_1 \wedge
\omega_2^* \in \Lambda^*(T^*_yM) \otimes {\cal B}^\infty$. Define
a ${\cal B}^\infty$-valued quadratic form $Q$ on $\Omega^*(M;
{\cal V}^\infty)$ by $$ Q(\alpha, \beta) = i^{- |\alpha| (n -
|\alpha|)} \int_M \alpha(y) \wedge {\beta}(y)^*. $$ It satisfies
$Q(\beta, \alpha) = {Q(\alpha, \beta)}^*$. Using the Hodge duality
operator $*$, define $\tau : \Omega^p(M; {\cal V}^\infty)
\rightarrow \Omega^{n-p}(M; {\cal V}^\infty)$ by $ \tau (\alpha) =
i^{- |\alpha| (n-|\alpha|)} * \alpha. $ Then $\tau^2 = 1$ and the
inner product $<\cdot,\cdot>$ on $\Omega^*(M; {\cal V}^\infty)$ is
given by $< \alpha, \beta> = Q(\alpha, \tau \beta)$. Define $D :
\Omega^*(M; {\cal V}^\infty) \rightarrow \Omega^{*+1}(M; {\cal
V}^\infty)$ by
\begin{equation}
D \alpha = i^{|\alpha|} d \alpha. \label{signed-diff}
\end{equation}
{\bf Warning}: in this subsection the differential $D$ should not
be confused with a Dirac-type operator.

\medskip
\n
 It satisfies $D^2 = 0$. Its dual $D^\prime$
with respect to $Q$, i.e., the operator $D^\prime$ such that
$Q(\alpha, D\beta) = Q(D^\prime \alpha, \beta)$, is given by
$D^\prime = -D$. The formal adjoint of $D$ with respect to
$<\cdot, \cdot>$ is $D^* = \tau D^\prime \tau = - \tau D \tau$.

\begin{definition}  \label{sigop} {\rm
If $n$ is even, the signature operator is
\begin{equation}
\Si_{(M,r)} = D + D^* = D - \tau D \tau. \label{evensign}
\end{equation}
It is formally self-adjoint and anticommutes with the
$\ZZ_2$-grading operator $\tau$. If $n$ is odd, the signature
operator is
\begin{equation}
\Si_{(M,r)} = -i (D \tau + \tau D). \label{oddsign}
\end{equation}
It is formally self-adjoint. }
\end{definition}

\subsubsection{{\bf More on the signature complex on manifolds with
boundary.}}
 Now suppose that $M$ is a compact oriented
manifold-with-boundary of dimension $n = 2m$. Let $r: M\rightarrow
B \Gamma$ be a reference map and let $\partial M$ denote the
boundary of $M$. We fix  a Riemannian metric on $M$ which is
isometrically a product in an (open) collar neighbourhood
${\cal{U}}\equiv (0,2)_x\times\pa M$ of $\partial M$.  Let ${\cal
V}^\infty_0$ denote the pullback of ${\cal V}^\infty$ from $M$ to
$\partial M$; there is a natural isomorphism $${\cal
V}^\infty|_{{\cal{U}}}\cong (0,2) \times {\cal V}^\infty_0 \,.$$

\n One can show that,
 up to  explicit isomorphisms, the signature
operator can be written near the boundary as $ \D^{{\rm sign},
+}_{(M,r)}=\pa_x + \Si_{(\pa M , r)} $.

\subsubsection{{\bf The perturbed signature complex.}}
Recall that  ${\cal B}^\infty$ denotes the Connes-Moscovici
sub-algebra of
 $C^*_r(\Gamma)$ and that on  $\partial M$ we have the bundles: $$ \Fl=C^*_r(\Gamma) \times_\Gamma
 \partial \widetilde{M} , \quad\quad \Fli=\Bi\times_\Gamma\partial \widetilde{M}. $$
 The following Proposition, from \cite{Lott 3}, is proved using
Assumption \ref{GreatLott} in a crucial way. See   \cite{LLP}.

\begin{proposition} There exists  a cochain complex
 $C^* = \bigoplus_{k=-1}^{2m} C^k$ where
$C^k=\Omega^k(\partial M ;\Fli)\oplus \whW^k$ and  $\whW^*$ is  a
 complex made by finitely generated left projective
$\Bi-$modules. There are two maps  $\widehat{f} :
\Omega^*(\partial M ;\Fli)\rightarrow \widehat{W}^*$ and
$\widehat{g}: \widehat{W}^*\rightarrow \Omega^*(\partial M ;\Fli)$
such that the following property is satisfied. For any real
$\epsilon
>0$ the  differential $D_C$  on $C^*$ defined by
\begin{equation}
D_C=
\begin{pmatrix} D_{\partial M} & \epsilon \widehat{g} \cr 0
& - D_{\whW} \end{pmatrix} \;\;{\rm if }\;\,*< m- {1 \over
2},\,\;\;\; D_C=
 \begin{pmatrix} D_{\partial M} & 0 \cr - \epsilon \widehat{f} & -
D_{\whW}
\end{pmatrix}
\;\;{\rm if }\,\;*> m-{1 \over 2}.
\end{equation}
is such that  $D_C^2=0$ and
  the complex $(C^*, D_C)$ has
vanishing cohomology.
\end{proposition}
Define a duality operator $\tau_C$ on $C^*$ by
\begin{equation}\tau_C =
 \begin{pmatrix} \tau_{\partial M} & 0 \cr 0 & \tau_{\whW}\end{pmatrix}  .
\label{tauc}
\end{equation}
 The
signature operator associated to the perturbed complex $(C^*,D_C)$
is defined to be $${\cal D}_{C\,,\pa}^{\rm sign} (\epsilon)
=-i(\tau_C D_C+ D_C \tau_C)$$ If $\epsilon >0$, it follows from
the vanishing of the cohomology of $C^*$ that
\begin{equation}\label{perturbed-invertible}
{\cal D}_{C,\pa}^{\rm
sign}(\epsilon)\;\;\;\text{is an invertible self-adjoint}\;\;
\Bi\text{-operator}.\end{equation} We shall set ${\cal
D}_{C,\pa}^{\rm sign}\,:=\,{\cal D}_{C,\pa}^{\rm sign}( 1)$

\subsection{Lott's higher eta invariant in the non-invertible
case.}\label{subsubsect:lott-higher}$\;$

\m
We are now in a position to recall the definition of the higher
eta invariant of Lott  for a closed $(2m-1)$-dimensional covering
$(N,s:N\to B\Gamma)$ satisfying Assumption \ref{GreatLott}. This
material comes from \cite{Lott 3} and \cite{LLP}. We shall
concentrate directly on the case $(N,s:N\to B\Gamma)=(\partial M,
r_\partial) $.

Let $$\nabla : \Omega^*(\partial M ; {\cal V}^\infty) \rightarrow
\Omega_1(\Bi) \otimes_{\Bi} \Omega^*(\partial M ; {\cal
V}^\infty)$$ be Lott's connection for the bundle $E=\Lambda^* (\pa
M)$, see (\ref{lott-connection}). As in \cite[(3.28)]{Lott 3}, let
\begin{equation}
\nabla^{\widehat{W}^*}
 : \widehat{W}^* \rightarrow
\Omega_1(\Bi) \otimes_{\Bi} \widehat{W}^*
\label{sdconnect}\end{equation} be a connection on $\widehat{W}^*$
which is invariant under the grading operator and preserves the
quadratic form of $\widehat{W}^*$. Set $\nabla^C = \nabla \oplus
\nabla^{\widehat{W}^*}$; thus  $$\nabla^C : \Omega^*(\partial M ;
{\cal V}^\infty) \oplus \whW ^* \rightarrow
\widehat{\Omega}_1(\Bi) \otimes_{\Bi} \left( \Omega^*(\partial M ;
{\cal V}^\infty)\oplus \whW ^* \right) \,.$$ Let ${\rm Cl(1)}$ be
the complex Clifford algebra of $\CC$ generated by 1 and $\sigma$,
with $\sigma^2=1$. Let $\epsilon \in C^\infty(0, \infty)$ now be a
nondecreasing function such that
$\epsilon(s)=0$ for $s\in (0,1]$ and $\epsilon(s)=1$ for $s\in
[2,+\infty)$. Consider the element in  $\widehat{\Omega}_{\rm
even}(\Bi)/\overline{[\widehat{\Omega}_*(\Bi),\widehat{\Omega}(\Bi)]}$
\begin{equation}
\widetilde{\eta}_{(\partial M,r_\pa)}(s)=  \; {2 \over \sqrt{\pi}}
\,  \, \STR_{{\rm Cl}(1)} \left( ({d\over ds} [\sigma s{\cal
D}_{C,\pa}^{\rm sign}(\epsilon (s))] ) 
\exp[-(\sigma s{\cal D}_{C,\pa}^{\rm sign}(\epsilon (s))+
\nabla^C)^2] \right)\,;
\label{integrand}
\end{equation}
here $ {\rm STR}_{{\rm Cl}(1)}$ is defined as in subsection
\ref{subsect:bismut-cheeger}. The higher eta invariant of
$(\partial M, r_\pa)$ is, by definition,
\begin{equation}
\widetilde{\eta}_{\partial M} =\int_0^\infty
\widetilde{\eta}_{\partial M}(s) ds. \label{heta}\end{equation}
Since $\epsilon(s)=0$ for $s\in (0,1]$, it follows that the
integral is convergent for $s\downarrow 0$ (in fact, the integrand
near $s=0$ is the same as the one  for the unperturbed operator
and for the latter we know that convergence is implied by  Lott's
heat-kernel proof of the higher index theorem). Since $\epsilon
(s)=1$ for $s>2$ and {\it since the perturbed signature operator
${\cal D}_{C,\pa}^{\rm sign}$ is invertible}, it follows that the
integral is also convergent as $s\uparrow \infty$.
 It is shown in
  \cite[Proposition 14]{Lott 3} that, {\it modulo exact forms},
  the higher eta invariant  $\widetilde{\eta}_{(\partial M,r_\pa)}$ is independent of
the particular choices of the function $\epsilon$, the perturbing
complex $\widehat{W}^*$ and the self-dual connection
$\nabla^{\widehat{W}}$ .

\medskip

\subsection{Homotopy invariant higher signatures on a manifold with boundary.}$\;$

\subsubsection{{\bf Conic and cylindrical higher index
classes.}}\label{subsubsect:conic-and-cylindrical} Having defined
the higher eta invariant under the more general hypothesis of
middle-degree invertibility, we would like to show that it enters
as a {\it boundary correction term} in a higher index theorem for
a homotopy-invariant index class on our covering with boundary. We
begin \cite{LLP} by recalling the construction of a perturbed
signature operator ${\cal D}_{C,M}^{{\rm sign},{\rm cone}},$ with
boundary operator equal to the {\it invertible} perturbed
signature operator ${\cal D}_{C,\pa}^{{\rm sign}}$ introduced in
the previous subsection, see (\ref{perturbed-invertible}).

We take an (open) collar neighborhood of $\pa M$ which is
diffeomorphic to $(0,2) \times \pa M$. Let $\varphi \in
C^\infty(0,2)$ be a nondecreasing function such that $\varphi(x) =
x$ if $x \le 1/2$ and $\varphi(x) = 1$ if $x \ge 3/2$. Given $t >
0$, consider a Riemannian metric on ${\rm {int }} ( M ) $ whose
restriction to $(0,2) \times \partial M$ is
\begin{equation}g_M  =  t^{-2}  dx^2  + \varphi^2(x)  g_{\partial M}. \label{metric}
\end{equation}

Consider the complex $\Omega^*_c(0,2) \,
\widehat{\otimes} \, \whW^*$. It is endowed with a natural 
differential $D_{\rm alg}.$ Then set:
$$C^*=\Omega^*_c(M ; \Fl^\infty) \oplus \left( \Omega^*_c(0,2) \,
\widehat{\otimes} \, \whW^* \right),$$
 $C^* $ is endowed with  a natural direct sum duality operator $\tau_C.
$

Let $\phi \in C^\infty(0,2)$ be a nonincreasing function
satisfying $\phi(x) = 1$ for $0 < x \leq {1 \over 4}$ and $\phi(x)
= 0$ for ${1  \over 2} \leq x < {2 }$.
 We extend $\widehat{f}$ and $\widehat{g}$ to act on
$\Omega^*_c(0,2) \, \widehat{\otimes} \, \Omega^*(\partial M;
\Fli_0)$ and $\Omega^*_c(0,2) \, \widehat{\otimes} \,
\widehat{W}^*$, respectively, by $$\widehat{f} (\omega_0 + dx
\wedge \omega_1) \, = \, \widehat{f} (\omega_0) - \, i \,  dx
\wedge \widehat{f}(\omega_1)$$ and $$\widehat{g} (w_0 + dx \wedge
w_1) \, = \, \widehat{g} (w_0) - \, i \,  dx \wedge
\widehat{g}(w_1).$$ Using the cutoff function $\phi$, it makes
sense to define an operator on $C^*$ by
\begin{equation}D^{{\rm cone}}_C =
\begin{cases}
\left(
\begin{pmatrix}D_M & \phi \widehat{g}  \cr
0 & D_{\rm alg} \end{pmatrix} \right) &\text{ if $\;* \le m - 1$,}
\cr & \cr \left(
\begin{pmatrix}D_M & 0 \cr
0 & D_{\rm alg} \end{pmatrix} \right) &\text{ if $\;* = m$, }\cr &
\cr \left(
\begin{pmatrix}D_M & 0 \cr
- \phi \widehat{f} & D_{\rm alg} \end{pmatrix} \right) & \text{ if
$\;* \ge m + 1$} .\cr
\end{cases}
\label{diffdef}\end{equation} Note that $(D_C^{{\rm cone}})^2 \ne
0$, as $\phi$ is nonconstant. We have thus defined an ``almost''
differential $D_C^{\rm cone}$  on the conic complex
$$C^*=\Omega^*_c(M ; \Fli) \oplus \left( \Omega^*_c(0,2) \,
\widehat{\otimes} \, \whW^* \right) .$$ The perturbed conic
signature operator ${\cal D}_{C}^{{\rm sign},{\rm cone}} =
D_C^{{\rm cone}} + (D_C^{{\rm cone}})^*$ satisfies $${\cal
D}_{C}^{{\rm sign},{\rm cone}} = D_C^{\rm cone} - \tau D_C^{\rm
cone} \tau.$$ By construction, the boundary signature operator
associated to ${\cal D}_{C}^{{\rm sign},{\rm cone}} $ is precisely
${\cal D}_{C,\pa}^{{\rm sign}}$, the perturbed signature operator
constructed in the previous section.\\ {\it Summarizing:} we have
defined a {\it perturbed signature complex} on $(M,r:M\to
B\Gamma)$ with the property  that the associated signature {\it
operator} has an invertible boundary operator.

Using this fundamental fact one can prove that $ {\cal
D}_{C}^{{\rm sign},{\rm cone},+}$
defines an index class
\begin{equation}\label{conic-index-class}
{\rm {Ind}}\, {\cal D}_{C}^{{\rm sign},{\rm cone}, +} \in
K_0(\Bi)\,.
\end{equation}
The proof, see  \cite{LLP}, employs in a crucial way elliptic
analysis on conic manifolds, see \cite{Cheeger}, Br\"uning-Seeley
\cite{BS}. We
shall see in a moment that the conic index class is  homotopy
invariant. This is a fundamental step in our strategy for defining
homotopy-invariant higher signatures. The last step will consist
in proving an index theorem. However, to do so it turns out that
the cylindrical, or $b$, picture is more convenient. Thus we
sketch briefly the construction of a $b-$signature operator ${\cal
D}_C^{{\rm sign},b}$  in an extended version of Melrose
$b-$calculus; the boundary operator will be  once again  ${\cal
D}_{C,\pa}^{\rm sign}$.

Thus, we consider a $b$-metric $g$ which is product like near the
boundary: $$g = \frac{d x^2}{x^2} + g_{\pa M},$$ for $0 <  x \leq
\frac{1} {2}.$ Recall that a $b$-differential form is locally of
the form $a(x,y)  \frac{d x}{x}\wedge d y^{I}.$ The space of
$b$-differential forms is usually denoted by $\ub\Omega^*$.

\medskip

We    consider a new differential $D_C$ on the perturbed complex
$C^*=\ub\Omega^*(M;\Fli)\oplus (\ub\Omega^* [0,2)
\widehat{\otimes} \whW)$; on the degree $j$-subspace we put
\begin{equation}D_C\equiv\left(\begin{pmatrix} D_M & 0 \cr
0 & D_{{\rm alg}}\end{pmatrix}
\right)+\begin{cases}\left(\begin{pmatrix}0 & \widehat{g}_b \cr 0
& 0 \end{pmatrix} \right)&\text{ if $\;j<m$}\cr
\left(\begin{pmatrix}0 & 0\cr -\widehat{f}_b & 0 \end{pmatrix}
\right)&\text{ if $\;j> m$}\cr
\end{cases} \label{differ}\end{equation}
where $\widehat{g}_b$ and $\widehat{f}_b$ are $b-$operators
associated in a natural way to $\phi  \widehat{g}$ and $\phi
\widehat{f} $ respectively.

Let ${\cal D}_{C}^{{\rm sign},b} =D_C + (D_C)^*$
be the $b$-signature operator associated to the $b$-complex
$(C^*,D_C)$. Then  ${\cal D}_{C}^{{\rm sign},b} =D_C- \tau_C D_C
\tau_C$ is odd with respect to the $\ZZ_2$-grading defined by the
Hodge duality operator $\tau_C$ on $C^*$. Since the boundary
operator is equal to ${\cal D}_{C,\pa}^{{\rm sign}} $  and is
therefore invertible, one can prove that  the {\it perturbed}
$b$-signature operator ${\cal D}_C^{{\rm sign},b,+} $ is $C^*_r
\Gamma$-Fredholm, i.e. invertible modulo $C^*_r \Gamma$-compact
operators. Thus there a well defined index class ${\rm Ind}\,
{\cal D}_{C}^{{\rm sign},b,+} \in K_0(\Bi)$. To prove these
statements an extended version of Melrose's $b$-calculus must be
used, see \cite{LLP}.

\n
The following theorem is proved in \cite{LLP}
\begin{theorem} \label{equality} The following equality holds
in $K_0(\Bi)= K_0(C^*_r \Gamma )$: $${\rm Ind} \,{\cal
D}_{C}^{{\rm sign},{\rm cone},+}= {\rm Ind}\, {\cal D}_{C}^{{\rm
sign},b,+}\,. $$
\end{theorem}

\begin{proof} (Sketch)
There is also a perturbed signature operator $ {\cal D}_C^{{\rm
sign}}$ with respect to an ordinary product-like metric on $M$
(meaning, of  type $dx^2 + g_{\pa M}$ near the boundary). Since
the associated boundary operator is still ${\cal D}_{C,\pa}^{{\rm
sign}} $, hence invertible, we can define the projection
$$\Pi_{\geq}\,=\,\ha \left( \Id + \frac{{\cal D}_{C,\pa}^{{\rm
sign}} } {| {\cal D}_{C,\pa}^{{\rm sign}} |} \right)\,$$ and a
 higher index class ${\rm Ind} ({\cal D}_C^{{\rm
sign},+},\Pi_{\geq})$ \`a la Atiyah-Patodi-Singer. One proves that
the following two equalities hold in  $K_0(\Bi)= K_0(C^*_r \Gamma
)$: $${\rm Ind} \,{\cal D}_{C}^{{\rm sign},{\rm cone},+} ={\rm
Ind} ({\cal D}_C^{{\rm sign},+},\Pi_{\geq})\,,\quad {\rm Ind}
({\cal D}_C^{{\rm sign},+},\Pi_{\geq})={\rm Ind}\, {\cal
D}_{C}^{{\rm sign},b,+}. $$
\end{proof}

\subsubsection{{\bf Homotopy invariance of the index class.}} We
can finally state the first crucial result toward a definition of
homotopy invariant higher signatures:

\begin{theorem} \label{theo:hominv-higher} Let $(M,r:M\to B\Gamma)$
be such that $(\pa M,r_\pa)$ satisfy the middle-degree assumption
\ref{GreatLott}. The index class $ {\rm Ind} {\cal D}_{C}^{{\rm
sign},{\rm cone},+}\in K_0(\Bi)$ is a homotopy invariant of the
pair $( M,
\partial M)$ and the classifying map $r: M \rightarrow B \Gamma$.
Consequently, the $b$-index class ${\rm Ind} {\cal D}_C^{{\rm
sign},b,+} $ is also a homotopy invariant.
\end{theorem}

\begin{proof} (Sketch) One observes that the resolvent of ${\cal D}_{C}^{{\rm sign},{\rm cone}} $ is
$C^*_r(\Gamma)$-compact and that $( D_C^{\rm cone} )^2 $ is small,
provided that the real $t>0$ is small (i.e. the length of the cone
is large). Then one can extend fundamental  results of
Hilsum-Skandalis \cite{H-S} for $t>0$ small enough, proving the
homotopy invariance of the index class. (We recall that Hilsum and
Skandalis  have proved the homotopy invariance of the index class
for a signature operator with coefficients in an almost flat
bundle of $C^*-$algebras). The details are somewhat of
a technical nature and can be found in \cite{LLP}.
\end{proof}

\subsubsection{{\bf The index theorem and the higher signature class $\widehat{\sigma}
(M,r)\in \widehat{H}_* (\Bi)$.}}
 Now we can state the following
theorem, proved in \cite{LLP}.
\begin{theorem} \label{HIT}
Under Assumption \ref{GreatLott} the following formula holds:
$$\ch\, {\rm Ind} {\cal D}_C^{{\rm sign},b,+} = \left[ \int_M L(M)
\wedge \omega - \frac {1} {2} \widetilde{\eta}_{(\partial
M,r_\pa)} \right] \;\;\;\text{in}\;\;\;\widehat{H}_* (\Bi)$$ where
$\omega_{(M,r)}$ is, once again, the bi-form appearing in Lott's
heat-kernel proof of the higher index theorem and
$\widetilde{\eta}_{(\partial M,r_\pa)}$ is the higher eta
invariant for the perturbed signature operator ${\cal D}^{{\rm
sign}}_{C,\pa}$.
\end{theorem}

Thus, under the middle-degree assumption \ref{GreatLott} on the
boundary covering $(\pa M,r_\pa:\pa M \to B\Gamma)$ we are finally
in the position of extending the definition of higher signature
class given in subsection \ref{def:higher-class}
\begin{equation}
 \widehat{\sigma} (M,r)\,:=\, \left[ \int_M
L(M) \wedge \omega_{(M,r)} - {1\over 2}
\widetilde{\eta}_{(\partial M,r_\pa)}\right] \in \widehat{H}_*
(\Bi)
\end{equation}

Using \ref{theo:hominv-higher} and \ref{HIT} we can finally state
one of the main results of \cite{LLP}:
\begin{theorem}
The class 
$\widehat{\sigma}(M,r)$ in $ \widehat{H}_* (\Bi)$
is a homotopy
invariant of the pair $(M,\pa M)$ and the classifying map $r:M\to B\Gamma$.
\end{theorem}

\subsubsection{{\bf Homotopy invariant higher signatures in the non-invertible case.}}
We are approaching the end of our journey. Let $\Gamma$ be a group
with the extension property. For example, $\Gamma$ is  Gromov
hyperbolic or virtually nilpotent. Let $c\in H^\ell(\Gamma, \CC)$
be a group cycle and let $\tau_c\in ZC^*(\CC \Gamma)$ be the
associated cyclic cocycle. We can assume $\tau_c$ to be extendable
and we still let $\tau_c \in ZC^*(\Bi)$ the extended cocycle.
\begin{definition}
 The complex number  $$ {\rm sign}(M,r;[c])=<\widehat{\sigma}(M,r),[\tau_c]>\in \CC$$ is called
 the higher signature associated to  $(M, r)$ and $[c]$.
\end{definition}

The following theorem gives an answer to {\bf Question 3} in
section \ref{sect:3questions}:

\begin{theorem} (\cite{LLP})
Let $(M,r:M\to B\Gamma)$ be a Galois covering with boundary $(\pa
M,r_\pa:\pa M\to B\Gamma)$ satisfying the middle-degree assumption
\ref{GreatLott}. Let $\Gamma$ be a finitely generated group with
the extension property. The higher signatures $$
 {\rm sign}(M,r;[c])=<\widehat{\sigma} (M,r),[\tau_c]>\,,\quad\widehat{\sigma}(M,r)\,:=\, [ \int_M
L(M) \wedge \omega_{(M,r)} - {1\over 2}
\widetilde{\eta}_{(\partial M,r_\pa)}] $$ are homotopy invariants
for each $[c]\in H^* (\Gamma,\CC)$.
\end{theorem}

\subsection{Cut-and-paste invariance of higher signatures: the  index theoretic approach.}\label{subsect:cut-higher-index}$\;$

\m

We now go back to the cut-and-paste invariance of Novikov's higher
signatures on a {\it closed} manifold. We are looking for
sufficient conditions ensuring that the higher signatures are
indeed cut-and-paste invariant. Recall that for the lower
signature we explained 3 approaches  to the problem: \\(i) index
theoretic, \\(ii) topological, \\(iii) via a spectral-flow
argument. \\The following theorem, from \cite{LLP}, extends to the
higher case the first
 of these approaches. We shall only treat the even-dimensional
 case, the odd-dimensional case being more complicated to state
 and to treat.
\medskip

  Let $M$ and $N$ be two compact oriented $2m$-dimensional manifolds with boundary.
  Let $\phi$ and $\psi$ two orientation
preserving diffeomorphisms from $\partial M$ onto $\partial N$.
Consider the closed manifolds $X_\phi:=M\cup_\phi N^-$ and $X_\psi:=M\cup_\psi N^- \,.$ 
Let $r:
M\cup_\phi N^- \rightarrow B \Gamma$ and $s:
M\cup_\psi N^- \rightarrow B \Gamma$ be two reference maps; we assume that
these two coverings are {\it cut-and-paste equivalent}.

\begin{theorem} \label{addd} (\cite{LLP}) Assume that $\Gamma$ has the extension property
and that $(\partial M, r_\pa:\pa M \to B\Gamma)$ satisfies
Assumption \ref{GreatLott}. Then, for every $[c] \in H^*(\Gamma,
\CC)=H^*(B\Gamma, \CC)$ one has:
\begin{equation}\label{first-equation}
 <L(M\cup_\phi N^-) \cup r^*[c],[M\cup_\phi N^-]> = <
\widehat{\sigma}(M,r_{| M}),[\tau_c]>\,+\, <
\widehat{\sigma}(N^-,r_{|N^-}),[\tau_c]> \end{equation}
\begin{equation}\label{second-equation}
  <L(M\cup_\phi
N^-) \cup r^*[c],[M\cup_\phi N^-]> = <L(M\cup_\psi N^-) \cup
s^*[c], [M\cup_\psi N^-]> \,.
\end{equation}
In particular under the stated assumptions the higher signatures
are cut-and-paste invariant.
\end{theorem}

\noindent {\bf Remark.} Notice that $(\partial M, r_\pa:\pa M \to
B\Gamma)$ satisfies Assumption \ref{GreatLott} iff $(\partial M,
s_\pa:\pa M \to B\Gamma)$ satisfies it.

\begin{proof} We begin by (\ref{first-equation}). As in subsection
\ref{subsect:index-proof}
  we write:
$$ M\cup_\phi N^- = M \cup_{{\rm Id}} {\rm Cyl}_\phi \cup_{{\rm
Id}} N^- $$ where  ${\rm Cyl}_\phi = ([-1, 0]\times(\partial M)^-
) \cup_{{\phi}} ([0, 1]\times \partial N)$ is isomorphic to ${\rm
Cyl}:= [-1, 1]\times \partial M$ via $\phi$.
Moreover
\begin{equation}\label{zero} 
 \widehat{\sigma}({\rm Cyl}_\phi,r_{|{\rm Cyl}_\phi})= \int_{{\rm
Cyl}_\phi } L({\rm Cyl}_\phi ) \wedge \omega + \frac {1} {2}
\widetilde{\eta}_{(\partial M,r_{|\pa M})} - \frac {1} {2}
\widetilde{\eta}_{(\partial N,r_{|\pa N})}
=0.
\end{equation}
since by the established  homotopy invariance
$\widehat{\sigma}({\rm Cyl}_\phi,r_{|{\rm Cyl}_\phi})=
\widehat{\sigma}({\rm Cyl},r_{|\pa M}\times\Id)$ and the latter is
zero for the usual orientation argument concerning the eta
invariant. By Lott's higher index theorem on closed manifolds $$
<L(M\cup_\phi N^-) \cup r^*[c],[M\cup_\phi N^- ]>\, =\, < [\tau_c]
;\int_{M\cup_\phi N^- } L(M\cup_\phi N^- ) \wedge \omega>$$  We
can rewrite the left hand side of (\ref{first-equation}) as $$<
[\tau_c] ;\int_{ } L(M) \wedge \omega -\frac {1} {2}
\widetilde{\eta}_{(\partial M,r_{|\pa M})}> +
 < [\tau_c] ;\int_{{\rm Cyl}_\phi } L({\rm Cyl}_\phi ) \wedge
\omega + \frac {1} {2} \widetilde{\eta}_{(\partial M,r_{|\pa M})}
+ \frac {1} {2} \widetilde{\eta}_{(\partial N^-, r_{|\pa N^-})}>$$
$$+< [\tau_c ];\int_{N^- } L(N^-) \wedge \omega- \frac {1} {2}
\widetilde{\eta}_{(\partial N^-,r_{|\pa N^-} )}>\,.$$
From (\ref{zero}) we immediately obtain (\ref{first-equation}).
Moreover, (\ref{second-equation})  is an immediate consequence of
(\ref{first-equation}).
\end{proof}

\section{{\bf The topological approach to the cut-and-paste problem for higher
signatures.}}\label{sect:topological-cut} 

In this section we shall describe a topological approach to the
study of cut and paste properties of higher signatures. This
material comes from Leichtnam-L\"uck-Kreck \cite{LLK} and should be seen as the higher
analogue of what we  presented in subsection
\ref{subsect:topology-proof}.
 Namely, assuming that $(\partial M, r_{\partial M})$ satisfies
Assumption \ref{GreatLott}, we shall define a symmetric signature
$\sigma( M,r) \in K_0(C^*_r(\Gamma))$ which is {\it both} a higher
generalization of the lower topological signature of $(M,
\partial M)$ {\it and} a generalization of the Mishchenko symmetric signature
when the boundary is empty. The
 properties of $\sigma(M,r)$, namely {\it additivity
 } and {\it  homotopy invariance}, will allow us to extend Theorem \ref{addd} to
the discrete finitely presently groups $\Gamma$ 
satisfying the Strong Novikov Conjecture.

\subsection{The symmetric signature on manifolds with boundary.}\label{subsect:llk}$\;$

\m
We shall follow the notation in \cite{LLK}; in particular we
denote by $\overline{M}\to M$ a Galois covering with base $M$.

Let $n=2m$ be an even integer and $M$ be an oriented compact
$n$-dimensional manifold possibly with boundary. Let $(M,r:M\to
B\Gamma)$ a Galois covering. Let $\overline{\partial M} \to
\partial M$ and
 $\overline{M} \to M$ be the $\Gamma$-coverings
associated to the maps $r|_{\partial M}: \partial M \to B\Gamma$
and $r: M \to B\Gamma$. Following Lott \cite[Section 4.7]{Lott II} and
\cite[Assumption 1 and Lemma 2.3]{LLP}, we make the following
assumption about $(\partial M,r|_{\partial M})$. \\

\begin{assumption} \label{assumption}
Recall that  $n = 2m$ . Let $C_*(\overline{\partial M})$ be the
cellular $\ZZ \Gamma$-chain complex. Then we assume that the
$C^*_r ( \Gamma)$-chain complex $C_*(\overline{\partial M})
\otimes_{\ZZ \Gamma} C^*_r ( \Gamma)$ is $C^*_r ( \Gamma)$-chain
homotopy equivalent to a $C^*_r ( \Gamma)$-chain complex $D_*$
whose $m$-th differential $d_m: D_m \to D_{m-1}$ vanishes.
\end{assumption}

\medskip
\n
Lemma 2.3 in Leichtnam-Lott-Piazza \cite{LLP} shows that this assumption is equivalent to
Assumption \ref{GreatLott}. Notice that Assumption
\ref{assumption} is equivalent to the assertion that the $m$-th
Novikov-Shubin invariant of $\overline{\partial M}$ is $\infty^+$
in the sense of Lott-L\"uck \cite[Definition 1.8, 2.1 and 3.1]{Lott-Lueck
(1995)}.

Under Assumption \ref{assumption} we shall now assign   to $(M,r)$
an element
\begin{eqnarray}
& \sigma ({M},r)  \in K_0 (C^*_r ( \Gamma)),& \label{sigma(M,r)}
\end{eqnarray}
Fix a chain homotopy equivalence $u: C_*(\overline{\partial M})
\otimes_{\ZZ \Gamma} C^*_r ( \Gamma) \rightarrow D_*$ as in
Assumption \ref{assumption}. Define $\overline{D}_*$ as the
quotient chain complex of $D_*$ such that $\overline{D}_i=D_i$ if
$0\leq i \leq m-1$ and $\overline{D}_i=0$ for $i\geq m$. One then
gets a Poincare pair $j_*: D_*\rightarrow \overline{D}_*$ whose
boundary is $D_*$.
  By glueing \cite{Ranicki(1981)}
$j_*: D_*\rightarrow \overline{D}_*$ with the Poincare pair $$
i_*:  C_*(\overline{\partial M}) \otimes_{\ZZ \Gamma} C^*_r (
\Gamma) \rightarrow C_*(\overline{         M}) \otimes_{\ZZ
\Gamma} C^*_r ( \Gamma) $$ with the help of $u$ (along the
boundary $ C_*(\overline{\partial M})$) one gets a true Poincare
complex whose signature in $L^0 (\CC \Gamma)$ is denoted
$\sigma_{\CC\Gamma} ({M},r)$. Our symmetric signature $\sigma
(M,r)\in K_0 (C^*_r \Gamma)$ is the image of this class under the
composition $$L^0 (\CC\Gamma)\to L^0 (C^*_r \Gamma)
\leftrightarrow K_0 (C^*_r \Gamma)\,.$$

This  construction of the invariant $\sigma({M},r)$ by glueing
algebraic Poincar\'e bordisms is motivated by and extends the one
of Weinberger \cite{Weinberger(1997)} (see also \cite[Appendix
A]{Lott 3}) who uses the more restrictive assumption that
$C_*(\overline{\partial M}) \otimes_{\zz \Gamma} C^*_r ( \Gamma)$ is
$C^*_r ( \Gamma)$-chain homotopy equivalent to a $C^*_r (
\Gamma)$-chain complex $D_*$ with $D_m = 0$. In fact, when $D_m =
0$ the invariant $\sigma({M},r)$ coincides with the one of
Weinberger \cite{Weinberger(1997)}. The relationship to symmetric
signatures of manifolds-with-boundary, and to the necessity of
Assumption \ref{assumption}, was pointed out by Weinberger (see
\cite[Section 4.1]{Lott 3}).

We will call $\sigma({M},r)\in K_0 (C^*_r \Gamma) $ the $C^*_r
\Gamma$-valued symmetric signature of $(M,r)$.
When $\partial M$ is empty, this element $\sigma({M},r) $ agrees
with the (Mischenko) symmetric signature we defined in
\ref{subsect:symmetric-closed}. See also \cite[page
26]{Ranicki(1981)} on this point.

\subsection{Properties of the symmetric  signature.}$\;$

\m
The main properties of this invariant will be that it occurs in a
glueing formula, is a homotopy invariant and is related to higher
signatures. More precisely:

\begin{theorem} \label{the: main properties of sigma(M,r)}

$\;$

\begin{enumerate}

\item Glueing formula
\label{the: main properties of sigma(M,r): glueing formula}
\\[1mm]
Let $M$ and $N$ be two oriented compact $2m$-dimensional manifolds
with boundary and let $\phi: \partial M \to \partial N$ be an
orientation preserving diffeomorphism. Let $r: M\cup_{\phi} N^-
\to B\Gamma$ be a reference map. Suppose that $(\partial
M,r|_{\partial M})$ satisfies Assumption \ref{assumption}. Then
$$\sigma({M\cup_{\phi}N^-},r) ~ = ~ \sigma({M},r|_M) -
\sigma({N},r|_N)\quad\text{in}\quad K_0 (C^*_r \Gamma);$$

\item Cut-and-Paste invariance \label{the: main properties of sigma(M,r): additivity}
\\[1mm]
Let $M$ and $N$ be two oriented compact $2m$-dimensional manifolds
with boundary and let $\phi,\psi: \partial M \to \partial N$ be
orientation preserving diffeomorphisms. Let $$(r: M\cup_{\phi} N^-
\to B\Gamma)\quad\text{and}\quad (s: M\cup_{\psi} N^- \to
B\Gamma)\;\;\text{be cut-and-paste equivalent}\,.$$ 
Suppose that $(\partial M,r|_{\partial M})$ satisfies
Assumption \ref{assumption}. Then $$\sigma( {M\cup_{\phi}N^-},r) ~
= ~ \sigma( {M\cup_{\psi}N^-},s)\quad\text{in}\quad K_0 (C^*_r
\Gamma) \,;$$

\item Homotopy invariance
\label{the: main properties of sigma(M,r): homotopy invariance}
\\[1mm]
Let $M_0$ and $M_1$ be two oriented compact $2m$-dimensional
manifolds possibly with boundaries together with reference maps
$r_i: M_i \to B\Gamma$ for $i = 0,1$. Let $(f,\partial f):
(M_0,\partial M_0) \to (M_1,\partial M_1)$ be an orientation
preserving homotopy equivalence of pairs with $r_1 \circ f \simeq
r_0$. Suppose that $(\partial M_0,r_0|_{\partial M_0})$ satisfies
Assumption \ref{assumption}. Then $$\sigma( {M_0},r_0) = \sigma(
{M_1}, r_1).$$

\end{enumerate}
\end{theorem}

The crux of the proof is Ranicki \cite[Proposition 1.8.2
ii)]{Ranicki(1981)} and the underlying philosophical idea is the
following: if $M,N,$ and $D$ are compact oriented manifolds with
boundary such that $\partial M=\partial N=\partial D$ then $M\cup
D^- - N \cup D^-$ is cobordant to $M\cup N^-$.

\subsection{On the cut-and-paste invariance of higher signatures
on closed manifolds.}\label{subsect:cut-higher-topological} $\;$

\m

From Theorem \ref{the: main properties of sigma(M,r)} \eqref{the:
main properties of sigma(M,r): additivity}, we obtain the
following corollary which extends \cite[Corollary 0.4]{LLP}, i.e
Theorem \ref{addd} above, to more general groups $\Gamma$.

\begin{corollary} \label{cor: Additivity of higher signatures} Recall that $n=2m$.
Let $M$ and $N$ be two oriented compact $n$-dimensional manifolds
with boundary
 and let
$\phi,\psi: \partial M \to \partial N$ be orientation preserving
diffeomorphisms. 
Let $$(r: M\cup_{\phi} N^- \to B\Gamma)\quad\text{and}\quad (s:
M\cup_{\psi} N^- \to B\Gamma)\;\;\text{be cut-and-paste
equivalent}\,.$$ Assume that the $\Gamma$-covering
associated to $r|_{\partial M}: \partial M \to
B\Gamma$ satisfies Assumption \ref{assumption}. Suppose
furthermore that the assembly map $\mu : K_n(B\Gamma) \to
K_n(C^*_r(\Gamma ) )$ is rationally injective. Then for all $c \in
H^*(B\Gamma,\QQ)$
\begin{eqnarray} \label{add}
\sign (M\cup_{\phi}N^-,r; [c] ) = \sign (M\cup_{\psi}N^-,s; [c]).
\end{eqnarray}
In words, under the stated assumptions the higher signatures are
cut-and-paste invariant.
\end{corollary}

\begin{proof} Since $\mu_\RR$ is assumed to be injective
we know that the equality of the symmetric signatures implies the
equality of all the higher signatures, see Proposition
\ref{from-mish-to-higher}. From Theorem \ref{the: main properties
of sigma(M,r)} \eqref{the: main properties of sigma(M,r):
additivity}
 we get immediately the
result.
\end{proof}

\n
{\bf Remark.} We have already remarked that
for groups having the extension
property the map $\mu_\RR$ is injective.
Thus Corollary \ref{cor: Additivity of higher signatures} is
indeed a generalization of Theorem \ref{addd}.

\section{{\bf Higher spectral flow and cut-and-paste
invariance.}}\label{sect:higher-spectral-flow}

In the subsections \ref{subsect:cut-higher-index},
\ref{subsect:cut-higher-topological} we have extended to the
higher context the index theoretic and topological proof of the
cut-and-paste invariance of the lower signature. The goal of this
Section is to (briefly) present the higher analogue of the third
and last approach, the one employing the notion of  spectral flow.
 Our strategy is to show, {\it
analytically}, that under the same assumptions of Theorem
\ref{the: main properties of sigma(M,r)} (\ref{the: main properties
of sigma(M,r): additivity}) above, the signature index classes of
two cut-and-paste equivalent coverings $(r: M\cup_{\phi} N^- \to
B\Gamma)$ and $ (s: M\cup_{\psi} N^- \to B\Gamma)$ are equal in
$K_* (C^*_r \Gamma)$. By Proposition
\ref{prop:from-index-to-higher} this will reprove Corollary
\ref{cor: Additivity of higher signatures}.\\ We shall follow
\cite{LPCUT}. Notice that Michel Hilsum has also obtained these
results by using the Kasparov intersection product and  a somewhat
different approach to boundary value problems in the
noncommutative context. See \cite{Hilsum}.

\subsection{Higher spectral flow.}$\;$

\m
 First of all we need a definition for the {\it higher
spectral flow}. This was defined in the family-case by Dai and
Zhang, \cite{Dai-Zhang 2}, and extended to the noncommutative
context by F. Wu \cite{Wu} and Leichtnam-Piazza \cite{LPGAFA}
\cite{LPCUT}. Let $(N,s:N\to B\Gamma)$ be an odd dimensional
Galois covering and let $\D_{(N,s)}$ a generalized $C^*_r
\Gamma$-linear Dirac operator. We assume that $\Ind \D^{{\rm
sign}}_{(N,s)}=0$ in $K_1 (C^*_r \Gamma)$. This is the case, for
example, if $(N,s:N\to B\Gamma)=(\pa M,r_\pa : \pa M\to B\Gamma)$,
with $(M,r:M\to B\Gamma)$ a Galois covering with boundary.
According to Theorem \ref{theo:main-theorem-MP} there exists
spectral sections for $\D_{(N,r)}$. Recall that given two spectral
section $\Q$ and $\P$, the difference class $[\P-\Q]\in K_0 (C^*_r
\Gamma)$ is well defined.
\\
Assume now that we have a continuous one-parameter family of such
operators, 
parametrized by a continuous family of  inputs (metrics,
connections, etc...); we denote by $({\D_u})_{u \in [0,1]}$ such a
family.
 Recall that for any $C^*$-algebra $\Lambda$ there
exists an isomorphism $
 {\cal U} : K_1( C^0( [0,1];\, \CC) \otimes \Lambda\,) \simeq \,
 K_1( \Lambda)
 $
 which is implemented by the evaluation map $f(\cdot)\otimes
 \lambda\rightarrow f(0)\lambda$.
Using the above isomorphism  $\,{\cal U}$ for $\Lambda=C^*_r
\Gamma$, one gets that the index class associated to the $C^0
[0,1]\otimes C^*_r \Gamma$-linear operator $({\D_u})_{u \in
[0,1]}$ vanishes in $K_1(C^0([0,1])\otimes \Lambda)$. Thus
according to Theorem \ref{theo:main-theorem-MP}
 the family $({\D_u})_{u \in
[0,1]}$
 admits a  (total) spectral section $\P=({\P_u})_{u \in [0,1]}$.

\begin{definition} \label{hsfeven}
If $\Q_0$
 (resp. $\Q_1$) is a spectral section associated with ${\D}_0 $
 (resp. ${\D}_1 $) then the  noncommutative (or higher) spectral flow
 ${\rm sf}( ({\D}_u)_{u \in [0,1]}; \Q_0, \Q_1)$ from $({\D}_0 , \Q_0)$
 to $({\D}_1 , \Q_1)$ through $ ({\D}_u)_{u \in [0,1]}$ is the
 $K_0(C^*_r \Gamma)$-class:
 $$
 {\rm sf}
( ({\D}_u)_{u \in [0,1]}; \Q_0, \Q_1)= [\Q_1-\P_1]-[\Q_0-\P_0] \in
K_0( C^*_r \Gamma).
 $$
\end{definition}

\smallskip

\n
This definition does not depend on the particular choice of
 the total spectral section $\P=({\P})_{u \in [0,1]}$.

\n
Theorem 1.4 in Dai-Zhang \cite{Dai-Zhang 2} proves that if $\Gamma$ is
trivial and $\Q_0=\Pi_{\geq}(0)$, $\Q_0=\Pi_{\geq} (1)$, then the
above definition agrees with the usual one (net number of
eigenvalues changing sign).

If  the family is {\it periodic} (i.e. $ {\D}_1= {\D}_0$) and if
we take $\Q_1= \Q_0  $ then
 the spectral flow $ {\rm sf}
( ({\D}_u)_{u \in [0,1]}; \Q_0, \Q_0)$ does not depend on the
choice of $\Q_1= \Q_0 $ and defines a $K$-theory class which is
intrinsically associated to the given periodic family; we shall
denote this class by $ {\rm sf} ( ({\D}_u)_{u \in S^1})$.

More generally we can consider a periodic family of operators
$({\D}_u)$ as above but acting on the fibers of a fiber bundle
$P\longrightarrow S^1$ with fibers diffeomorphic to our manifold
$M$. Also in this case there is a well-defined noncommutative
spectral flow $ {\rm sf} ( ({\D}_u)_{u \in S^1})\in K_0( C^*_r
\Gamma)$. We shall encounter an example of this more general
situation in the coming subsections.

\subsection{The defect formula for cut-and-paste equivalent
coverings.} $\;$

\m

The higher spectral flow fits into a variational formula for  APS
index classes; this formula is the analogue of formula
(\ref{variational}) in subsection \ref{subsubsect:variational}.
Thus let $(\D_{(M,r)} (u))_{u\in [0,1]}$ be a 1-paramater family
of $C^*_r \Gamma$-linear operator on a covering with boundary. Let
$(\D_{(\pa M,r_\pa)} (u))_{u\in [0,1]}$ be the associated boundary
family. Fix a spectral section $\Q_0$ for $\D_{(\pa M,r_\pa)} (0)$
and a spectral section for $\D_{(\pa M,r_\pa)} (1)$. Then the APS
index classes $\Ind (\D_{(M,r)} (1), \Q_1)$ and $\Ind (\D_{(M,r)}
(0),\Q_0)$ are well defined in $K_0 (C^*_r \Gamma)$ and the
following formula holds:
\begin{equation}\label{higher-variational}
\Ind (\D_{(M,r)} (1), \Q_1) \,-\,\Ind (\D_{(M,r)} (0),\Q_0)= {\rm
sf}(  (\D_{(\pa M,r_\pa)} (u))_{u\in [0,1]}  ;\Q_1,\Q_0)\quad\text{in}\;\;\;K_0(C^*_r \Gamma)
\end{equation}
Next, the gluing   formula (\ref{additivity-simple}) given for the
numeric indeces in subsubsection \ref{subsubsect:additive} can be
extended to index classes. We state it directly for the signature
operator: if $$X=M\cup_F N^-\,,\quad\text{with}\quad F=\pa M=-\pa
N^-\,$$ and $r:X\to B\Gamma$ is a classifying map, then
\begin{equation}
\Ind (\D^{{\rm sign}}_{(X,r)})= \Ind( \D^{{\rm
sign}}_{(M,r_{|M})},\P)+ \Ind( \D^{{\rm
sign}}_{(N^-,r_{|N^-})},\Id-\P)\,,\;\;\;\text{in}\;\;\;K_0 (C^*_r
\Gamma)
\end{equation}
with $\P$ a spectral section for $\D^{{\rm sign}}_{(\pa M,r_{|\pa
M})}$. This formula can be extended to $X_\phi=M\cup_\phi N^-$
with $\phi:\pa M\to \pa N$ an oriented diffeomorphism. Using these
two formulae and proceeding  as in the numeric case one can prove
a {\it defect formula} for the difference $ \Ind (\D^{{\rm
sign}}_{(X_\phi,r)})\,-\,\Ind (\D^{{\rm
sign}}_{(X_\psi,s)})\;,\;\text{in}\;\;\;K_0 (C^*_r \Gamma)\;, $
associated to two cut-and-paste equivalent coverings $r:
X_\phi:=M\cup_\phi N^-\,\to\,B\Gamma\quad\text{and}\quad s:
X_\psi:=M\cup_\psi N^-\,\to\, B\Gamma$:

\begin{theorem}\label{theo:sf=0} There exists a periodic family of twisted signature
operators on $F=\pa M$, $\{{\cal D}_F(\theta)\}_{\theta\in S^1} $,
such that
\begin{equation}\label{c-p-6}
\Ind{\cal D}^{{\rm sign}}_{(X_\phi,r)} -\Ind{\cal D}^{{\rm
sign}}_{(X_\psi,s)}= {\rm sf}(\{{\cal D}_F(\theta)\}_{\theta\in
S^1}) \;\;\;\text{in}\;\;\;K_0 (C^*_r \Gamma)
\end{equation}
\end{theorem}
The family appearing on the right hand side of (\ref{c-p-6})
 is a $S^1$-family acting on the fibers of the mapping torus
$M(F,\phi^{-1}\circ \psi)\rightarrow S^1$.

\subsection{Vanishing higher spectral flow and the cut-and-paste
invariance.} $\;$

\m
The equality of the index class with the Mishchenko symmetric
signature, and the example given in section
\ref{sect:cut-example}, show together that the right hand side of
formula (\ref{c-p-6}) is in general different from zero. This is
in contrast with the numeric case. The following result is proved
by making use of the {\it symmetric spectral sections} we alluded
to in subsection \ref{subsect:middle-degree}.

\begin{theorem}
Let $M$ and $N$ be two oriented compact $2m$-dimensional manifolds
with boundary and let $\phi,\psi: \partial M \to \partial N$ be
orientation preserving diffeomorphisms. We let $F=\pa M$. Let
$$(r: M\cup_{\phi} N^- \to B\Gamma)\quad\text{and}\quad (s:
M\cup_{\psi} N^- \to B\Gamma)\;\;\text{be cut-and-paste equivalent
coverings}\,.$$ Suppose that $(\partial M,r|_{\partial M})$
satisfies Assumption \ref{GreatLott}. Then
\begin{equation}\label{c-p-7}
{\rm sf}(\{{\cal D}_F(\theta)\}_{\theta\in S^1})=0
\;\;\;\text{in}\;\;\;K_0 (C^*_r \Gamma)
\end{equation}
Consequently, by \ref{theo:sf=0}, the signature index classes of
$(r: M\cup_{\phi} N^- \to B\Gamma)$ and $(s: M\cup_{\psi} N^- \to
B\Gamma)$ coincide. Thus, by Proposition
\ref{prop:from-index-to-higher}, if the assembly map is rationally
injective then  for all $c \in H^*(B\Gamma,\CC)$
\begin{eqnarray}
\sign (M\cup_{\phi}N^-,r; [c] ) = \sign (M\cup_{\psi}N^-,s; [c]).
\end{eqnarray}
\end{theorem}

\section{{\bf Open problems.}}\label{sect:open-problems} 

\medskip
\n
{\bf I.} Let $(M,r)$ be an even dimensional oriented manifold
with boundary such that Assumption \ref{GreatLott} (or
\ref{assumption}) is satisfied. Then one observes that 
 the $C^*_r \Gamma$-valued symmetric signature
 class $\sigma (M ,r)$ constructed in \cite{LLK} (see
 Subsection \ref{subsect:llk})
 and the signature index class of \cite{LLP} $\Ind \cal{D}^
{{\rm sign},b,+}_C $ (see Subsection
\ref{subsubsect:conic-and-cylindrical}) have the same gluing and
homotopy invariance properties. Moreover, when $\partial
M=\emptyset$, these two classes coincide: see Theorem
\ref{theo:sign=ind} . Therefore it is natural to {\it conjecture}
that
\begin{equation}\sigma
 (M ,r)=\Ind \cal{D}^ {{\rm sign},b,+}_C \quad\text{in}\quad K_0 (C^*_r \Gamma) \,.
 \end{equation}

\medskip
\n
{\bf II.}  Let $(M, \mathcal{F})$ and
 $(N,\mathcal{F}^\prime)$           be two
foliated manifolds with boundary such that the leaves are
even-dimensional oriented and transverse to the boundary. Then $
\mathcal{F} $ has a product structure near $\partial M$. One
should try to formulate for $(\partial M,  \mathcal{F}_{|\partial
M} )$ an assumption analogous to \ref{GreatLott} and then define
for $(M, \mathcal{F})$ a signature index class which should be a
leafwise homotopy invariant (see Baum-Connes \cite{BaC} for the boundaryless
case). Now let $\phi$ and $\psi$ be two diffeomorphisms from
$\partial M$ to $\partial N$ sending a leaf of $
\mathcal{F}_{|\partial M}$ onto
   a leaf of
$ \mathcal{F}^\prime_{|\partial N}$ and preserving the
orientation. Then one gets two closed foliated manifolds
$(M\cup_\phi N^-, \mathcal{F}_\phi )$ and $(M\cup_\psi N^-,
\mathcal{F}_\psi )$. Let $q$ denote the common codimension of $
\mathcal{F}_\phi $ and $\mathcal{F}_\psi $ and consider the two
corresponding Haefliger classifying maps (see \cite[page 11]{BaC}):
$$ h_\phi: M\cup_\phi N^- \rightarrow B\Gamma_q,\; \; h_\psi:
M\cup_\psi N^- \rightarrow B\Gamma_q. $$ Then for each $\alpha \in
H^*(B\Gamma_q,\QQ)$ one should try to compare $$ \int_{M\cup_\phi
N^-} \,L(M\cup_\phi N^- ) \cup h_\phi^*(\alpha)\; {\rm and}\;
\int_{M\cup_\psi N^-} \,L(M\cup_\psi N^- ) \cup h_\psi^*(\alpha).
$$
\medskip

\n
{\bf Remark.} For
the particular case of foliated bundles see the recent paper \cite{LP04}.



\end{document}